\documentclass[moor,sglanonrev]{informs5}
\usepackage{eqndefns-left} 
\RequirePackage{tgtermes}
\RequirePackage{newtxtext}
\RequirePackage{newtxmath}
\RequirePackage{bm}
\RequirePackage{endnotes}

\OneAndAHalfSpacedXII 

\usepackage{algorithm}
\usepackage{algpseudocode}
\usepackage{tikz}
\usepackage{mathtools}
\usepackage[hidelinks]{hyperref}
\usepackage{enumitem}
\usepackage[dvipsnames]{xcolor}
\usepackage[normalem]{ulem}
\usepackage{booktabs}
\usepackage{comment}
\usepackage{graphicx}

\usepackage[sort&compress]{natbib}
 \bibpunct[, ]{[}{]}{,}{n}{}{,}%
 \def\bibfont{\small}%
\EquationsNumberedThrough    

\TheoremsNumberedThrough     
\ECRepeatTheorems  %

\MANUSCRIPTNO{MOOR-0001-2024.00}

\newcommand{\dd}[1]{\text{\ttfamily #1}}

\begin{document}


\RUNAUTHOR{Qiu and Davarnia}

\RUNTITLE{Intersection Cuts via Decision-Diagram-Free Sets}

\TITLE{Intersection Cuts for Mixed-Integer Nonlinear Programs via Decision-Diagram-Free Sets}

\ARTICLEAUTHORS{%
\AUTHOR{Junyuan Qiu, Danial Davarnia}
\AFF{Edwardson School of Industrial Engineering,
Purdue University,
West Lafayette, Indiana 47907, USA,
\EMAIL{qiu311@purdue.edu},
\EMAIL{ddavarn@purdue.edu}}
} 

\ABSTRACT{%
We introduce a novel decision-diagram-based framework for generating intersection cuts to strengthen polyhedral relaxations of bounded mixed-integer nonlinear programs (MINLPs) with general nonlinear constraints. Building on decision diagrams (DDs) that provide outer approximations of MINLP feasible regions, we propose the concept of DD-free sets, which are convex regions whose interiors do not intersect the sets represented by the DDs. By exploiting the underlying graphical structure, we develop efficient methods for constructing these regions without requiring restrictive assumptions on the functional form of nonlinear constraints or access to distance oracles. We investigate the maximality properties of DD-free sets and devise algorithms for obtaining maximal forbidden regions within several structured classes. The resulting intersection cuts exploit geometric information retained in DD representations, enabling further strengthening of relaxations beyond what can be achieved through conventional DD-based convexification. We also introduce interval DDs to facilitate the representation of general nonlinear constraints and establish convergence results for intersections of constraint-specific DDs, providing a theoretical foundation for their integration into global optimization frameworks. Computational experiments on challenging open instances from MINLPLib, encompassing a broad range of nonlinear and nonconvex structures, demonstrate that the proposed approach can substantially improve dual bounds and achieve greater optimality-gap closure than state-of-the-art global optimization solvers. These findings highlight the potential of exploiting graphical representations of infeasible regions to develop stronger cutting planes for general MINLPs.
}%

\FUNDING{This work was supported in part by the U.S. Air Force Office of Scientific Research grant FA9550261B176, and the National Science Foundation grant CMMI-2547883.}



\KEYWORDS{Decision diagrams; Mixed-integer nonlinear programs; $S$-free sets; Intersection cuts; Cutting planes; Global optimization;} 

\maketitle


\section{Introduction}
\label{sec:introduction}

Intersection cuts provide a classical geometric approach to strengthening linear programming (LP) relaxations of mixed-integer linear programs
\cite{balas1971intersection,conforti2011corner}.
These cuts are derived from convex forbidden regions, or $S$-free sets, whose interiors do not intersect the relevant feasible set $S$.
Intersection cuts have also been developed for mixed-integer nonlinear programs (MINLPs), including problems with polynomial
\cite{bienstock2020outer} and signomial constraints
\cite{xu2025cutting}.
For general MINLP constraints, however, a central challenge is to construct $S$-free sets in a computationally tractable form and to exploit them efficiently in cut generation.

Several existing approaches construct such forbidden regions by exploiting specific properties of nonlinear constraint functions, including factorability \cite{serrano2019intersection} and Lipschitz continuity with known constants \cite{schmidt2022cost,goss2026norm}, or by assuming access to geometric information about the feasible set, such as a distance oracle \cite{bienstock2020outer}. These requirements limit their applicability to more general MINLPs that do not satisfy such structural assumptions.

Recently, \cite{davarnia2026graphical} showed that general bounded MINLPs can be modeled using decision diagrams (DDs), providing a graphical representation of their feasible regions. By employing cut-generation methods, this representation can be used to construct polyhedral outer approximations of MINLPs. In this paper, we show that such DD representation offers an additional advantage beyond generating cutting planes for the feasible region: it also provides graph-based structural information about the infeasible region of the problem. This observation motivates the notion of a decision-diagram-free (DD-free) set, defined as a convex set whose interior does not intersect the feasible region represented by the DD. Since this DD contains the original feasible set, every DD-free set is also $S$-free and can therefore serve as a forbidden region for cut generation. For points in the relaxation that lie outside the DD feasible region, we exploit the graphical structure of the DD to efficiently construct DD-free regions and develop procedures for attaining set maximality within different classes of DD-free sets.



This DD-free construction also provides a mechanism for strengthening DD-based outer approximations beyond the convexification of DDs representing individual nonlinear constraints of the MINLP. Existing DD-based cut-generation methods derive inequalities that are valid for the convex hulls of the sets represented by the DDs \cite{tjandraatmadja2019target,davarnia2021outer,davarnia2026graphical}. For fixed diagrams, however, such cuts cannot separate points that lie in the intersection of these convex hulls, even when they are excluded by an individual DD. By combining the nonconvex geometric information retained by the DDs with information from the current polyhedral relaxation, the proposed intersection cuts can separate such points when they are extreme points of that relaxation. This approach avoids the potentially costly construction of a joint DD representing the intersection of the solution sets associated with multiple nonlinear constraints.

To demonstrate the effectiveness of the proposed methods for constructing DD-free sets and generating intersection cuts, we conduct computational experiments on several classes of MINLPs from MINLPLib that are classified as `open' problems, meaning that a gap remains between the dual bounds obtained by existing solvers and the best-known primal bounds. These instances span a broad range of nonlinear structures, including quadratic, polynomial, fractional, exponential, logarithmic, and other forms. Our results show that the dual bounds obtained using the proposed DD-based cutting planes can substantially outperform those produced by state-of-the-art global solvers, highlighting the potential of the proposed approach.

\subsection{Review of Relevant Work}
\label{subsec:literature-review}

Balas \cite{balas1971intersection} introduced intersection cuts for integer programming; see also \cite{conforti2011corner} for a survey of corner polyhedra and intersection cuts.
For MINLPs, intersection cuts have been developed using $S$-free sets constructed by exploiting nonlinear constraint structures or geometric information about feasible regions.
For factorable nonlinear constraints, Serrano 
\cite{serrano2019intersection} constructs concave underestimators that are tight at the point to be separated and uses the corresponding nonnegative superlevel sets as convex $S$-free regions for generating intersection cuts.
Related constructions include outer-product-free sets for polynomial optimization
\cite{bienstock2020outer}, quadratic-free sets for quadratic constraints
\cite{munoz2022maximal}, and signomial term-free sets for signomial programming
\cite{xu2025cutting}.
Beyond these specialized constructions,
Bienstock et al.\
\cite{bienstock2020outer} also develop a cutting-plane method for strengthening polyhedral outer approximations using $S$-free Euclidean balls, assuming access to a distance oracle for a closed set $S$.
In a different approach,
Schmidt et al.\
\cite{schmidt2022cost} construct nonconvex outer approximations by excluding open $\ell_\infty$-balls around points violating implicit nonlinear equality constraints, with radii determined by constraint residuals and known global Lipschitz constants.
G\"o{\ss} et al.\
\cite{goss2026norm} then generalize this idea to vector-valued nonlinear inequality constraints and a broader class of norms.
In our approach, valid lower bounds on nonlinear constraint functions over box domains form the basis for compiling DD relaxations, enabling efficient graph-based construction of DD-free regions in a general bounded MINLP setting.

For intersection cuts, enlarging an $S$-free set can strengthen the resulting cut, motivating the study of inclusion-wise maximal $S$-free sets
\cite{bienstock2020outer,conforti2015cut}.
For $S=\mathbb{Z}^n$, full-dimensional maximal $S$-free sets are polyhedra with an integer point in the relative interior of each facet
\cite{lovasz1989geometry,basu2010maximal}.
Maximality results have also been developed for outer-product-free sets in polynomial optimization
\cite{bienstock2020outer} and for quadratic-free sets
\cite{munoz2022maximal,munoz2025characterization,munoz2026characterization}.
In our setting, we first analyze maximality for general convex sets and then specialize the analysis to structured families inspired by norm-induced DD-free balls, for which tractable results can be obtained.
We also develop efficient procedures for enlarging these sets and achieve maximality.

More broadly, spatial branch-and-bound (B\&B) methods for globally solving MINLPs typically rely on convex relaxations and bound tightening; see
\cite{kronqvist202650} for a recent review.
Tawarmalani and Sahinidis
\cite{tawarmalani2005polyhedral} develop a polyhedral branch-and-cut approach that exploits convex subexpressions to generate relaxations and cutting planes.
Further developments in strengthening polyhedral relaxations for specific problem classes include learning-based partitioning of variable domains to improve approximations for quadratic problems \cite{kannan2026strong}, convex hull characterization for bilinear programs \cite{Davarnia2017,KhademniaDavarnia2025,GuDeyRichard2023}, and convexification of multilinear terms in binary polynomial optimization \cite{DelPiaKhajavirad2017,DelPiaKhajavirad2018}.
DDs provide a complementary graphical approach to developing relaxations and cutting planes.

In their standard use for discrete optimization, DDs represent sets of variable assignments using layered directed acyclic graphs, with each root-to-terminal path encoding an assignment; see
\cite{bergman2016decision,van2024introduction}.
Early work developed binary decision diagrams (BDDs) for representing and manipulating Boolean functions 
\cite{akers1978binary,bryant1986graph}.
BDDs have also been used to construct feasible solutions and primal bounds for binary optimization
\cite{bergman2014bdd}.
Multivalued decision diagrams (MDDs) extend BDDs from binary variables to discrete variables with general finite domains 
\cite{srinivasan1990algorithms}.
Andersen et al.\ 
\cite{andersen2007constraint} introduced relaxed MDD representations for constraint programming by limiting the diagram width to control its size.
Cir\'e and van Hoeve 
\cite{cire2013multivalued} developed limited-width MDD relaxations for sequencing problems.
Cappart et al.\
\cite{cappart2022improving} use deep reinforcement learning to learn variable orderings for improving the bounds obtained from relaxed and restricted DDs.
Serra and Hooker
\cite{serra2020compact} develop compact weighted DD representations of near-optimal integer programming solutions for postoptimality analysis.
For discrete optimization, relaxed DDs have subsequently been integrated into B\&B frameworks to provide dual bounds 
\cite{gillard2021improving,kowalczyk2024flow}.
Methods for parallelizing DD-based B\&B algorithms have also been investigated to improve the scalability of the solution process \cite{bergman2014parallel,kataru:da:ja:2026}.

For continuous nonlinear programs, Davarnia
\cite{davarnia2021strong} extends DD-based modeling to continuous domains through arc-reduced DD relaxations.
For mixed-integer optimization, Salemi and Davarnia
\cite{salemi2023structure} establish necessary and sufficient conditions for DD-representability of compact mixed-integer sets using rectangular decompositions and develop a specialized DD-based Benders decomposition for bounded MILPs.
This extension has enabled the application of DD-based solution methods to broader problem classes such as 
unit commitment 
\cite{salemi2023structure}, 
unsplittable network-flow problems 
\cite{salemi2023solving}, 
and sparse regression with nonconvex regularization 
\cite{davarnia:ki:2025}.

DDs have also been used for cut generation; 
see \cite{castro2022decision} for a survey.
Tjandraatmadja and van Hoeve 
\cite{tjandraatmadja2019target} exploit the connection between DDs and polar sets to generate target cuts for the convex hull of the set represented by a relaxed DD.
Davarnia and van Hoeve 
\cite{davarnia2021outer} develop a cut-generating linear program and a subgradient-type method for separating points from such convex hulls and incorporate the resulting cuts into outer-approximation methods for integer nonlinear programs.
More recently, Davarnia et al.\ 
\cite{davarnia2026graphical} extend DD-based convexification and cut generation to general bounded MINLPs within a spatial B\&B framework that guarantees convergence to global optimal solutions.
In this paper, we build on these results by exploiting the DD structure used to model MINLPs to construct DD-free regions, which are then used to derive intersection cuts.


\subsection{Contributions}
\label{subsec:contributions}

\paragraph{Contributions to the MINLP literature}
This work introduces a new framework for constructing $S$-free sets for bounded MINLPs with constraints of general functional form, thereby addressing a gap in the literature for problems that do not satisfy the structural assumptions or fall within the problem classes covered by existing frameworks for deriving intersection cuts. We establish maximality properties for the proposed sets and develop efficient algorithms for attaining maximality in different settings. 
We discuss how these DD-free sets can be used to derive intersection cuts that strengthen MINLP relaxations. We further demonstrate the effectiveness of the proposed framework through computational experiments on a collection of open benchmark instances from the MINLPLib. Our approach achieves greater gap closure than state-of-the-art global solvers, even when those solvers are configured to aggressively exploit their full suite of relaxation-strengthening techniques and tools.

\paragraph{Contributions to the DD literature}
This work advances the recent line of research on extending DD-based solution methods to MINLPs by introducing a new framework for deriving an important class of cutting planes, namely intersection cuts, directly from DD structures to strengthen problem relaxations. 
In particular, we show that the proposed intersection cuts can strengthen previously developed DD cuts, providing a systematic mechanism for extracting additional strength from a given DD representation. We also establish convergence results for the intersection of DDs associated with multiple constraints, which provides a theoretical foundation for combining constraint-specific DD representations within a unified global optimization framework.
This work further extends the graphical DD framework through the introduction of interval DDs, which compactly represent state values for arc reduction when modeling general MINLP constraints. This representation facilitates the translation of forbidden regions identified in the DD space.

\medskip

The remainder of this paper is organized as follows.
Section~\ref{sec:preliminaries} presents the notation, problem formulation, and basic DD concepts.
Section~\ref{sec:interval-DD} introduces interval DDs, describes their construction as relaxations of individual nonlinear constraints, and reviews their use in outer approximation.
Section~\ref{sec:DD-free-IC} introduces DD-free sets, establishes convergence results for the associated cut closures, and describes how these sets are used to derive intersection cuts.
Section~\ref{sec:maximal-DD-free} establishes general maximality results for DD-free sets, characterizes maximality within the families of DD-free boxes and $\ell_p$-balls, and develops procedures for enlarging these sets to strengthen the resulting cuts.
Section~\ref{sec:computational-results} presents computational experiments on benchmark MINLP instances to evaluate the performance of the proposed methods against state-of-the-art solvers.
Section~\ref{sec:conclusion} concludes the paper.

\section{Preliminaries}
\label{sec:preliminaries}

In this section, we introduce the notation used throughout the paper,
formulate the class of MINLPs under consideration, and review the basic
concepts of decision diagrams.

\subsection{Notation}

We use bold letters for tuples, vectors, and matrices.
For a vector $\boldsymbol{x}\in\mathbb{R}^n$, $x_i$ denotes its $i$-th component.
We generally use calligraphic letters for sets (e.g., $\mathcal{P}$) and uppercase Roman letters for index sets (e.g., $I$),  with occasional exceptions following standard notation. For an index set $I$, $|I|$ denotes its cardinality. Symbols associated with decision diagrams (DDs) are typeset in typewriter font (e.g., $\dd{D}$).

We let $\mathbb{N}\coloneqq\{1,2,\dots\}$ and $\mathbb{W}\coloneqq\{0,1,2,\dots\}$ denote the sets of positive and nonnegative integers, respectively.
For $k\in\mathbb{W}$, we write $[k]\coloneqq\{1,2,\dots,k\}$, with $[0]\coloneqq\emptyset$.
For $\delta\in\mathbb{R}$, let 
$(\delta)_+
\coloneqq
\max\{\delta,0\}$ denote its positive part.


For a set $\mathcal{S}\subseteq\mathbb{R}^n$,
$\operatorname{conv}(\mathcal{S})$,
$\operatorname{cl}(\mathcal{S})$,
$\operatorname{int}(\mathcal{S})$,
$\operatorname{relint}(\mathcal{S})$, and
$\operatorname{bd}(\mathcal{S})$
denote its convex hull, closure, interior, relative interior, and
boundary, respectively.
For a convex set $\mathcal{C}$, 
$\operatorname{ext}(\mathcal{C})$ denotes its set of extreme points. 
For a nonempty compact set $\mathcal{S}\subseteq\mathbb{R}$, we write $\mathcal{S}\downarrow \coloneqq \min\{x:x\in\mathcal{S}\}$ and $\mathcal{S}\uparrow\coloneqq\max\{x:x\in\mathcal{S}\}$. 
Unless stated otherwise, an interval is understood to be nonempty and compact.
A \textit{hyperrectangle} is a Cartesian product of intervals.

For $p\in[1,\infty]$, let $\lVert\cdot\rVert_p$ denote the standard $\ell_p$-norm on $\mathbb{R}^n$. Specifically, $\lVert\boldsymbol{x}\rVert_p\coloneqq \left(\sum_{i \in [n]}|x_i|^p\right)^{1/p}$ for $p \in [1, \infty)$, and $\lVert\boldsymbol{x}\rVert_\infty\coloneqq\max_{i\in[n]}|x_i|$.
For $\boldsymbol{x}\in\mathbb{R}^n$ and $r\ge 0$, let
\[
    \mathcal{B}_p(\boldsymbol{x},r)
    \coloneqq
    \{\boldsymbol{y}\in\mathbb{R}^n:
    \lVert\boldsymbol{y}-\boldsymbol{x}\rVert_p\le r\}
\]
denote the closed $\ell_p$-ball of radius $r$ centered at $\boldsymbol{x}$.
For $\boldsymbol{x}\in\mathbb{R}^n$ and a nonempty set $\mathcal{S}\subseteq\mathbb{R}^n$, we define
\[
    d_p(\boldsymbol{x},\mathcal{S})
    \coloneqq
    \inf_{\boldsymbol{s}\in\mathcal{S}}
    \lVert\boldsymbol{x}-\boldsymbol{s}\rVert_p.
\]
When $\mathcal{S}$ is closed, the infimum is attained.

\subsection{Problem Formulation} 
\label{subsec:problem}

Consider an MINLP of the form 
\begin{subequations}
\label{eq:problem-def}
\begin{align}
\zeta^\ast \coloneqq \max_{\boldsymbol{x}}\quad
    & \boldsymbol{c}^\top\boldsymbol{x}
    \label{eq:problem-def-obj}\\
\text{s.t.}\quad
    & g^k(\boldsymbol{x}) \le b_k,
    \qquad \forall k\in K,
    \label{eq:problem-def-NL-constr}\\
    & x_i \in \mathcal{D}_i,
    \qquad \forall i\in I\cup J,
    \label{eq:problem-domain}
\end{align}
\end{subequations}
Here, $I$ and $J$ are disjoint index sets of integer and continuous variables, respectively, and $I\cup J=[n]$ indexes all decision variables. For each $i\in [n]$, $\mathcal{D}_i$ denotes the domain of variable $x_i$, where 
$\mathcal{D}_i=[\mathcal{D}_i\downarrow,\mathcal{D}_i\uparrow]\cap\mathbb{Z}$ 
if $i\in I$, and
$\mathcal{D}_i=[\mathcal{D}_i\downarrow,\mathcal{D}_i\uparrow]$
if $i\in J$.
We write $\mathcal{D}\coloneqq\prod_{i\in [n]}\mathcal{D}_i$ for the domain of $\boldsymbol{x}$. 
Let $K$ be a nonempty finite index set for the nonlinear constraints. For each $k\in K$,  we assume that $g^k$ is well defined and bounded below over $\mathcal{D}$. We define the feasible region associated with the $k$-th nonlinear constraint by
\[
    \mathcal{G}^k
    \coloneqq
    \{\boldsymbol{x}\in\mathcal{D}: g^k(\boldsymbol{x})\le b_k\}.
\]

\subsection{Decision Diagram Overview}
\label{subsec:DD-overview}

A DD $\dd{D}=(\dd{U},\dd{A},\dd{l}(\cdot))$ over variables $x_1,\dots,x_n$ is a layered directed acyclic graph.
Its node set $\dd{U}$ is partitioned into layers $\dd{U}_1,\dots,\dd{U}_{n+1}$, and its arc set $\dd{A}$ is partitioned into layers $\dd{A}_1,\dots,\dd{A}_n$. Arc layer $\dd{A}_i$ corresponds to variable $x_i$ for each
$i\in[n]$.
Layer $\dd{U}_1$ contains a single \emph{root} node $\dd{r}$, and layer $\dd{U}_{n+1}$ contains a single \emph{terminal} node $\dd{t}$. For each $i\in[n]$, every arc $\dd{a}\in\dd{A}_i$ is directed from its tail node $\dd{tail}(\dd{a})\in\dd{U}_i$ to its head node $\dd{head}(\dd{a})\in\dd{U}_{i+1}$. 
Each node $\dd{u}\in\dd{U}$ is associated with a \emph{state}, denoted by $\dd{s}(\dd{u})$, and the mapping $\dd{l}:\dd{A} \to \mathbb{R}$ assigns a label to each arc $\dd{a}\in\dd{A}$. 
For a node $\dd{u}\in\dd{U}$, let $\dd{e}^+(\dd{u})$ and $\dd{e}^-(\dd{u})$ denote the sets of its outgoing and incoming arcs, respectively. 
The width of $\dd{D}$ is the maximum number of nodes in any node layer, i.e., $\dd{width}(\dd{D})\coloneqq \max_{j\in [n+1]}|\dd{U}_j|$.

An $\dd{r}$--$\dd{t}$ \emph{path} $\dd{p}$ of $\dd{D}$ is an arc sequence $\dd{p}=(\dd{a}_1,\dots,\dd{a}_n)$ with $\dd{a}_i\in \dd{A}_i$ for all $i\in[n]$ such that $\dd{tail}(\dd{a}_1)=\dd{r}$, $\dd{head}(\dd{a}_n)=\dd{t}$, and $\dd{head}(\dd{a}_i)=\dd{tail}(\dd{a}_{i+1})$ for all $i\in[n-1]$. 
When the label mapping $\dd{l}:\dd{A}\to\mathbb{R}$ is scalar-valued, the path $\dd{p}$ encodes the assignment $x_i=\dd{l}(\dd{a}_i)$ for all $i\in[n]$. In particular, bounded subsets of $\mathbb{Z}^n$ can be represented by DDs with integer-valued arc labels. 

Throughout, we take the node set of a DD to consist of its root, its terminal, and the endpoints of its arcs.
For any such DD containing an $\dd{r}$--$\dd{t}$ path, we have
\begin{equation}
\label{eq:DD-size-bounds}
    |\dd{U}|\le 2|\dd{A}|
    \qquad\text{and}\qquad
    n\le|\dd{A}|.
\end{equation}
Indeed, every node is incident to at least one arc, and 
every $\dd{r}$--$\dd{t}$ path contains exactly $n$ arcs.
Consequently, in the subsequent running-time analyses,
we express the complexity in terms of $|\dd{A}|$,
as both $|\dd{U}|$ and $n$ are
$\mathcal{O}(|\dd{A}|)$.

In Section~\ref{sec:interval-DD}, we extend this DD representation by allowing arcs to carry intervals of values. The resulting \textit{interval decision diagrams} (IDDs) will be used to model MINLP structures that contain both discrete and continuous variables.

\section{Interval Decision Diagrams}
\label{sec:interval-DD}

IDDs extend standard DDs by allowing arcs to carry interval-valued labels; see, e.g., \cite{strehl1998symbolic}.
Here, we use IDDs to represent relaxations of nonlinear feasible regions as finite unions of hyperrectangles. 

Consider an IDD $\dd{D}=(\dd{U},\dd{A},\dd{l}(\cdot))$ with $n$ arc layers, where arc layer $\dd{A}_i$ corresponds to variable $x_i$ for each $i\in[n]$. 
Each arc $\dd{a}\in\dd{A}$ carries an interval-valued label $\dd{l}(\dd{a})\subseteq\mathbb{R}$.
For a node $\dd{u} \in \dd{U}$, let $\dd{P}(\dd{u})$ denote the set of all $\dd{r}$--$\dd{u}$ paths in $\dd{D}$. In particular, $\dd{P}(\dd{t})$ denotes the set of all $\dd{r}$--$\dd{t}$ paths. 

Consider a node $\dd{u}\in\dd{U}_i$ for some $i\in\{2,\dots,n+1\}$. Each path $\dd{p}=(\dd{a}_1,\dots,\dd{a}_{i-1})\in\dd{P}(\dd{u})$ induces the hyperrectangle
\[
    \mathcal{S}_{\dd{p}}
    \coloneqq
    \prod_{j\in[i-1]}\dd{l}(\dd{a}_j)
    =
    \left\{
        \boldsymbol{x}\in\mathbb{R}^{i-1}:
        x_j\in\dd{l}(\dd{a}_j)
        \quad \forall j\in[i-1]
    \right\}.
\]
Thus, $\mathcal{S}_{\dd{p}}$ is the set of assignments to $x_1,\dots,x_{i-1}$ that are consistent with the interval labels encountered along $\dd{p}$.
For such a node $\dd{u}$, define
\[
    \mathcal{S}^{\cup}(\dd{u})
    \coloneqq
    \bigcup_{\dd{p}\in\dd{P}(\dd{u})}
    \mathcal{S}_{\dd{p}}
    \subseteq\mathbb{R}^{i-1}.
\]
When $\dd{u}=\dd{t}$, each path in $\dd{P}(\dd{t})$ specifies all $n$ variables.
Therefore, we can characterize the entire solution set represented by $\dd{D}$ as
\[
    \operatorname{Sol}(\dd{D})
    \coloneqq
    \mathcal{S}^{\cup}(\dd{t})
    =
    \bigcup_{\dd{p}\in\dd{P}(\dd{t})}
    \mathcal{S}_{\dd{p}}
    \subseteq\mathbb{R}^{n}.
\]
Since $\dd{P}(\dd{t})$ is finite, $\operatorname{Sol}(\dd{D})$ is a finite union of compact hyperrectangles and is therefore compact.
For later use, we define $\mathcal{S}^{\cup}(\dd{r})$ to consist solely of the empty assignment.

We say that an IDD $\dd{D}$ is an \emph{IDD relaxation} of a set
$\mathcal{P}\subseteq\mathbb{R}^n$ if $\mathcal{P}\subseteq\operatorname{Sol}(\dd{D})$. 
Equivalently, every point in $\mathcal{P}$ belongs to a hyperrectangle
induced by some $\dd{r}$--$\dd{t}$ path of $\dd{D}$. In particular, if
$\dd{D}$ is an IDD relaxation of $\mathcal{P}$, then
\[
    \operatorname{conv}(\mathcal{P})
    \subseteq
    \operatorname{conv}(\operatorname{Sol}(\dd{D})).
\]





The IDD representation is closely related to the arc-reduced DDs used in DD-based approaches for nonlinear optimization \cite{davarnia2021strong,davarnia2026graphical}. Indeed, replacing each interval-labeled arc by arcs carrying its two endpoints, followed by standard arc reduction, preserves $\operatorname{conv}(\operatorname{Sol}(\dd{D}))$. Our development instead retains the interval-valued labels because the path-induced hyperrectangles themselves are used to characterize DD-free regions.

In the remainder of this section, we construct an IDD relaxation of $\mathcal{G}^k$ associated with a single nonlinear constraint $g^k(\boldsymbol{x}) \le b_k$, for $k \in K$. Following existing DD-based outer-approximation frameworks, we construct a separate IDD for each nonlinear constraint \cite{davarnia2021outer,davarnia2026graphical}.
Whenever the analysis applies to an arbitrary constraint $k\in K$, we suppress the constraint index for notational simplicity whenever no ambiguity arises and write $g\coloneqq g^k$, $b\coloneqq b_k$, and $\mathcal{G}\coloneqq\mathcal{G}^k$.

Section~\ref{subsec:IDD-components} introduces three components used in the IDD construction. Section~\ref{subsec:construct-IDD} then combines these components in a top-down compilation algorithm that yields an IDD relaxation of $\mathcal{G}$.
Section~\ref{subsec:outer-approximation} reviews how cuts derived from IDDs can be incorporated into an outer-approximation framework.
Finally, Section~\ref{subsec:illustrative-example} illustrates a limitation of these cuts for fixed diagrams and motivates the intersection cuts proposed in this paper.

\subsection{IDD Components}
\label{subsec:IDD-components}

This subsection introduces three components of the IDD construction: a subdomain partition, a look-ahead node state, and a node-merging oracle.

\subsubsection{Subdomain Partition}
\label{subsubsec:subdomain-partition}

Consider an IDD $\dd{D}=(\dd{U},\dd{A},\dd{l}(\cdot))$.
For each $i\in[n]$, let $m_i\in\mathbb{N}$ and define $L_i\coloneqq[m_i]$.
Let $\{\mathcal{D}_i^j\}_{j\in L_i}$ be a collection of nonempty subdomains of $\mathcal{D}_i$.
If $i\in I$, each $\mathcal{D}_i^j$ is a set of consecutive integers, i.e., $\mathcal{D}_i^j=[\mathcal{D}_i^j\downarrow,\mathcal{D}_i^j\uparrow]\cap\mathbb{Z}$, where $\mathcal{D}_i^j\downarrow,\mathcal{D}_i^j\uparrow\in\mathbb{Z}$.
If $i\in J$, each $\mathcal{D}_i^j$ is an interval, i.e., $\mathcal{D}_i^j = [\mathcal{D}_i^j\downarrow,\mathcal{D}_i^j\uparrow]$.
The subdomains cover $\mathcal{D}_i$ and are ordered consecutively so that, for every $j<m_i$,
\[
    \mathcal{D}_i^{j+1}\downarrow
    =
    \begin{cases}
        \mathcal{D}_i^j\uparrow+1, & i\in I,\\
        \mathcal{D}_i^j\uparrow,   & i\in J.
    \end{cases}
\]
Thus, the subdomains are disjoint for integer variables, whereas adjacent
subdomains for continuous variables may share a common endpoint.

\subsubsection{Look-ahead Node State}
\label{subsubsec:look-ahead-state}

Consider an IDD $\dd{D}=(\dd{U},\dd{A},\dd{l}(\cdot))$. 
For a node $\dd{u}\in\dd{U}_i$ with $i\in\{2,\dots,n+1\}$, define $\mathcal{S}_{\mathcal D}^{\cup}(\dd{u})
\coloneqq
\mathcal{S}^{\cup}(\dd{u})
\cap
\prod_{j=1}^{i-1}\mathcal{D}_j$. 
For the root node, we set 
$\mathcal{S}_{\mathcal D}^{\cup}(\dd{r})
\coloneqq
\mathcal{S}^{\cup}(\dd{r})$. 
Thus, $\mathcal{S}_{\mathcal D}^{\cup}(\dd{u})$ consists of the partial assignments reaching $\dd{u}$ that lie in the variable domains.
Since $\mathcal{S}_{\mathcal D}^{\cup}(\dd{r})$ contains only the empty assignment, we identify
$\mathcal{S}_{\mathcal D}^{\cup}(\dd{r})
\times
\prod_{j=1}^{n}\mathcal{D}_j$
with $\mathcal{D}$.

\begin{definition}[Look-ahead state]
\label{def:look-ahead-state}
Let $\dd{u}\in\dd{U}_i$ for some $i\in[n]$. A \emph{look-ahead state} of $\dd{u}$ with respect to $g(\boldsymbol{x})$ is any scalar $\dd{s}(\dd{u})$ satisfying
\[
    \dd{s}(\dd{u})
    \;\le\;
    \inf\big\{ 
        g(\boldsymbol{x}) : \boldsymbol{x}\in 
        \mathcal{S}_{\mathcal{D}}^{\cup}(\dd{u})
        \times
        \prod_{j=i}^{n}\mathcal{D}_j
    \big\}.
\]
Equivalently, $\dd{s}(\dd{u})$ is a valid lower bound on $g$ over all completions of the partial assignments reaching $\dd{u}$ that are consistent with the variable domains.
\end{definition}




The techniques in \cite{davarnia2026graphical} can be adapted to efficiently compute valid look-ahead states for general nonlinear functions, including non-separable functions. In particular, coordinatewise bounds on the partial assignments reaching a node can be updated recursively and used to construct a bounding box over which a lower bound on $g$ is computed. Since the adaptation is straightforward, we omit the details and refer interested readers to Section~3.4 of \cite{davarnia2026graphical} for a detailed treatment.

\subsubsection{Merging Oracle}
\label{subsubsec:merging-oracle} 

Consider an IDD $\dd{D}=(\dd{U},\dd{A},\dd{l}(\cdot))$.
Let $\omega\in\mathbb{N}$ denote a prescribed width limit.
For each intermediate node layer $\dd{U}_j$, $j=2,\dots,n$, we apply the merging oracle $\dd{Merge}(\omega,\dd{U}_j)$, which selects and merges nodes in $\dd{U}_j$ so that the resulting layer contains at most $\omega$ nodes.
Specifically, the oracle partitions $\dd{U}_j$ into at most $\omega$ nonempty subsets and replaces the nodes in each subset $\dd{M}$ by a single node $\bar{\dd{v}}$. All incoming arcs of nodes in $\dd{M}$ are redirected to $\bar{\dd{v}}$, whose look-ahead state is defined as $\dd{s}(\bar{\dd{v}})\coloneqq\min\{\dd{s}(\dd{v}):\dd{v}\in\dd{M}\}$.
Since $\mathcal{S}_{\mathcal{D}}^{\cup}(\bar{\dd{v}})=\bigcup_{\dd{v}\in\dd{M}}\mathcal{S}_{\mathcal{D}}^{\cup}(\dd{v})$, all partial assignments reaching the original nodes that satisfy the
variable domains are retained, and the minimum of their look-ahead states remains a valid look-ahead state for $\bar{\dd{v}}$.

The particular choice of node groups determines the merging policy.
Two general-purpose merging oracles for MINLPs are given in \cite[Algorithms~4 and 5]{davarnia2026graphical}. Since these oracles can be directly incorporated into our IDD construction, we treat $\dd{Merge}(\omega,\cdot)$ as a black-box oracle.

\subsection{IDD Construction}
\label{subsec:construct-IDD} 

We now present the top-down compilation procedure used to construct an IDD relaxation of $\mathcal{G}$ in Algorithm~\ref{alg:IDD-construction}. The procedure follows the standard top-down DD compilation framework in \cite{davarnia2026graphical}. 
For each $i\in[n]$ and each arc $\dd{a}\in\dd{A}_i$, the interval label satisfies $\dd{l}(\dd{a})=\operatorname{conv}(\mathcal{D}_i^j)=[\mathcal{D}_i^j\downarrow,\mathcal{D}_i^j\uparrow]$ for some $j\in L_i$. At each intermediate node layer, nodes whose look-ahead states exceed $b$ are pruned in Line~\ref{alg-line:IDD-construction-prune}. Here, pruning a node means deleting the node together with all of its incident arcs. The merging oracle is then applied in Line~\ref{alg-line:IDD-construction-merge} to enforce the width limit $\omega$.

\begin{algorithm}
\caption{IDD construction for a general nonlinear constraint}
\label{alg:IDD-construction}
\begin{algorithmic}[1]
    \Require nonlinear constraint 
    $g(\boldsymbol{x})\le b$ over $\mathcal{D}$; 
    subdomain partitions 
    $\{\mathcal{D}_i^j\}_{j\in L_i}$ for $i\in[n]$; 
    width limit $\omega$; 
    merging oracle $\dd{Merge}(\omega,\cdot)$
    \Ensure an IDD $\dd{D}=(\dd{U},\dd{A},\dd{l}(\cdot))$
    \State \textbf{Initialize root:} create the root node $\dd{r}\in\dd{U}_1$ with look-ahead state $\dd{s}(\dd{r})$
    \State \textbf{Initialize terminal:} create the terminal node $\dd{t}\in\dd{U}_{n+1}$
    \For{$i\in[n-1]$}
        \For{each node $\dd{u}\in\dd{U}_i$ and each $j\in L_i$}
        \label{alg-line:IDD-construction-expansion-start}
            \State create a node $\dd{v}\in\dd{U}_{i+1}$
            \State add an arc from $\dd{u}$ to $\dd{v}$ with label $\operatorname{conv}(\mathcal{D}_i^j)$
            \State compute a look-ahead state $\dd{s}(\dd{v})$ satisfying
            \Statex \hspace{\algorithmicindent}
            $\displaystyle
            \dd{s}(\dd{v})
            \le
            \inf\left\{
                g(\boldsymbol{x}) :
                \boldsymbol{x}\in
                \mathcal{S}_{\mathcal D}^{\cup}(\dd{v})
                \times
                \prod_{h=i+1}^{n}\mathcal{D}_h
            \right\}$
        \EndFor
        \label{alg-line:IDD-construction-expansion-end}
        \State prune every node $\dd{v}\in\dd{U}_{i+1}$ with $\dd{s}(\dd{v})> b$
        \label{alg-line:IDD-construction-prune} 
        \State enforce the width limit $\omega$ at layer $\dd{U}_{i+1}$ by applying $\dd{Merge}(\omega,\dd{U}_{i+1})$
        \label{alg-line:IDD-construction-merge}
    \EndFor
    \For{each $\dd{u}\in\dd{U}_n$ and each $j\in L_n$}\label{alg-line:IDD-construction-last-start}
        \State compute $\eta^\ast$ such that $
        \eta^\ast 
        \le 
        \inf\left\{
            g(\boldsymbol{x}): 
            \boldsymbol{x}
            \in 
            \mathcal{S}_{\mathcal{D}}^{\cup}(\dd{u})
            \times 
            \mathcal{D}_n^{j}
        \right\}$
        \label{alg-line:IDD-construction-last-state-calculation} 
        \If{$\eta^\ast \le b$}
        \label{alg-line:IDD-construction-last-condition}
            \State add an arc from $\dd{u}$ to $\dd{t}$ with label $\operatorname{conv}(\mathcal{D}_n^j)$
        \EndIf
    \EndFor
    \label{alg-line:IDD-construction-last-end}
    \State\Return $\dd{D}$
\end{algorithmic}
\end{algorithm}

We next establish the correctness of Algorithm~\ref{alg:IDD-construction} and characterize the size of the resulting IDD.
In particular, Theorem~\ref{thm:IDD-valid} shows that the resulting diagram is an IDD
relaxation of $\mathcal{G}$, i.e., $\mathcal{G}\subseteq\operatorname{Sol}(\dd{D})$, and hence $\operatorname{conv}(\mathcal{G})\subseteq\operatorname{conv}(\operatorname{Sol}(\dd{D}))$.

\begin{theorem}
\label{thm:IDD-valid}
Consider $\mathcal{G}\coloneqq\{\boldsymbol{x}\in \mathcal{D}: g(\boldsymbol{x})\le b\}$, and let $\dd{D}=(\dd{U},\dd{A},\dd{l}(\cdot))$ be the IDD constructed by Algorithm~\ref{alg:IDD-construction} for $\mathcal{G}$ using subdomain
partitions $\{\mathcal{D}_i^j\}_{j\in L_i}$, for $i\in[n]$, width limit $\omega$, and merging oracle $\dd{Merge}(\omega, \cdot)$. Then $\dd{D}$ is an IDD relaxation of $\mathcal{G}$, i.e., $\mathcal{G}\subseteq\operatorname{Sol}(\dd{D})$, and hence $\operatorname{conv}(\mathcal{G})\subseteq\operatorname{conv}(\operatorname{Sol}(\dd{D}))$.
\end{theorem}

\proof{Proof.}

Fix an arbitrary $\bar{\boldsymbol{x}}\in\mathcal{G}$. 
Since the subdomains cover the variable domains, for each $i\in[n]$, choose $j_i\in L_i$ such that $\bar{x}_i\in\mathcal{D}_i^{j_i}$.

If $n=1$, processing $(\dd{r},j_1)$ in the final loop yields $\eta^\ast\le \inf\{g(x_1):x_1\in\mathcal{D}_1^{j_1}\}\le g(\bar{\boldsymbol{x}})\le b$. Hence, the algorithm adds an arc $\dd{a}_1$ from $\dd{r}$ to $\dd{t}$ labeled $\operatorname{conv}(\mathcal{D}_1^{j_1})$. 
Since 
$\bar{x}_1
\in
\mathcal{D}_1^{j_1}
\subseteq
\operatorname{conv}(\mathcal{D}_1^{j_1})=\mathcal{S}_{\dd{p}}$, 
where $\dd{p}=(\dd{a}_1)$, we have 
$\bar{\boldsymbol{x}}
\in
\mathcal{S}_{\dd{p}}
\subseteq
\operatorname{Sol}(\dd{D})$. 
Thus, $\bar{\boldsymbol{x}}\in\operatorname{Sol}(\dd{D})$ when $n=1$. Henceforth, suppose that $n\ge 2$.

We prove by induction that, for every $i\in[n-1]$, after completing iteration $i$, including the pruning and merging steps, there exist a node $\dd{v}^{(i)}\in\dd{U}_{i+1}$ and an $\dd{r}$--$\dd{v}^{(i)}$ path $\dd{p}^{(i)}=(\dd{a}_1,\ldots,\dd{a}_i)$ such that $\dd{l}(\dd{a}_k)=\operatorname{conv}(\mathcal{D}_k^{j_k})$ for all $k\in[i]$ and
$(\bar{x}_1,\ldots,\bar{x}_i)
\in
\mathcal{S}_{\dd{p}^{(i)}}
\subseteq
\mathcal{S}^{\cup}(\dd{v}^{(i)})$.

\emph{Base case ($i=1$).}
When the algorithm processes $(\dd{r},j_1)$, it creates a node $\tilde{\dd{v}}^{(1)}\in\dd{U}_2$ and an arc $\dd{a}_1$ from $\dd{r}$ to $\tilde{\dd{v}}^{(1)}$ labeled $\operatorname{conv}(\mathcal{D}_1^{j_1})$. 
Let $\dd{p}^{(1)}=(\dd{a}_1)$. Since 
$\bar{x}_1
\in
\mathcal{D}_1^{j_1}
\subseteq
\operatorname{conv}(\mathcal{D}_1^{j_1})$, 
we have 
$\bar{x}_1
\in
\mathcal{S}_{\dd{p}^{(1)}}
\subseteq
\mathcal{S}^{\cup}(\tilde{\dd{v}}^{(1)})$. 
Moreover, since $\bar{\boldsymbol{x}}\in\mathcal{D}$, we have $\bar{x}_1\in\mathcal{S}_{\mathcal D}^{\cup}(\tilde{\dd{v}}^{(1)})$ and hence 
$\bar{\boldsymbol{x}}
\in
\mathcal{S}_{\mathcal D}^{\cup}(\tilde{\dd{v}}^{(1)})
\times
\prod_{h=2}^{n}\mathcal{D}_h$. 
Definition~\ref{def:look-ahead-state} therefore yields $\dd{s}(\tilde{\dd{v}}^{(1)})\le g(\bar{\boldsymbol{x}})\le b$, so $\tilde{\dd{v}}^{(1)}$ is not pruned. 
Let $\dd{v}^{(1)}$ denote the head of $\dd{a}_1$ after the merging step. 
Since merging only redirects incoming arcs without changing their labels, $\dd{p}^{(1)}$ remains an $\dd{r}$--$\dd{v}^{(1)}$ path satisfying the claimed properties.

\emph{Inductive step.}
If $n\ge 3$, fix $i\in\{2,\dots,n-1\}$ and suppose that the claim holds for $i-1$.
When the algorithm processes $(\dd{v}^{(i-1)},j_i)$, it creates a node $\tilde{\dd{v}}^{(i)}\in\dd{U}_{i+1}$ and an arc $\dd{a}_i$ from $\dd{v}^{(i-1)}$ to $\tilde{\dd{v}}^{(i)}$ labeled $\operatorname{conv}(\mathcal{D}_i^{j_i})$. 
The path $\dd{p}^{(i)}=(\dd{a}_1,\ldots,\dd{a}_{i-1},\dd{a}_i)$ therefore has the required labels. 
Since 
$\bar{x}_k
\in
\mathcal{D}_k^{j_k}
\subseteq
\operatorname{conv}(\mathcal{D}_k^{j_k})$ for all $k\in[i]$, we have 
$(\bar{x}_1,\ldots,\bar{x}_i)
\in
\mathcal{S}_{\dd{p}^{(i)}}
\subseteq
\mathcal{S}^{\cup}(\tilde{\dd{v}}^{(i)})$. 
Because $\bar{\boldsymbol{x}}\in\mathcal{D}$, it follows that 
$(\bar{x}_1,\ldots,\bar{x}_i)
\in
\mathcal{S}_{\mathcal D}^{\cup}(\tilde{\dd{v}}^{(i)})$, 
and hence 
$\bar{\boldsymbol{x}}
\in
\mathcal{S}_{\mathcal D}^{\cup}(\tilde{\dd{v}}^{(i)})
\times
\prod_{h=i+1}^{n}\mathcal{D}_h$. 
Thus, $\dd{s}(\tilde{\dd{v}}^{(i)})\le g(\bar{\boldsymbol{x}})\le b$, so $\tilde{\dd{v}}^{(i)}$ is not pruned. 
Let $\dd{v}^{(i)}$ denote the head of $\dd{a}_i$ after the merging step. 
Again, because merging only redirects incoming arcs, $\dd{p}^{(i)}$ remains an $\dd{r}$--$\dd{v}^{(i)}$ path satisfying the claimed properties. This completes the induction.

Applying the induction result with $i=n-1$, there exist a node $\dd{v}^{(n-1)}\in\dd{U}_n$ and an $\dd{r}$--$\dd{v}^{(n-1)}$ path $\dd{p}^{(n-1)}=(\dd{a}_1,\ldots,\dd{a}_{n-1})$ such that $\dd{l}(\dd{a}_k)=\operatorname{conv}(\mathcal{D}_k^{j_k})$ for all $k\in[n-1]$ and 
$(\bar{x}_1,\ldots,\bar{x}_{n-1})
\in
\mathcal{S}_{\dd{p}^{(n-1)}}
\subseteq
\mathcal{S}^{\cup}(\dd{v}^{(n-1)})$. 
Since $\bar{\boldsymbol{x}}\in\mathcal{D}$, we further have 
$(\bar{x}_1,\ldots,\bar{x}_{n-1})
\in
\mathcal{S}_{\mathcal D}^{\cup}(\dd{v}^{(n-1)})$.

When the final loop processes $(\dd{v}^{(n-1)},j_n)$, the inclusion $\bar{x}_n\in\mathcal{D}_n^{j_n}$ implies that 
$\bar{\boldsymbol{x}}
\in
\mathcal{S}_{\mathcal D}^{\cup}(\dd{v}^{(n-1)})\times\mathcal{D}_n^{j_n}$. 
Therefore, the computed value satisfies $\eta^\ast\le g(\bar{\boldsymbol{x}})\le b$, and the algorithm adds an arc $\dd{a}_n$ from $\dd{v}^{(n-1)}$ to $\dd{t}$ labeled $\operatorname{conv}(\mathcal{D}_n^{j_n})$.

Consequently, $\dd{p}=(\dd{a}_1,\ldots,\dd{a}_n)$ is an $\dd{r}$--$\dd{t}$ path satisfying $\dd{l}(\dd{a}_i)=\operatorname{conv}(\mathcal{D}_i^{j_i})$ for all $i\in[n]$. 
Since $\bar{x}_i\in\mathcal{D}_i^{j_i}\subseteq\dd{l}(\dd{a}_i)$ for every $i\in[n]$, we have 
$\bar{\boldsymbol{x}}
\in
\mathcal{S}_{\dd{p}}
\subseteq
\mathcal{S}^{\cup}(\dd{t})=\operatorname{Sol}(\dd{D})$. 
Since $\bar{\boldsymbol{x}}\in\mathcal{G}$ was arbitrary, together with the case $n=1$, we conclude that $\mathcal{G}\subseteq\operatorname{Sol}(\dd{D})$. 
Taking convex hulls for both sides yields $\operatorname{conv}(\mathcal{G})\subseteq\operatorname{conv}(\operatorname{Sol}(\dd{D}))$.
\Halmos
\endproof

The following result bounds the number of arcs in the IDD constructed by Algorithm~\ref{alg:IDD-construction}, showing that it is quadratic in the prescribed width limit and the subdomain partition size.

\begin{lemma}
\label{lemma:num-nodes-arcs-bounds}
Let $\dd{D}=(\dd{U},\dd{A},\dd{l}(\cdot))$ be the IDD constructed by
Algorithm~\ref{alg:IDD-construction} using width limit $\omega$ and
subdomain index sets $L_i$, $i\in[n]$. 
Then 
$|\dd{A}|
\le 
|L_1|+\omega\sum_{i=2}^n|L_i|$.
\end{lemma}

\proof{Proof.}

By construction, pruning can only delete arcs, whereas merging only redirects existing arcs.
For each $i\in[n]$, the algorithm adds at most one arc in
$\dd{A}_i$ for each pair
$(\dd{u},j)\in\dd{U}_i\times L_i$.
Hence, $|\dd{A}_i|\le|\dd{U}_i||L_i|$.
Since $|\dd{U}_1|=1$ and 
$|\dd{U}_i|\le\omega$ for every
$i\in\{2,\dots,n\}$, summing over the arc layers gives
the stated bound.
\Halmos
\endproof

By \eqref{eq:DD-size-bounds}, both $|\dd{U}|$ and $n$ are $\mathcal{O}(|\dd{A}|)$ for any IDD with a nonempty solution set.
In the remainder of the manuscript, we express the running-time complexity of the algorithms in terms of the number of arcs in the IDD, as Lemma~\ref{lemma:num-nodes-arcs-bounds} provides a corresponding bound in terms of the prescribed parameters.

\subsection{IDD Cuts}
\label{subsec:outer-approximation}

In the literature, DDs have been used within an outer-approximation framework to globally solve MINLPs \cite{davarnia2021outer, davarnia2021strong, davarnia2026graphical}. Inspired by this approach, we can leverage the correspondence between DDs and IDDs to design an outer-approximation framework for globally solving \eqref{eq:problem-def}. Below, we provide a high-level overview of the main components of this framework and refer the reader to the aforementioned references for full details.

Once an IDD relaxation $\dd{D}^k$ is constructed for each set $\mathcal{G}^k$ for $k \in K$, we can use cut-generation methods developed in \cite{davarnia2021outer} to efficiently generate valid inequalities for $\operatorname{conv}(\operatorname{Sol}(\dd{D}^k))$. By analogy with the \textit{DD cuts} introduced in the DD literature, we refer to these inequalities as \textit{IDD cuts}. These cuts can be incorporated into an outer-approximation framework to produce a linear relaxation of the feasible region of \eqref{eq:problem-def}. The resulting linear relaxation is then solved repeatedly, with additional IDD cuts generated from individual IDDs to progressively strengthen the relaxation (dual) bound.

Once the improvement in the dual bound becomes negligible, a spatial B\&B procedure is invoked to tighten variable domains through spatial partitioning. The process is then repeated by constructing stronger IDDs over the reduced variable domains. As shown in \cite{davarnia2026graphical}, this combination of outer-approximation and spatial B\&B guarantees convergence to the global optimal solution of the MINLP \eqref{eq:problem-def}.

Our goal in this paper, however, is to strengthen the aforementioned outer-approximation methodology by introducing a new tool for improving the dual bounds of the linear relaxation: intersection cuts derived from characterizing \textit{infeasible} regions encoded by IDDs. This tool complements the previously developed IDD cuts that form the outer-approximation. Specifically, these cuts are generated from feasible regions represented in the IDD structure, whereas the proposed intersection cuts are produced from infeasible regions induced by the IDD structure.

\subsection{Illustrative Example}
\label{subsec:illustrative-example}

\begin{figure}[tb]
    \centering
    \includegraphics[width=0.65\linewidth]{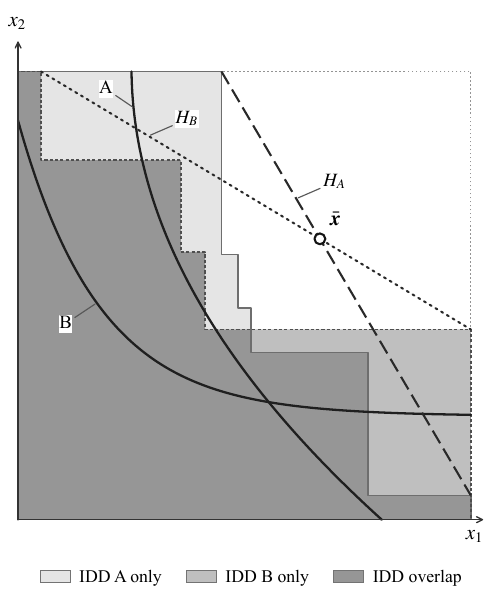}
    \caption{Limitation of IDD cuts for fixed diagrams.
        The thick solid curves labeled $A$ and $B$ depict the nonlinear boundaries of $\mathcal{G}^A$ and $\mathcal{G}^B$, respectively, with each feasible region lying on the side containing the origin.
        The light and medium gray regions belong only to
        $\operatorname{Sol}(\dd{D}^A)$ and
        $\operatorname{Sol}(\dd{D}^B)$, respectively, whereas the dark gray region is their intersection.
        The dashed line $H_A$ and the dotted line $H_B$ each define a facet of the corresponding convex hull.
        The marked point $\bar{\boldsymbol{x}}$ belongs to both convex hulls but neither IDD solution set.
        }
    \label{fig:idd-cuts-limitation}
\end{figure}

We illustrate the limitation of IDD cuts for fixed diagrams using a two-dimensional example, which will also be used in subsequent sections.
Consider two nonlinear constraints, labeled $A$ and $B$, in continuous variables $x_1$ and $x_2$ over a bounded rectangular domain $\mathcal{D}$.
Let $\dd{D}^A$ and $\dd{D}^B$ be IDD relaxations of the corresponding constraint sets $\mathcal{G}^A$ and $\mathcal{G}^B$, respectively.
Figure~\ref{fig:idd-cuts-limitation} depicts the nonlinear boundaries of these constraint sets and the solution sets represented by the IDDs.

For these fixed diagrams, adding all IDD cuts to the box relaxation yields the polytope
\[
    \mathcal{P}_{\mathrm{IDD}}
    \coloneqq
    \operatorname{conv}\bigl(\operatorname{Sol}(\dd{D}^A)\bigr)
    \cap
    \operatorname{conv}\bigl(\operatorname{Sol}(\dd{D}^B)\bigr).
\]
Together with the remaining axis-aligned supporting halfspaces, the halfspaces bounded by $H_A$ and $H_B$ exactly describe
$\operatorname{conv}(\operatorname{Sol}(\dd{D}^A))$ and
$\operatorname{conv}(\operatorname{Sol}(\dd{D}^B))$, respectively.
Since $\bar{\boldsymbol{x}}
\in
\mathcal{P}_{\mathrm{IDD}}$
lies at the intersection of the nonparallel lines
$H_A$ and $H_B$, it is an extreme point of
$\mathcal{P}_{\mathrm{IDD}}$.

Every IDD cut generated from either diagram is valid for its associated convex hull and therefore cannot separate $\bar{\boldsymbol{x}}$.
Nevertheless, $\bar{\boldsymbol{x}}$ lies outside both
$\operatorname{Sol}(\dd{D}^A)$ and
$\operatorname{Sol}(\dd{D}^B)$.
Since the corresponding constraint sets are contained in
their respective IDD solution sets,
$\bar{\boldsymbol{x}}$ violates both nonlinear constraints.
Thus, the fixed IDDs already exclude $\bar{\boldsymbol{x}}$,
but convexifying their solution sets separately does not retain
this exclusion.

In Section~\ref{sec:DD-free-IC},
we return to this example and combine a DD-free set from either IDD with the corner relaxation at $\bar{\boldsymbol{x}}$ to derive an intersection cut that separates this point.


\section{DD-Free Sets and Intersection Cuts}
\label{sec:DD-free-IC}

In this section, we develop a framework for deriving intersection cuts from an IDD relaxation. 
For each $k\in K$, let $\dd{D}^k=(\dd{U}^k,\dd{A}^k,\dd{l}^k(\cdot))$ denote the IDD relaxation of $\mathcal{G}^k$,
such as one constructed by Algorithm~\ref{alg:IDD-construction}.
Throughout the remainder of the paper, we assume that $\operatorname{Sol}(\dd{D}^k)\neq\emptyset$ for all $k\in K$. Indeed, if $\operatorname{Sol}(\dd{D}^k)=\emptyset$ for some $k\in K$, then $\mathcal{G}^k=\emptyset$, and hence \eqref{eq:problem-def} is infeasible.
Following the earlier convention, when considering a single constraint for an arbitrary $k\in K$, we suppress the index $k$ and write $\mathcal{G}\coloneqq\mathcal{G}^k$ and $\dd{D}\coloneqq\dd{D}^k$. We retain the index when referring to the family
$\{\dd{D}^k\}_{k\in K}$.

\subsection{$S$-Free and DD-Free Sets}
\label{subsec:DD-free-set}

Let $S\coloneqq \operatorname{cl}(\mathcal{G})$. We adopt the standard notion of an $S$-free set, which is the basic geometric object used to generate intersection cuts; see, e.g., \cite{balas1971intersection, conforti2015cut, bienstock2020outer}. 

\begin{definition}[$S$-free set]
\label{def:S-free-set}
    A convex set $\mathcal{C}\subseteq\mathbb{R}^n$ is \emph{$S$-free} if $\operatorname{int}(\mathcal{C})\cap S=\emptyset$.
\end{definition}

Replacing $S$ in Definition~\ref{def:S-free-set} by $\operatorname{Sol}(\dd{D})$ yields the following notion.

\begin{definition}[Decision-diagram-free set]
\label{def:DD-free-sets}
    Given an IDD $\dd{D}$, a convex set $\mathcal{C}\subseteq \mathbb{R}^{n}$ is \emph{DD-free with respect to $\dd{D}$} if $\operatorname{int}(\mathcal{C})\cap \operatorname{Sol}(\dd{D})=\emptyset$.
\end{definition}

In the remainder of this paper, whenever a set $\mathcal{C}$ is referred to as DD-free, it is understood to be convex unless stated otherwise.
The next proposition shows that DD-free sets are $S$-free and hence are valid forbidden regions for intersection-cut constructions for $\mathcal{G}$.

\begin{proposition}
\label{prop:DD-free-also-S-free}
    Let $\dd{D}$ be an IDD relaxation of $\mathcal{G}$.
    Then every DD-free set is $S$-free.
\end{proposition}

\proof{Proof.}

Let $\mathcal{C}$ be a DD-free set.
Because $\dd{D}$ is an IDD relaxation of
$\mathcal{G}$, we have $\mathcal{G}\subseteq\operatorname{Sol}(\dd{D})$. 
Since $\mathcal{C}$ is DD-free, $\operatorname{int}(\mathcal{C})\cap\mathcal{G}=\emptyset$. Hence $\mathcal{G}\subseteq\mathbb{R}^n\setminus\operatorname{int}(\mathcal{C})$.
Because the latter set is closed, $S=\operatorname{cl}(\mathcal{G})\subseteq\mathbb{R}^n\setminus\operatorname{int}(\mathcal{C})$.
Thus, $\operatorname{int}(\mathcal{C})\cap S=\emptyset$, and hence $\mathcal{C}$ is $S$-free by Definition~\ref{def:S-free-set}.
\Halmos
\endproof

To abstract the pointwise construction of DD-free sets from IDDs, we introduce the following oracle model.

\begin{definition}[Decision-diagram-free-set oracle]
\label{def:DD-free-oracle}
    Given an IDD $\dd{D}$, a set-valued mapping $\mathcal{C}(\cdot):\mathbb{R}^n \rightrightarrows \mathbb{R}^n$ is a \emph{DD-free-set oracle} for $\dd{D}$ if it satisfies the following properties:
    \begin{enumerate}[label=\textup{(\roman*)}, ref=\textup{(\roman*)}]
        \item $\mathcal{C}(\boldsymbol{x})=\emptyset$ for all $\boldsymbol{x}\in\operatorname{Sol}(\dd{D})$;
        \item for every $\boldsymbol{x}\notin \operatorname{Sol}(\dd{D})$, the set $\mathcal{C}(\boldsymbol{x})$ is DD-free with respect to $\dd{D}$ and $\boldsymbol{x}\in \operatorname{int}(\mathcal{C}(\boldsymbol{x}))$.
    \end{enumerate}
\end{definition}

An important component in designing a DD-free-set oracle is verifying that a convex set $\mathcal{C}(\boldsymbol{x})$ is DD-free, as required by condition~(ii) of Definition~\ref{def:DD-free-oracle}. The following result establishes the computational difficulty of this verification task by showing that the associated decision problem is coNP-complete even when $\mathcal{C}(\boldsymbol{x})$ is restricted to be a rational polyhedron.

\begin{proposition}
\label{prop:DD-free-intersection-hard}
    Given an IDD $\dd{D}$ with rational interval endpoints and a full-dimensional polyhedron $\mathcal{C}$ described by rational linear inequalities, deciding whether $\mathcal{C}$ is DD-free with respect to $\dd{D}$ is coNP-complete.
    Equivalently, deciding whether
    $\operatorname{int}(\mathcal{C})
    \cap
    \operatorname{Sol}(\dd{D})
    \neq
    \emptyset$
    is NP-complete.
    Both conclusions hold even when $\dd{D}$ has width one and every arc label is either $[0,0]$ or $[1,1]$.
\end{proposition}

\proof{Proof.}

It suffices to prove that deciding whether
$\operatorname{int}(\mathcal{C})
\cap
\operatorname{Sol}(\dd{D})
\neq
\emptyset$
is NP-complete,
since the two problems are complements.

We first show that this problem belongs to NP.
Write
$\mathcal{C}
=
\{\boldsymbol{x}\in\mathbb{R}^n:
\boldsymbol{A}\boldsymbol{x}
\le\boldsymbol{\beta}\}$,
where $\boldsymbol{A}$ has no zero rows.
Since $\mathcal{C}$ is full-dimensional, its interior is
$\{\boldsymbol{x}\in\mathbb{R}^n:
\boldsymbol{A}\boldsymbol{x}
<
\boldsymbol{\beta}\}$.
If
$\operatorname{int}(\mathcal{C})
\cap
\operatorname{Sol}(\dd{D})
\neq
\emptyset$,
then there exists an $\dd{r}$--$\dd{t}$ path $\dd{p}$ such that 
$\mathcal{S}_{\dd{p}}
\cap
\operatorname{int}(\mathcal{C})
\neq
\emptyset$, which can be considered as the \textit{certificate} for this assumption.
For any given path $\dd{p}\in\dd{P}(\dd{t})$, 
the condition 
$\mathcal{S}_{\dd{p}}
\cap
\operatorname{int}(\mathcal{C})
\neq
\emptyset$
holds if and only if the linear program
\[
    \max_{\boldsymbol{x},\varepsilon}
    \left\{
        \varepsilon:
        \boldsymbol{A}\boldsymbol{x}
        +\varepsilon\boldsymbol{1}
        \le
        \boldsymbol{\beta},\
        \boldsymbol{x}\in\mathcal{S}_{\dd{p}},\
        0\le\varepsilon\le 1
    \right\}
\]
has a positive optimal value.
Since $\dd{p}$ consists of $n$ arcs and the linear program has rational data and encoding length polynomial in the input size, the intersection condition can be verified in polynomial time.
Hence the problem belongs to NP.

We now establish NP-hardness by a reduction from
\textsc{Partition} \cite{karp1972reducibility}.
Consider an instance of \textsc{Partition} problem specified by positive integers $a_1,\ldots,a_n$, where the goal is to determine whether they can be divided into two subsets of equal sum.
Construct an IDD with $n$ arc layers, one node in each node layer, and two parallel arcs, labeled $[0,0]$ and $[1,1]$, between each pair of consecutive node layers.
This IDD has $n+1$ nodes, $2n$ arcs, and width one, with
$\operatorname{Sol}(\dd{D})=\{0,1\}^n$.
Define
\[
    \mathcal{C}
    \coloneqq
    \left\{
        \boldsymbol{x}\in\mathbb{R}^n:
        -\frac{1}{2}
        \le
        \sum_{i=1}^n a_i(2x_i-1)
        \le
        \frac{1}{2}
    \right\}.
\]
The set $\mathcal{C}$ is a polyhedron described by rational linear inequalities and is full-dimensional,
since $(1/2,\ldots,1/2)
\in
\operatorname{int}(\mathcal{C})$.
For every
$\boldsymbol{x}
\in
\operatorname{Sol}(\dd{D})$,
the sum 
$\sum_{i=1}^n a_i(2x_i-1)$ is an integer.
Thus, 
$\boldsymbol{x}
\in
\operatorname{int}(\mathcal{C})$
if and only if
$\sum_{i=1}^n a_i(2x_i-1)=0$,
or equivalently,
$\sum_{i:x_i=1}a_i=\sum_{i:x_i=0}a_i$.
Consequently,
$\operatorname{int}(\mathcal{C})
\cap
\operatorname{Sol}(\dd{D})
\neq
\emptyset$
if and only if the \textsc{Partition} instance has a solution.
This construction can be carried out in polynomial time.
Therefore, the problem of deciding whether
$\operatorname{int}(\mathcal{C})
\cap
\operatorname{Sol}(\dd{D})
\neq
\emptyset$
is NP-complete.
The complementary problem of deciding whether $\mathcal{C}$ is DD-free with respect to $\dd{D}$ is therefore coNP-complete.
\Halmos
\endproof

Although deciding whether a general convex set, such as a polyhedron, is DD-free is computationally difficult, certain families of convex sets induced by $\ell_p$-norm functions admit efficient DD-freeness verification by exploiting the recursive structure of DDs. In what follows, we focus on these families to construct DD-free-set oracles and show that both their construction and DD-freeness verification can be carried out efficiently.

For each $p\in[1,\infty]$, the distance to $\operatorname{Sol}(\dd{D})$ naturally induces the following \emph{$\ell_p$-ball oracle}:
\begin{equation}
\label{eq:lp-ball-oracle}
    \mathcal{C}_p(\boldsymbol{x})
    \coloneqq
    \begin{cases}
        \mathcal{B}_p\bigl(
            \boldsymbol{x},
            d_p(\boldsymbol{x},\operatorname{Sol}(\dd{D}))
        \bigr),
        & \text{if } \boldsymbol{x}\notin\operatorname{Sol}(\dd{D}),\\
        \emptyset,
        & \text{otherwise}.
    \end{cases}
\end{equation}

\begin{proposition}
\label{prop:lp-ball-DD-free}
    For each $p\in[1,\infty]$, the mapping $\mathcal{C}_p(\cdot)$ defined in \eqref{eq:lp-ball-oracle} is a DD-free-set oracle for $\dd{D}$.
\end{proposition}

\proof{Proof.}
Fix $p\in[1,\infty]$.
Property~\textup{(i)} of Definition~\ref{def:DD-free-oracle} follows directly from \eqref{eq:lp-ball-oracle}.
Now let $\boldsymbol{x}\notin\operatorname{Sol}(\dd{D})$ and set $r\coloneqq d_p(\boldsymbol{x},\operatorname{Sol}(\dd{D}))$.
Since $\operatorname{Sol}(\dd{D})$ is nonempty and closed and
$\boldsymbol{x}\notin\operatorname{Sol}(\dd{D})$, we have
$r\in(0,\infty)$. Thus $\mathcal{C}_p(\boldsymbol{x})=\mathcal{B}_p(\boldsymbol{x},r)$ is convex and contains $\boldsymbol{x}$ in its interior.
Furthermore, every $\boldsymbol{s}\in\operatorname{Sol}(\dd{D})$ satisfies $\lVert\boldsymbol{s}-\boldsymbol{x}\rVert_p\ge r$, whereas every $\boldsymbol{y}\in\operatorname{int}(\mathcal{C}_p(\boldsymbol{x}))$ satisfies $\lVert\boldsymbol{y}-\boldsymbol{x}\rVert_p<r$.
Hence $\operatorname{int}(\mathcal{C}_p(\boldsymbol{x}))\cap\operatorname{Sol}(\dd{D})=\emptyset$. Therefore $\mathcal{C}_p(\boldsymbol{x})$ is DD-free, and property~\textup{(ii)} follows.
\Halmos
\endproof

\begin{figure}[tb]
    \centering
    \includegraphics[width=0.65\linewidth]{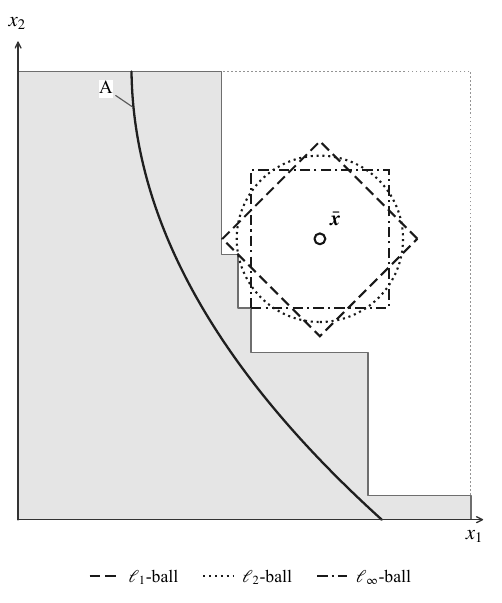}
    \caption{DD-free $\ell_p$-balls centered at       
        $\bar{\boldsymbol{x}}$
        for $p\in\{1,2,\infty\}$.
        The thick solid curve labeled $A$ depicts the nonlinear boundary of $\mathcal{G}^A$, with the feasible region lying on the side containing the origin.
        The shaded region represents
        $\operatorname{Sol}(\dd{D}^A)$.
        Each ball has radius
        $d_p(\bar{\boldsymbol{x}},\operatorname{Sol}(\dd{D}^A))$.
        }
    \label{fig:lp-ball-example}
\end{figure}

We return to the illustrative example in 
Section~\ref{subsec:illustrative-example},
focusing on constraint $A$ with the same IDD relaxation $\dd{D}^A$ and point $\bar{\boldsymbol{x}}$.
As illustrated in Figure~\ref{fig:lp-ball-example}, varying $p$ need not produce a nested family of DD-free sets, because both the geometry of the $\ell_p$-ball and its radius $d_p(\boldsymbol{x},\operatorname{Sol}(\dd{D}))$ can change for different $p$. Consequently, DD-free sets obtained from different norms may be incomparable under set inclusion.

To evaluate the $\ell_p$-ball oracle at a query point $\bar{\boldsymbol{x}}\in\mathbb{R}^n$, we need to compute $d_p(\bar{\boldsymbol{x}},\operatorname{Sol}(\dd{D}))$.
For each $i\in[n]$ and $\dd{a}\in\dd{A}_i$, define the \emph{arc deviation} as the distance from $\bar{x}_i$ to the interval $\dd{l}(\dd{a})$. We can write
\begin{equation}
\label{eq:arc-deviation}
    \xi(\dd{a};\bar{\boldsymbol{x}})
    \coloneqq
    \min_{y\in\dd{l}(\dd{a})}
    |\bar{x}_i-y|
    =
    \max\bigl\{
        \dd{l}(\dd{a})\downarrow-\bar{x}_i,\,
        \bar{x}_i-\dd{l}(\dd{a})\uparrow,\,
        0
    \bigr\}.
\end{equation}
The next proposition shows that $d_p(\bar{\boldsymbol{x}},\operatorname{Sol}(\dd{D}))$ can be computed
through path optimization problems on the IDD using arc deviation values.

\begin{proposition}
\label{prop:distance-to-path}
    Let $\dd{D}=(\dd{U},\dd{A},\dd{l}(\cdot))$ be an IDD with $\operatorname{Sol}(\dd{D})\neq\emptyset$. 
    For every $\boldsymbol{x}\in \mathbb{R}^n$ and every $p\in[1,\infty)$,
    \begin{equation}
    \label{eq:dp-path-finite}
        d_p(\boldsymbol{x}, \operatorname{Sol}(\dd{D}))^p
        =
        \min_{\dd{q}=(\dd{a}_1,\dots,\dd{a}_n)\in\dd{P}(\dd{t})} \sum_{i=1}^n \xi(\dd{a}_i;\boldsymbol{x})^p,
    \end{equation}
    and, for $p=\infty$,
    \begin{equation}
    \label{eq:dp-path-infinite}
        d_{\infty}(\boldsymbol{x}, \operatorname{Sol}(\dd{D}))
        =
        \min_{\dd{q}=(\dd{a}_1,\dots,\dd{a}_n)\in\dd{P}(\dd{t})} \max_{i\in [n]} \ \xi(\dd{a}_i;\boldsymbol{x}).
    \end{equation}
\end{proposition}

\proof{Proof.}
Fix $\boldsymbol{x}\in\mathbb{R}^n$ and an $\dd{r}$--$\dd{t}$ path $\dd{q}=(\dd{a}_1,\dots,\dd{a}_n)\in\dd{P}(\dd{t})$ with induced hyperrectangle $\mathcal{S}_{\dd{q}}=\prod_{i=1}^n \dd{l}(\dd{a}_i)$.
By \eqref{eq:arc-deviation}, for each $i\in[n]$, choose $y_i^\ast\in\dd{l}(\dd{a}_i)$ such that $|x_i-y_i^\ast|=\xi(\dd{a}_i;\boldsymbol{x})$.
Let $\boldsymbol{y}^\ast\coloneqq(y_1^\ast,\dots,y_n^\ast) \in\mathcal{S}_{\dd{q}}$. 
For every $\boldsymbol{y}\in\mathcal{S}_{\dd{q}}$ and $i\in[n]$, it follows from \eqref{eq:arc-deviation} that $|x_i-y_i|\ge\xi(\dd{a}_i;\boldsymbol{x})=|x_i-y_i^\ast|$.
Hence, for every $p\in[1,\infty)$,
\[
    d_p(\boldsymbol{x},\mathcal{S}_{\dd{q}})^p
    =
    \sum_{i=1}^n \xi(\dd{a}_i;\boldsymbol{x})^p,
\]
and
\[
    d_\infty(\boldsymbol{x},\mathcal{S}_{\dd{q}})
    =
    \max_{i\in[n]} \xi(\dd{a}_i;\boldsymbol{x}).
\]

Since $\dd{q}$ was arbitrary, these relations hold for every $\dd{q}\in\dd{P}(\dd{t})$.
Using the fact that $\operatorname{Sol}(\dd{D})=\bigcup_{\dd{q}\in\dd{P}(\dd{t})}\mathcal{S}_{\dd{q}}$ and the finiteness of $\dd{P}(\dd{t})$, we obtain, for every $p\in[1,\infty)$,
\[
    d_p(\boldsymbol{x},\operatorname{Sol}(\dd{D}))^p
    =
    \min_{\dd{q}\in\dd{P}(\dd{t})}
    d_p(\boldsymbol{x},\mathcal{S}_{\dd{q}})^p
    =
    \min_{\dd{q}=(\dd{a}_1,\dots,\dd{a}_n)\in\dd{P}(\dd{t})}
    \sum_{i=1}^n \xi(\dd{a}_i;\boldsymbol{x})^p,
\]
which proves \eqref{eq:dp-path-finite}. Similarly,
\[
    d_\infty(\boldsymbol{x},\operatorname{Sol}(\dd{D}))
    =
    \min_{\dd{q}\in\dd{P}(\dd{t})}
    d_\infty(\boldsymbol{x},\mathcal{S}_{\dd{q}})
    =
    \min_{\dd{q}=(\dd{a}_1,\dots,\dd{a}_n)\in\dd{P}(\dd{t})}
    \max_{i\in[n]}\xi(\dd{a}_i;\boldsymbol{x}),
\]
which proves \eqref{eq:dp-path-infinite}.
\Halmos
\endproof

Proposition~\ref{prop:distance-to-path} reduces the evaluation of the $\ell_p$-ball oracle at a query point $\bar{\boldsymbol{x}}\in\mathbb{R}^n$ to an optimization problem over the paths of the IDD, which can be solved efficiently, as discussed next.
Let $r\coloneqq d_p(\bar{\boldsymbol{x}},\operatorname{Sol}(\dd{D}))$.
For $p\in[1,\infty)$, assign each arc $\dd{a}\in\dd{A}$ the weight $\dd{w}(\dd{a})\coloneqq\xi(\dd{a};\bar{\boldsymbol{x}})^p$. Then $r^p$ is the shortest-path distance from $\dd{r}$ to $\dd{t}$ under these arc weights.
For $p=\infty$, instead set $\dd{w}(\dd{a})\coloneqq\xi(\dd{a};\bar{\boldsymbol{x}})$. In this case, $r$ is the value of a minimax $\dd{r}$--$\dd{t}$ path, i.e., the minimum, over all $\dd{r}$--$\dd{t}$ paths, of the largest arc weight along the path.
Since $\dd{D}$ is layered and acyclic, each path problem can be solved by a single forward pass over the layers.
Algorithm~\ref{alg:lp-ball-DD-free} combines these two cases and computes $\mathcal{C}_p(\bar{\boldsymbol{x}})$ as defined in \eqref{eq:lp-ball-oracle}.

\begin{algorithm}
\caption{Computing the $\ell_p$-ball oracle for an IDD relaxation}
\label{alg:lp-ball-DD-free}
\begin{algorithmic}[1]
    \Require an IDD 
    $\dd{D}=(\dd{U},\dd{A},\dd{l}(\cdot))$ with 
    $\operatorname{Sol}(\dd{D})\neq\emptyset$;
    a point 
    $\bar{\boldsymbol{x}}\in\mathbb{R}^n$;
    and a norm order $p\in[1,\infty]$
    \Ensure $\mathcal{C}_p(\bar{\boldsymbol{x}})$ as defined in
    \eqref{eq:lp-ball-oracle}

    \State initialize $\dd{dist}(\dd{u})\gets+\infty$
    for all $\dd{u}\in\dd{U}\setminus\{\dd{r}\}$ and
    $\dd{dist}(\dd{r})\gets 0$

    \For{$i\in[n]$}
        \For{each arc
        $\dd{a}=(\dd{u},\dd{v})\in\dd{A}_i$}
            \If{$p<\infty$}
                \State
                $\dd{dist}(\dd{v})
                \gets
                \min\bigl\{
                    \dd{dist}(\dd{v}),
                    \dd{dist}(\dd{u})
                    +\xi(\dd{a};\bar{\boldsymbol{x}})^p
                \bigr\}$
            \Else
                \State
                $\dd{dist}(\dd{v})
                \gets
                \min\bigl\{
                    \dd{dist}(\dd{v}),
                    \max\{
                        \dd{dist}(\dd{u}),
                        \xi(\dd{a};\bar{\boldsymbol{x}})
                    \}
                \bigr\}$
            \EndIf
        \EndFor
    \EndFor

    \If{$p<\infty$}
        \State $r\gets\dd{dist}(\dd{t})^{1/p}$
    \Else
        \State $r\gets\dd{dist}(\dd{t})$
    \EndIf

    \If{$r>0$}
        \State\Return $\mathcal{B}_p(\bar{\boldsymbol{x}},r)$
    \Else
        \State\Return $\emptyset$
    \EndIf
\end{algorithmic}
\end{algorithm}

Theorem~\ref{thm:lp-ball-DD-free-valid} establishes that Algorithm~\ref{alg:lp-ball-DD-free} exactly computes the $\ell_p$-ball oracle in~\eqref{eq:lp-ball-oracle} for every $p\in[1,\infty]$, while
Proposition~\ref{prop:lp-ball-DD-free-time-complexity}
shows that the algorithm has running-time complexity linear in the number of arcs in the IDD.

\begin{theorem}
\label{thm:lp-ball-DD-free-valid}
    Let $\dd{D}=(\dd{U},\dd{A},\dd{l}(\cdot))$ be an IDD with $\operatorname{Sol}(\dd{D})\neq\emptyset$. 
    For any $\bar{\boldsymbol{x}}\in\mathbb{R}^n$ and $p\in[1,\infty]$, let $r$ denote the radius computed by Algorithm~\ref{alg:lp-ball-DD-free}. 
    Then $r=d_p(\bar{\boldsymbol{x}},\operatorname{Sol}(\dd{D}))$, and Algorithm~\ref{alg:lp-ball-DD-free} returns $\mathcal{C}_p(\bar{\boldsymbol{x}})$ as defined in \eqref{eq:lp-ball-oracle}.
\end{theorem}

\proof{Proof.}

We first show that the forward pass computes the appropriate optimal
path value. For an $\dd{r}$--$\dd{u}$ path, define its cost as the sum of the $p$-th powers of its arc deviations when $p<\infty$, and as its largest arc deviation when $p=\infty$. The trivial path at $\dd{r}$ has cost zero.

A layer-by-layer induction shows that, before processing $\dd{A}_i$, $\dd{dist}(\dd{u})$ equals the minimum cost of an $\dd{r}$--$\dd{u}$ path for every $\dd{u}\in\dd{U}_i$.
The claim holds for $i=1$ by initialization.
For the induction step, every $\dd{r}$--$\dd{v}$ path to a node $\dd{v}\in\dd{U}_{i+1}$ is obtained by appending an arc $\dd{a}=(\dd{u},\dd{v})\in\dd{A}_i$ to an $\dd{r}$--$\dd{u}$ path. The update in Algorithm~\ref{alg:lp-ball-DD-free} minimizes over
exactly these choices, using addition when $p<\infty$ and the maximum operator when $p=\infty$.
Thus, after processing $\dd{A}_n$, $\dd{dist}(\dd{t})$ equals the
minimum cost of an $\dd{r}$--$\dd{t}$ path.

For $p<\infty$, Proposition~\ref{prop:distance-to-path} therefore gives $\dd{dist}(\dd{t})=d_p(\bar{\boldsymbol{x}},\operatorname{Sol}(\dd{D}))^p$.
For $p=\infty$, the same proposition gives $\dd{dist}(\dd{t})=d_\infty(\bar{\boldsymbol{x}},\operatorname{Sol}(\dd{D}))$.
Hence, the radius assignment yields $r=d_p(\bar{\boldsymbol{x}},\operatorname{Sol}(\dd{D}))$ for all $p\in[1,\infty]$.

Finally, since $\operatorname{Sol}(\dd{D})$ is closed, $r=0$ if and only if $\bar{\boldsymbol{x}}\in\operatorname{Sol}(\dd{D})$. Hence, the algorithm returns $\emptyset$ if $\bar{\boldsymbol{x}}\in\operatorname{Sol}(\dd{D})$ and $\mathcal{B}_p(\bar{\boldsymbol{x}},r)$ otherwise.
Together with $r=d_p(\bar{\boldsymbol{x}},\operatorname{Sol}(\dd{D}))$, this is precisely $\mathcal{C}_p(\bar{\boldsymbol{x}})$ as defined in \eqref{eq:lp-ball-oracle}.
\Halmos
\endproof

\begin{proposition}
\label{prop:lp-ball-DD-free-time-complexity}
    Let $\dd{D}=(\dd{U},\dd{A},\dd{l}(\cdot))$ be an IDD 
    with
    $\operatorname{Sol}(\dd{D})\neq\emptyset$.
    For every fixed $p\in[1,\infty]$ and every
    $\bar{\boldsymbol{x}}\in\mathbb{R}^n$,
    Algorithm~\ref{alg:lp-ball-DD-free} computes
    $\mathcal{C}_p(\bar{\boldsymbol{x}})$ in
    $\mathcal{O}(|\dd{A}|)$ time.
\end{proposition}

\proof{Proof.}

Initializing the distance labels requires
$\mathcal{O}(|\dd{U}|)$ time.
The forward pass processes each arc exactly once.
For each arc, computing its deviation using
\eqref{eq:arc-deviation} and performing the corresponding update take constant time.
Including the iteration over the $n$ arc layers, the total running time is therefore
$\mathcal{O}(|\dd{U}|+|\dd{A}|+n)$.
By \eqref{eq:DD-size-bounds},
$|\dd{U}|\le 2|\dd{A}|$ and $n\le|\dd{A}|$.
Therefore, Algorithm~\ref{alg:lp-ball-DD-free} runs in
$\mathcal{O}(|\dd{A}|)$ time.
\Halmos
\endproof

The above results show that, for any query point $\bar{\boldsymbol{x}}\notin\operatorname{Sol}(\dd{D})$, a DD-free set containing $\bar{\boldsymbol{x}}$ in its interior can be constructed efficiently by evaluating an $\ell_p$-ball oracle.
By Proposition~\ref{prop:DD-free-also-S-free}, every such set is also $S$-free and can therefore serve as a forbidden region in cut-generating procedures for $\mathcal{G}$.
The remainder of this section considers the valid inequalities associated with DD-free sets from two complementary perspectives. We first study the convergence properties of the resulting cut closure and then focus on intersection cuts, a subclass of these inequalities.

\subsection{Convergence of Cut Closures}
\label{subsec:convergence-cut-closures}

Let $\mathcal{P}\subseteq\operatorname{conv}(\mathcal{D})$ be a polytope representing the current outer approximation in the procedure described in Section~\ref{subsec:outer-approximation}.
We first consider a single pair $(\dd{D},\mathcal{C}(\cdot))$ and
show how DD-free sets induce valid separating inequalities for extreme
points of $\mathcal{P}$ lying outside $\operatorname{Sol}(\dd{D})$.
We then consider the oracle family $\{\mathcal{C}^k(\cdot)\}_{k\in K}$ associated with the individual constraint IDDs $\{\dd{D}^k\}_{k\in K}$ and define the cut closure induced jointly by these DD-free-set oracles.
Our main result establishes a sufficient local covering condition on the oracle family under which the infinite-rank cut closure coincides with $\operatorname{conv}\bigl(\mathcal{P}\cap\bigcap_{k\in K}\operatorname{Sol}(\dd{D}^k)\bigr)$.
We then verify this condition for the $\ell_p$-ball oracle family, yielding a concrete convergence guarantee for every fixed $p\in[1,\infty]$.

Consider $\bar{\boldsymbol{x}}\in\operatorname{ext}(\mathcal{P})\setminus\operatorname{Sol}(\dd{D})$.
Since $\mathcal{C}(\bar{\boldsymbol{x}})$ is DD-free, we have
$\operatorname{conv}\bigl(
\mathcal{P}\cap\operatorname{Sol}(\dd{D})
\bigr)
\subseteq
\operatorname{conv}\bigl(
\mathcal{P}\setminus\operatorname{int}(\mathcal{C}(\bar{\boldsymbol{x}}))
\bigr)$.
Moreover, since $\bar{\boldsymbol{x}}$ is an extreme point of $\mathcal{P}$ and belongs to $\operatorname{int}(\mathcal{C}(\bar{\boldsymbol{x}}))$, it does not belong to $\operatorname{conv}\bigl(
\mathcal{P}\setminus\operatorname{int}(\mathcal{C}(\bar{\boldsymbol{x}}))
\bigr)$.
When this set is nonempty, any inequality separating $\bar{\boldsymbol{x}}$ from it is valid for $\operatorname{conv}\bigl(\mathcal{P}\cap\operatorname{Sol}(\dd{D})\bigr)$; see, e.g., \cite[Section~3.1]{bienstock2020outer}.
This separation construction motivates the sequence of cut closures induced jointly by the DD-free-set oracles $\{\mathcal{C}^k(\cdot)\}_{k\in K}$, which we analyze next.

Set $\mathcal{P}_0\coloneqq\mathcal{P}$. For $t\in\mathbb{W}$, we define the \emph{rank-$t$ cut closure} (see \cite{chvatal1973edmonds}) induced by $\{\mathcal{C}^k(\cdot)\}_{k\in K}$ recursively as
\begin{equation}
\label{eq:rank-t-cut-closure}
    \mathcal{P}_{t+1}
    \coloneqq
    \bigcap_{\boldsymbol{x}\in\operatorname{ext}(\mathcal{P}_t)}
    \bigcap_{k\in K}
    \operatorname{conv}\left(
        \mathcal{P}_t
        \setminus
        \operatorname{int}(\mathcal{C}^k(\boldsymbol{x}))
    \right).
\end{equation}
The corresponding \emph{infinite-rank cut closure} is $\mathcal{P}_\infty\coloneqq\bigcap_{t\in\mathbb{W}}\mathcal{P}_t$.

Theorem~\ref{thm:convergence-cut-closure} 
shows that the local covering condition is sufficient for the infinite-rank cut closure to recover the convex hull obtained by intersecting $\mathcal{P}$ with the solution sets of all IDD relaxations for individual constraints.
Proposition~\ref{prop:lp-ball-local-covering-condition} 
shows that the $\ell_p$-ball oracle family satisfies this condition for every fixed $p\in[1,\infty]$.
Corollary~\ref{cor:lp-ball-enlargement-cut-closure-convergence}
extends the resulting convergence guarantee to any DD-free-set oracle family that contains the corresponding $\ell_p$-ball oracles at every query point.

\begin{theorem}[Global convergence of cut closures under a local covering condition]
\label{thm:convergence-cut-closure}
    Let $\{\dd{D}^k\}_{k\in K}$ be IDD relaxations of 
    $\{\mathcal{G}^k\}_{k\in K}$, and, for each $k\in K$, let $\mathcal{C}^k(\cdot)$ be a DD-free-set oracle for $\dd{D}^k$.
    Assume that the following local covering condition holds:
    for every $k\in K$ and every
    $\bar{\boldsymbol{x}}
    \notin\operatorname{Sol}(\dd{D}^k)$, there exists an open
    neighborhood $\mathcal{U}(\bar{\boldsymbol{x}})$ of $\bar{\boldsymbol{x}}$ such that
    \begin{equation}
    \label{eq:local-covering-condition}
        \mathcal{U}(\bar{\boldsymbol{x}})
        \subseteq
        \operatorname{int}
        \bigl(\mathcal{C}^k(\boldsymbol{y})\bigr),
        \qquad
        \forall \boldsymbol{y}\in\mathcal{U}.
    \end{equation}
    Then
    \[
        \mathcal{P}_{\infty}
        =
        \operatorname{conv}\left(
            \mathcal{P}
            \cap
            \bigcap_{k\in K}
            \operatorname{Sol}(\dd{D}^k)
        \right).
    \]
\end{theorem}

\proof{Proof.} 

Let 
$\mathcal{Q}
\coloneqq 
\mathcal{P}\cap \bigcap_{k\in K}\operatorname{Sol}(\dd{D}^k)$. 
We first show that 
$\operatorname{conv}(\mathcal{Q})\subseteq\mathcal{P}_t$ for every $t\in\mathbb{W}$. 
The claim holds for $t=0$ since 
$\mathcal{Q}\subseteq\mathcal{P}_0=\mathcal{P}$ and $\mathcal{P}$ is convex.
Suppose the inclusion holds for some $t\in\mathbb{W}$. For every $\boldsymbol{x}\in\operatorname{ext}(\mathcal{P}_t)$ and every $k\in K$, Definition~\ref{def:DD-free-oracle} implies that 
$\operatorname{int}(\mathcal{C}^k(\boldsymbol{x})) 
\cap
\operatorname{Sol}(\dd{D}^k)=\emptyset$. 
Hence 
$\mathcal{Q}
\subseteq
\mathcal{P}_t\setminus\operatorname{int}(\mathcal{C}^k(\boldsymbol{x}))$, 
and therefore 
$\operatorname{conv}(\mathcal{Q})
\subseteq
\operatorname{conv}\bigl(
    \mathcal{P}_t\setminus\operatorname{int}(\mathcal{C}^k(\boldsymbol{x}))
\bigr)$.
Since this holds for every $\boldsymbol{x}\in\operatorname{ext}(\mathcal{P}_t)$ and every $k\in K$, it follows from \eqref{eq:rank-t-cut-closure} that 
$\operatorname{conv}(\mathcal{Q})\subseteq\mathcal{P}_{t+1}$.
The claim follows by induction, and therefore 
$\operatorname{conv}(\mathcal{Q})\subseteq\mathcal{P}_\infty$.

We next prove the reverse inclusion $\mathcal{P}_\infty\subseteq\operatorname{conv}(\mathcal{Q})$. 
The result holds trivially for the case with $\mathcal{P}_\infty=\emptyset$.
Suppose therefore that $\mathcal{P}_\infty\neq\emptyset$. 
Then every $\mathcal{P}_t$ is nonempty. 
By induction, every $\mathcal{P}_t$ is compact and convex, with
$\mathcal{P}_{t+1}\subseteq\mathcal{P}_t$. 
Indeed, for every $\boldsymbol{x}\in\operatorname{ext}(\mathcal{P}_t)$ and $k\in K$, the set 
$\mathcal{P}_t\setminus\operatorname{int}(\mathcal{C}^k(\boldsymbol{x}))$ is compact, and its convex hull is a compact convex subset of $\mathcal{P}_t$. Thus, $\mathcal{P}_\infty$ is a nonempty compact convex set. 
By Minkowski's theorem \cite[Corollary~1.4.5]{alma99169524416601081}, it suffices to show that 
$\operatorname{ext}(\mathcal{P}_\infty)\subseteq\mathcal{Q}$. 
Suppose, to the contrary, that there exists
$\bar{\boldsymbol{x}}\in 
\operatorname{ext}(\mathcal{P}_\infty)
\setminus\mathcal{Q}$. 
Since $\mathcal{P}_\infty\subseteq\mathcal{P}_0=\mathcal{P}$, there exists $k\in K$ such that $\bar{\boldsymbol{x}}\notin\operatorname{Sol}(\dd{D}^k)$.
Let $\mathcal{U}(\bar{\boldsymbol{x}})$ be the open neighborhood of $\bar{\boldsymbol{x}}$ satisfying the local covering condition \eqref{eq:local-covering-condition}. 
By the cap characterization of extreme points in \cite[Lemma~1.4.6]{alma99169524416601081}, there exists a closed halfspace 
$H^-\coloneqq 
\{\boldsymbol{z}: \boldsymbol{a}^{\top}\boldsymbol{z}\le\beta\}$ 
such that $\boldsymbol{a}^{\top}\bar{\boldsymbol{x}}<\beta$ and $\mathcal{P}_\infty\cap H^-\subseteq\mathcal{U}(\bar{\boldsymbol{x}})$. 
Since $\mathcal{P}_\infty\cap H^-\subseteq\mathcal{U}(\bar{\boldsymbol{x}})$ and the set sequence $(\mathcal{P}_t)_{t\in\mathbb{W}}$ is decreasing with
$\bigcap_{t\in\mathbb{W}}\mathcal{P}_t=\mathcal{P}_\infty$, there exists $t_0\in\mathbb{W}$ such that $\mathcal{P}_{t_0}\cap H^-\subseteq\mathcal{U}(\bar{\boldsymbol{x}})$.
Otherwise, the decreasing compact sets
$(\mathcal{P}_t\cap H^-)\setminus\mathcal{U}(\bar{\boldsymbol{x}})$
would all be nonempty and hence have nonempty intersection, contradicting
$\mathcal{P}_\infty\cap H^-\subseteq\mathcal{U}(\bar{\boldsymbol{x}})$.
Because $\bar{\boldsymbol{x}}\in \mathcal{P}_{t_0}\cap\operatorname{int}(H^-)$, we must have $\min\{\boldsymbol{a}^{\top}\boldsymbol{z} : \boldsymbol{z} \in \mathcal{P}_{t_0}\} < \beta$. 
Since $\mathcal{P}_{t_0}$ is compact, the minimum is attained, and the minimizing face contains an extreme point of $\mathcal{P}_{t_0}$.
Let $\boldsymbol{x}_{t_0}$ be such an extreme point.
Then $\boldsymbol{x}_{t_0}\in \operatorname{ext}(\mathcal{P}_{t_0}) \cap\operatorname{int}(H^-) \subseteq\mathcal{U}(\bar{\boldsymbol{x}})$. 
The local covering condition therefore gives $\mathcal{U}(\bar{\boldsymbol{x}})\subseteq \operatorname{int} (\mathcal{C}^k(\boldsymbol{x}_{t_0}))$. 
Let 
$H^+\coloneqq 
\{\boldsymbol{z}: \boldsymbol{a}^{\top}\boldsymbol{z}\ge\beta\}$. 
Since 
$\mathcal{U}(\bar{\boldsymbol{x}}) \subseteq
\operatorname{int}(\mathcal{C}^k(\boldsymbol{x}_{t_0}))$
and
$\mathcal{P}_{t_0}\cap H^-\subseteq\mathcal{U}(\bar{\boldsymbol{x}})$,
we have 
$\mathcal{P}_{t_0}
\setminus
\operatorname{int}(\mathcal{C}^k(\boldsymbol{x}_{t_0}))
\subseteq
\mathcal{P}_{t_0}\setminus\mathcal{U}(\bar{\boldsymbol{x}})
\subseteq
H^+$.
Since $H^+$ is convex, the relation \eqref{eq:rank-t-cut-closure} implies $\mathcal{P}_{t_0+1}\subseteq H^+$. 
However, 
$\bar{\boldsymbol{x}}
\in 
\mathcal{P}_\infty
\subseteq
\mathcal{P}_{t_0+1}$ and $\bar{\boldsymbol{x}}\in\operatorname{int}(H^-)$, yielding
a contradiction. 
Therefore, $\operatorname{ext}(\mathcal{P}_\infty)\subseteq\mathcal{Q}$, and Minkowski's theorem yields $\mathcal{P}_\infty\subseteq \operatorname{conv}(\mathcal{Q})$. Combining this with the fact that
$\operatorname{conv}(\mathcal{Q})\subseteq\mathcal{P}_\infty$ gives the desired equality.
\Halmos 
\endproof

\begin{proposition}
\label{prop:lp-ball-local-covering-condition}
    For every fixed $p\in[1,\infty]$, the family of $\ell_p$-ball
    oracles $\{\mathcal{C}_p^k(\cdot)\}_{k\in K}$ associated with
    $\{\dd{D}^k\}_{k\in K}$ satisfies the local covering condition
    \eqref{eq:local-covering-condition}.
\end{proposition}

\proof{Proof.}

Fix $p\in[1,\infty]$, $k\in K$, and 
$\bar{\boldsymbol{x}}\notin\operatorname{Sol}(\dd{D}^k)$, and write
$\mathcal{T}\coloneqq\operatorname{Sol}(\dd{D}^k)$.
Set $r\coloneqq d_p(\bar{\boldsymbol{x}},\mathcal{T})$.
Since $\mathcal{T}$ is nonempty and closed and $\bar{\boldsymbol{x}}\notin\mathcal{T}$, we have $r\in(0,\infty)$.
Let $\rho\coloneqq r/3$ and 
$\mathcal{U}(\bar{\boldsymbol{x}}) \coloneqq
\operatorname{int}(\mathcal{B}_p(\bar{\boldsymbol{x}},\rho))$, 
which is an open neighborhood of $\bar{\boldsymbol{x}}$.

For any $\boldsymbol{y}\in\mathcal{U}(\bar{\boldsymbol{x}})$, let $\boldsymbol{t}(\boldsymbol{y}) \in \mathcal{T}$ be the point that yields the minimum distance $d_p(\boldsymbol{y},\mathcal{T})$. Then, we can write  
$d_p(\boldsymbol{y},\mathcal{T})
\ge \lVert\boldsymbol{t}(\boldsymbol{y})-\bar{\boldsymbol{x}}\rVert_p - \lVert\boldsymbol{y}-\bar{\boldsymbol{x}}\rVert_p \geq
r-\lVert\boldsymbol{y}-\bar{\boldsymbol{x}}\rVert_p
>
r-\rho
=
2\rho
>
0$,
where the first inequality follows from the triangle inequality, the second inequality holds by definition of $r$, and the third inequality is due to $\boldsymbol{y}\in\mathcal{U}(\bar{\boldsymbol{x}})$.
Hence $\boldsymbol{y}\notin\mathcal{T}$, and \eqref{eq:lp-ball-oracle} yields 
$\mathcal{C}_p^k(\boldsymbol{y})
=
\mathcal{B}_p\bigl(
\boldsymbol{y},
d_p(\boldsymbol{y},\mathcal{T})
\bigr)$.
Now consider any $\boldsymbol{z}\in\mathcal{U}(\bar{\boldsymbol{x}})$. We can write
$\lVert\boldsymbol{z}-\boldsymbol{y}\rVert_p
\leq 
\lVert\boldsymbol{z}-\bar{\boldsymbol{x}}\rVert_p + \lVert\boldsymbol{y}-\bar{\boldsymbol{x}}\rVert_p
<
2\rho
<
d_p(\boldsymbol{y},\mathcal{T})$,
where the first inequality follows from the triangle inequality, the second inequality is due to $\boldsymbol{y}\in\mathcal{U}(\bar{\boldsymbol{x}})$ and $\boldsymbol{z}\in\mathcal{U}(\bar{\boldsymbol{x}})$, and the last inequality holds because of the relation derived above. 
Therefore,
$\boldsymbol{z}\in
\operatorname{int}(\mathcal{C}_p^k(\boldsymbol{y}))$.

Since $\boldsymbol{y}\in\mathcal{U}(\bar{\boldsymbol{x}})$ was arbitrary and the preceding argument applies to every $\boldsymbol{z}\in\mathcal{U}(\bar{\boldsymbol{x}})$, we obtain
$\mathcal{U}(\bar{\boldsymbol{x}})\subseteq
\operatorname{int}(\mathcal{C}_p^k(\boldsymbol{y}))$
for every $\boldsymbol{y}\in\mathcal{U}(\bar{\boldsymbol{x}})$.
Thus, \eqref{eq:local-covering-condition} holds.
\Halmos
\endproof

\begin{corollary}
\label{cor:lp-ball-enlargement-cut-closure-convergence}
    For every fixed $p\in[1,\infty]$, let
    $\{\widehat{\mathcal{C}}^k(\cdot)\}_{k\in K}$ be a family of
    DD-free-set oracles associated with
    $\{\dd{D}^k\}_{k\in K}$ satisfying
    \[
        \mathcal{C}_p^k(\boldsymbol{x})
        \subseteq
        \widehat{\mathcal{C}}^k(\boldsymbol{x}),
        \qquad
        \forall k\in K,\;
        \boldsymbol{x}\in\mathbb{R}^n.
    \]
    Then
    $\{\widehat{\mathcal{C}}^k(\cdot)\}_{k\in K}$
    satisfies the local covering condition
    \eqref{eq:local-covering-condition}.
    Consequently, if
    $\{\widehat{\mathcal{P}}_t\}_{t\in\mathbb{W}}$
    denotes the cut closures in
    \eqref{eq:rank-t-cut-closure} induced by this oracle family, then
    \[
        \widehat{\mathcal{P}}_\infty
        =
        \operatorname{conv}\left(
            \mathcal{P}
            \cap
            \bigcap_{k\in K}
            \operatorname{Sol}(\dd{D}^k)
        \right).
    \]
\end{corollary}

\proof{Proof.}

By Proposition~\ref{prop:lp-ball-local-covering-condition}, the family
$\{\mathcal{C}_p^k(\cdot)\}_{k\in K}$ satisfies the local covering
condition \eqref{eq:local-covering-condition}.
The inclusion assumption implies
$\operatorname{int}(\mathcal{C}_p^k(\boldsymbol{x}))
\subseteq
\operatorname{int}(\widehat{\mathcal{C}}^k(\boldsymbol{x}))$
for every $k\in K$ and $\boldsymbol{x}\in\mathbb{R}^n$.
Hence,
$\{\widehat{\mathcal{C}}^k(\cdot)\}_{k\in K}$
also satisfies \eqref{eq:local-covering-condition}.
The desired equality then follows from
Theorem~\ref{thm:convergence-cut-closure}.
\Halmos
\endproof


An important implication of the global convergence property established above is the ability to derive strong cuts that can further tighten the relaxations obtained from DD-based outer approximation methods for MINLPs, even when individual constraints are treated separately.
To illustrate, a direct approach to deriving valid inequalities for $\operatorname{conv}\bigl(\bigcap_{k\in K}\operatorname{Sol}(\dd{D}^k)\bigr)$ is to explicitly represent the intersection of the solution sets of the individual IDDs.
Specifically, one may construct an intersection IDD, denoted by
$\dd{D}^{\cap}$, satisfying
$\operatorname{Sol}(\dd{D}^{\cap})
=
\bigcap_{k\in K}\operatorname{Sol}(\dd{D}^k)$
and then apply DD-based cut-generation methods, as described in
Section~\ref{subsec:outer-approximation}, to generate IDD cuts valid for $\operatorname{conv}(\operatorname{Sol}(\dd{D}^{\cap}))$.
However, explicitly constructing $\dd{D}^{\cap}$ can be computationally expensive. Under a common variable ordering, the width of an intersection IDD can be as large as the product of the widths of the individual IDDs \cite[Section~9.3.3]{bergman2016decision}. This worst-case width growth can substantially increase the size of the resulting diagram and the computational cost of subsequent cut generation.

In contrast, Theorem~\ref{thm:convergence-cut-closure} shows that, under the local covering condition, cuts induced by DD-free-set oracles for the individual IDDs can be generated separately and, when combined through the cut closure procedure, collectively recover the convexification characterized above without explicitly constructing an intersection IDD.
The computational results in Section~\ref{sec:computational-results} complement the preceding theoretical analysis.
In particular, they show that the outer approximation augmented with the proposed intersection cuts, generated from DD-free sets associated with individual IDDs, can yield stronger bounds than those obtained by global solvers.


\subsection{Classical Intersection Cuts Derivation}
\label{subsec:IC-derivation}

The cut-closure analysis above does not prescribe how the separating inequalities are generated. Exact separation over
$\operatorname{conv}\bigl(
    \mathcal{P}\setminus\operatorname{int}(\mathcal{C})
\bigr)$ 
can be computationally demanding; in particular, the associated separation problem is NP-hard even when $\mathcal{C}$ is a ball \cite{freund1985complexity,bienstock2020outer}.
We therefore use the classical intersection-cut construction introduced by Balas \cite{balas1971intersection}, applied to the corner relaxation induced by a fixed LP basis; see also \cite[Section~6]{confortiIntegerProgramming2014}.

Let 
$\bar{\boldsymbol{x}}
\in
\mathcal{P}\setminus\operatorname{Sol}(\dd{D})$ be a basic feasible
solution, and fix an associated LP basis.
Let 
$\mathcal{P}'
\coloneqq
\bar{\boldsymbol{x}}
+
\operatorname{cone}\{\boldsymbol{v}^j:j\in R\}$
denote the resulting simplicial cone relaxation of $\mathcal{P}$, where $R$ indexes the extreme rays of $\mathcal{P}'$ and $\boldsymbol{v}^j$ is a direction vector for ray $j$.
By construction, $\mathcal{P}\subseteq\mathcal{P}'$.
Since $\mathcal{P}'$ is simplicial, every 
$\boldsymbol{x}\in\mathcal{P}'$ admits a unique representation
$\boldsymbol{x}
=
\bar{\boldsymbol{x}}
+
\sum_{j\in R}
\lambda_j(\boldsymbol{x})\boldsymbol{v}^j$,
where $\lambda_j(\boldsymbol{x})\ge 0$ for all $j\in R$.

Let $\mathcal{C}\subseteq\mathbb{R}^n$ be a DD-free set with respect to
$\dd{D}$ such that
$\bar{\boldsymbol{x}}\in\operatorname{int}(\mathcal{C})$.
For each $j\in R$, define the \emph{step length}
\[
    \alpha_j
    \coloneqq
    \sup
    \bigl\{
        t\ge 0:
        \bar{\boldsymbol{x}}
        +
        t\boldsymbol{v}^j
        \in\mathcal{C}
    \bigr\}
    \in(0,+\infty].
\]
Thus, $\alpha_j=+\infty$ when the entire ray
$\{\bar{\boldsymbol{x}}+t\boldsymbol{v}^j:t\ge 0\}$
is contained in $\mathcal{C}$.
With the convention $1/(+\infty)\coloneqq 0$, the classical \textit{intersection cut} is given by
\[
    \sum_{j\in R}
    \frac{\lambda_j(\boldsymbol{x})}{\alpha_j}
    \ge 1.
\]
This inequality is valid for
$\mathcal{P}'\setminus\operatorname{int}(\mathcal{C})$.
Since $\mathcal{C}$ is DD-free and
$\mathcal{P}\subseteq\mathcal{P}'$, we have
$\mathcal{P}\cap\operatorname{Sol}(\dd{D})
\subseteq
\mathcal{P}'
\setminus
\operatorname{int}(\mathcal{C})$.
The cut is therefore valid for
$\mathcal{P}\cap\operatorname{Sol}(\dd{D})$ and for
$\mathcal{P}\cap\mathcal{G}$ as $\mathcal{G}\subseteq\operatorname{Sol}(\dd{D})$.

For later use, it is convenient to express the cut coefficients through a gauge function. 
Define
$\mathcal{K}\coloneqq\mathcal{C}-\bar{\boldsymbol{x}}$, so that
$\boldsymbol{0}\in\operatorname{int}(\mathcal{K})$, and let
$\gamma_{\mathcal{K}}(\boldsymbol{w})
\coloneqq
\inf\{t>0:\boldsymbol{w}\in t\mathcal{K}\}$
denote the gauge of $\mathcal{K}$.
The definition of $\alpha_j$ gives
$\gamma_{\mathcal{K}}(\boldsymbol{v}^j)=1/\alpha_j$ for every $j\in R$.
Hence, the intersection cut can equivalently be written as
\begin{equation}
\label{eq:IC-gauge}
    \sum_{j\in R}
    \gamma_{\mathcal{K}}(\boldsymbol{v}^j)
    \lambda_j(\boldsymbol{x})
    \ge 1.
\end{equation}

We finally specialize \eqref{eq:IC-gauge} to the DD-free set generated by the $\ell_p$-ball oracle. 
Let
$r\coloneqq
d_p(\bar{\boldsymbol{x}},\operatorname{Sol}(\dd{D}))$.
Since
$\bar{\boldsymbol{x}}\notin\operatorname{Sol}(\dd{D})$, 
we have $r>0$, and
$\mathcal{C}_p(\bar{\boldsymbol{x}})
=
\mathcal{B}_p(\bar{\boldsymbol{x}},r)$
by~\eqref{eq:lp-ball-oracle}.
It follows that
$\mathcal{K}=\mathcal{B}_p(\boldsymbol{0},r)$ and
$\gamma_{\mathcal{K}}(\boldsymbol{w})
=
\lVert\boldsymbol{w}\rVert_p/r$.
Substituting into~\eqref{eq:IC-gauge} yields the
$\ell_p$-ball intersection cut
\begin{equation}
\label{eq:lp-ball-IC}
    \sum_{j\in R}
    \frac{\lVert\boldsymbol{v}^j\rVert_p}{r}
    \lambda_j(\boldsymbol{x})
    \ge 1.
\end{equation}

\begin{figure}[tb]
    \centering
    \includegraphics[width=\linewidth]{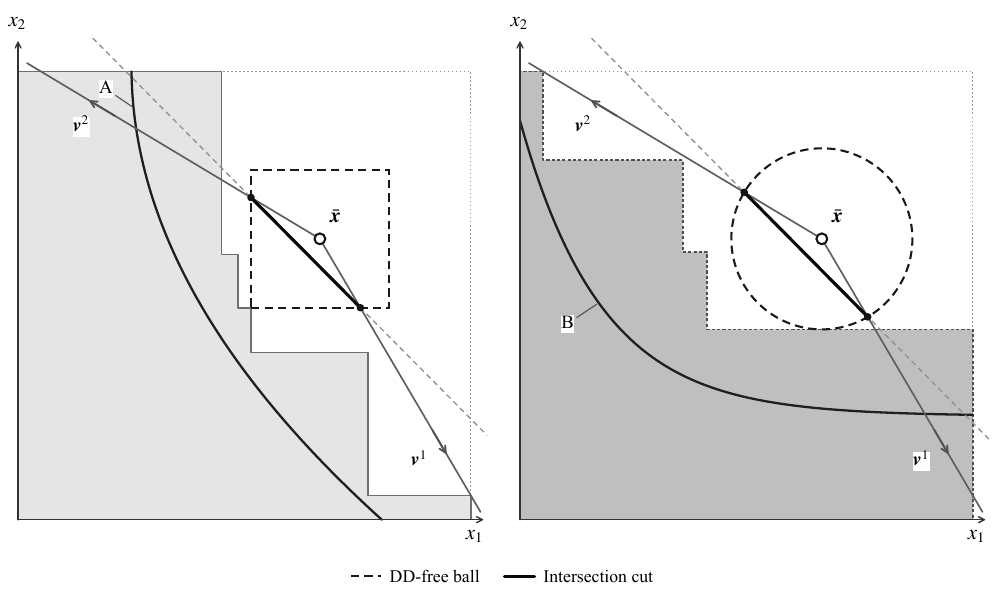}
    \caption{Intersection cuts generated separately from DD-free
        $\ell_\infty$-ball for constraints $A$ (left) and $\ell_2$-ball for constraint $B$ (right).
        The solid curves labeled $A$ and $B$ depict the nonlinear boundaries of the corresponding constraint sets, with each feasible region lying on the side containing the origin.
        The shaded regions represent
        $\operatorname{Sol}(\dd{D}^A)$ (left) and
        $\operatorname{Sol}(\dd{D}^B)$ (right).
        The gray rays with directions $\boldsymbol{v}^1$
        and $\boldsymbol{v}^2$ are the extreme rays of the common corner relaxation $\mathcal{P}'$.
        In each panel, the black dashed lines/curves bound the
        DD-free balls centered at $\bar{\boldsymbol{x}}$,
        and the two filled dots mark the intersections of the rays with the ball boundary.
        The thick black segment and its gray dashed extensions show the cut boundary, with the retained halfspace lying on the side opposite to $\bar{\boldsymbol{x}}$.
    }
    \label{fig:idd-intersection-cuts}
\end{figure}

To illustrate, consider the example in
Section~\ref{subsec:illustrative-example}, with the same IDD relaxations $\dd{D}^A$ and $\dd{D}^B$ and point $\bar{\boldsymbol{x}}$.
Take $\mathcal{P}=\mathcal{P}_{\mathrm{IDD}}$ and form the corner relaxation $\mathcal{P}'$ from the two active inequalities whose boundaries are $H_A$ and $H_B$.
Figure~\ref{fig:idd-intersection-cuts} illustrates the
$\ell_\infty$-ball and $\ell_2$-ball intersection cuts obtained separately from the two IDDs on this common corner relaxation.
For constraint $A$, we use the DD-free $\ell_\infty$-ball
$\mathcal{B}_\infty\left(
    \bar{\boldsymbol{x}},
    d_\infty(\bar{\boldsymbol{x}},\operatorname{Sol}(\dd{D}^A))
\right)$, and for constraint $B$, we use the DD-free $\ell_2$-ball
$\mathcal{B}_2\left(
    \bar{\boldsymbol{x}},
    d_2(\bar{\boldsymbol{x}},\operatorname{Sol}(\dd{D}^B))
\right)$.
The corresponding cut given by \eqref{eq:lp-ball-IC} is valid for
$\mathcal{P}_{\mathrm{IDD}}
\cap
\operatorname{Sol}(\dd{D}^k)$,
and hence for $\mathcal{G}^A\cap\mathcal{G}^B$.
Although no IDD cut generated from either fixed diagram can separate $\bar{\boldsymbol{x}}$, each of the two intersection cuts separates this point and thereby strictly tightens $\mathcal{P}_{\mathrm{IDD}}$.

\section{Maximal DD-Free Sets}
\label{sec:maximal-DD-free}

The preceding section developed intersection cuts from DD-free sets.
For a fixed corner relaxation $\mathcal{P}'$ with apex $\bar{\boldsymbol{x}}$, let $\mathcal{C}_1\subseteq\mathcal{C}_2$ be DD-free sets with respect to
the same IDD $\dd{D}$, and suppose that
$\bar{\boldsymbol{x}}\in\operatorname{int}(\mathcal{C}_1)$.
Define $\mathcal{K}_1\coloneqq \mathcal{C}_1-\bar{\boldsymbol{x}}$ and $\mathcal{K}_2\coloneqq \mathcal{C}_2-\bar{\boldsymbol{x}}$. 
The set inclusion implies that $\gamma_{\mathcal{K}_2}
(\boldsymbol{v}^j)
\le
\gamma_{\mathcal{K}_1}
(\boldsymbol{v}^j)$ for every $j\in R$.
Hence, by \eqref{eq:IC-gauge}, the intersection cut generated from $\mathcal{C}_2$ is at least as strong on $\mathcal{P}'$ as that generated from $\mathcal{C}_1$.
This monotonicity motivates the study of inclusion-wise maximal DD-free sets, following the general notion of maximality for $S$-free sets; see, e.g., \cite{conforti2015cut,bienstock2020outer}.

Section~\ref{subsec:general-maximal-DD-free} establishes general maximality results for DD-free sets, analyzes their computational complexity, and discusses how these derivations can be facilitated by focusing on structured families based on $\ell_p$-balls, building on the results of the previous section. We then focus on two structured families: boxes, which include $\ell_{\infty}$-balls, in  Section~\ref{subsec:box-maximal-DD-free}, and $\ell_p$-balls for $p\in[1,\infty)$ in Section~\ref{subsec:lp-ball-maximal-DD-free}. For each family, we characterize maximality and develop corresponding enlargement procedures.
For the remainder of this section, we assume that the IDDs under consideration have nonempty solution sets.


\subsection{General Maximality Results}
\label{subsec:general-maximal-DD-free}

We begin by defining maximality for DD-free sets.
Given an IDD $\dd{D}$, we say that a convex set $\mathcal{C}$ is \textit{maximal with respect to} $\dd{D}$ if there is no DD-free set $\mathcal{C}'$ such that $\mathcal{C}\subsetneq\mathcal{C}'$. 
We first show that every full-dimensional maximal DD-free set is a polyhedron.

\begin{proposition}
\label{prop:polyhedral-maximal-DD-free}
    Let $\dd{D}$ be an IDD.
    Every full-dimensional maximal DD-free set with respect to
    $\dd{D}$ is a polyhedron that can be described by at most
    $|\dd{P}(\dd{t})|$ linear inequalities.
\end{proposition}

\proof{Proof.}

Let $\mathcal{C}$ be a full-dimensional maximal DD-free set with respect to $\dd{D}$.
For each 
$\dd{p}\in\dd{P}(\dd{t})$, the sets
$\operatorname{int}(\mathcal{C})$ and $\mathcal{S}_{\dd{p}}$
are nonempty, convex, and disjoint.
By the separating hyperplane theorem, there exists a closed halfspace $H_{\dd{p}}$ such that
$\operatorname{int}(\mathcal{C})
\subseteq 
H_{\dd{p}}$ and
$\operatorname{int}(H_{\dd{p}})
\cap
\mathcal{S}_{\dd{p}}
=\emptyset$.
Since $\mathcal{C}$ is full-dimensional and convex, we have
$\mathcal{C}
\subseteq
\operatorname{cl}(\operatorname{int}(\mathcal{C}))
\subseteq 
H_{\dd{p}}$.
Let
$\mathcal{Q}
\coloneqq
\bigcap_{
    \dd{p}\in\dd{P}(\dd{t})
}H_{\dd{p}}$.
Because $\dd{P}(\dd{t})$ is finite, $\mathcal{Q}$ is a
polyhedron containing $\mathcal{C}$ and is described by at most $|\dd{P}(\dd{t})|$ linear inequalities.
Moreover,
$\operatorname{int}(\mathcal{Q})
\subseteq
\operatorname{int}(H_{\dd{p}})$
for every $\dd{p}\in\dd{P}(\dd{t})$, so
$\operatorname{int}(\mathcal{Q})
\cap
\mathcal{S}_{\dd{p}}
=\emptyset$.
Since
$\operatorname{Sol}(\dd{D})
=\bigcup_{
    \dd{p}\in\dd{P}(\dd{t})
}\mathcal{S}_{\dd{p}}$,
the polyhedron $\mathcal{Q}$ is DD-free.
Maximality of $\mathcal{C}$ implies that
$\mathcal{C}=\mathcal{Q}$.
\Halmos
\endproof

For a closed DD-free set $\mathcal{C}$, a point
$\boldsymbol{z}
\in
\operatorname{Sol}(\dd{D})
\cap
\operatorname{bd}(\mathcal{C})$
is called a \emph{contact point} of $\mathcal{C}$.
We next characterize maximality for full-dimensional DD-free polyhedra through a facet-contact condition.

\begin{proposition}
\label{prop:polyhedral-maximal-characterization}
    Let $\dd{D}$ be an IDD, and let
    $\mathcal{C}\subseteq\mathbb{R}^n$ be a full-dimensional DD-free polyhedron with respect to $\dd{D}$.
    Then $\mathcal{C}$ is maximal if and only if every facet $F$ of $\mathcal{C}$ satisfies
    $\operatorname{relint}(F)
    \cap
    \operatorname{Sol}(\dd{D})
    \neq
    \emptyset$.
\end{proposition}

\proof{Proof.}

Write
$\mathcal{C}
=
\{\boldsymbol{x}\in\mathbb{R}^n:
\boldsymbol{A}\boldsymbol{x}
\le
\boldsymbol{\beta}\}$,
where the description is irredundant and
$\boldsymbol{A}$ has $m$ rows.
For each $j\in[m]$, 
let $\boldsymbol{a}_j^\top$ denote the $j$-th row of $\boldsymbol{A}$ and define
$F_j
\coloneqq
\{\boldsymbol{x}\in\mathcal{C}:
\boldsymbol{a}_j^\top\boldsymbol{x}=\beta_j\}$.

Suppose first that $\mathcal{C}$ is maximal.
Assume, to the contrary, that
$\operatorname{relint}(F_i)
\cap
\operatorname{Sol}(\dd{D})
=
\emptyset$
for some $i\in[m]$, and choose
$\boldsymbol{z}
\in
\operatorname{relint}(F_i)$.
We construct a DD-free polyhedron $\mathcal{Q}$ containing $\mathcal{C}$ such that
$\boldsymbol{z}
\in
\operatorname{int}(\mathcal{Q})$.
Since $\operatorname{Sol}(\dd{D})$ is a finite union of path-induced hyperrectangles, it suffices to construct,
for each $\dd{p}\in\dd{P}(\dd{t})$,
a closed halfspace $H_{\dd{p}}$ satisfying
$\mathcal{C}\subseteq H_{\dd{p}}$,
$\boldsymbol{z}
\in
\operatorname{int}(H_{\dd{p}})$,
and
$\operatorname{int}(H_{\dd{p}})
\cap
\mathcal{S}_{\dd{p}}
=
\emptyset$.
Fix $\dd{p}\in\dd{P}(\dd{t})$.
If
$\mathcal{S}_{\dd{p}}
\cap
\{
    \boldsymbol{x}:
    \boldsymbol{a}_i^\top\boldsymbol{x}
    \le 
    \beta_i
\}
=\emptyset$,
compactness of $\mathcal{S}_{\dd{p}}$ implies that
$\delta_{\dd{p}}
\coloneqq
\min_{
    \boldsymbol{x}\in\mathcal{S}_{\dd{p}}
}\boldsymbol{a}_i^\top\boldsymbol{x}>\beta_i$.
Then
$H_{\dd{p}}
\coloneqq
\{
    \boldsymbol{x}\in\mathbb{R}^n:
    \boldsymbol{a}_i^\top\boldsymbol{x}
    \le
    \delta_{\dd{p}}
\}$
has the required properties.
If
$\mathcal{S}_{\dd{p}}
\cap
\{
    \boldsymbol{x}:
    \boldsymbol{a}_i^\top\boldsymbol{x}
    \le 
    \beta_i
\}
\neq\emptyset$,
then we write
$\mathcal{S}_{\dd{p}}
=
\{
    \boldsymbol{x}\in\mathbb{R}^n:
    \boldsymbol{l}_{\dd{p}}
    \le
    \boldsymbol{x}
    \le
    \boldsymbol{u}_{\dd{p}}
\}$
and consider the linear program
\[
    \max_{\boldsymbol{x},\tau}
    \left\{
        \tau:
        \boldsymbol{x}
        \in
        \mathcal{S}_{\dd{p}},\
        \boldsymbol{a}_i^\top\boldsymbol{x}
        \le
        \beta_i,\
        \boldsymbol{a}_j^\top\boldsymbol{x}
        +\tau
        \le
        \beta_j
        \quad
        \forall j\in[m]\setminus\{i\}
    \right\}.
\]
In this case, $m\ge2$; otherwise, 
$\mathcal{S}_{\dd{p}}
\cap
\{
    \boldsymbol{x}:
    \boldsymbol{a}_i^\top\boldsymbol{x}
    \le 
    \beta_i
\}
\neq\emptyset$
would contradict the fact that $\mathcal{C}$ is DD-free or the assumption that $\operatorname{relint}(F_i)
\cap
\operatorname{Sol}(\dd{D})
=
\emptyset$.
The above linear program is feasible because
$\mathcal{S}_{\dd{p}}$ intersects
$\{
    \boldsymbol{x}:
    \boldsymbol{a}_i^\top\boldsymbol{x}
    \le
    \beta_i
\}$
and $\tau$ can be chosen sufficiently small.
Its optimal value is at most zero, since a feasible
solution with $\tau>0$ would imply
$\boldsymbol{x}
\in
\operatorname{int}(\mathcal{C})
\cup
\operatorname{relint}(F_i)$,
contradicting DD-freeness of $\mathcal{C}$ or the assumption that $\operatorname{relint}(F_i)
\cap
\operatorname{Sol}(\dd{D})
=
\emptyset$.
The dual linear program is
\[
\begin{aligned}
    \min_{
        \boldsymbol{\mu}^+,
        \boldsymbol{\mu}^-,
        \boldsymbol{\lambda}
    }\quad&
    (\boldsymbol{\mu}^+)^\top\boldsymbol{u}_{\dd{p}}
    -(\boldsymbol{\mu}^-)^\top\boldsymbol{l}_{\dd{p}}
    +\sum_{j\in[m]}\lambda_j\beta_j\\
    \text{s.t.}\quad&
    \boldsymbol{\mu}^+-\boldsymbol{\mu}^-
    +\sum_{j\in[m]}\lambda_j\boldsymbol{a}_j
    =\boldsymbol{0},\\
    &\sum_{j\in[m]\setminus\{i\}}\lambda_j=1,\\
    &\boldsymbol{\mu}^+,\boldsymbol{\mu}^-
    \ge\boldsymbol{0},\
    \lambda_j\ge0\quad\forall j\in[m].
\end{aligned}
\]
By strong duality, the dual admits an optimal solution with the same nonpositive objective value.
Let
$(
    \boldsymbol{\mu}^+,
    \boldsymbol{\mu}^-,
    \boldsymbol{\lambda}
)$
be such a solution.
For every 
$\boldsymbol{s}
\in
\mathcal{S}_{\dd{p}}$,
we have
\begin{equation}
\label{eq:dual-inequality}
    \sum_{j\in[m]}\lambda_j
    \boldsymbol{a}_j^\top\boldsymbol{s}
    =
    (\boldsymbol{\mu}^-)^\top\boldsymbol{s}
    -(\boldsymbol{\mu}^+)^\top\boldsymbol{s}
    \ge
    (\boldsymbol{\mu}^-)^\top\boldsymbol{l}_{\dd{p}}
    -(\boldsymbol{\mu}^+)^\top\boldsymbol{u}_{\dd{p}}
    \ge
    \sum_{j\in[m]}\lambda_j\beta_j,
\end{equation}
where the equality follows from the first dual equality constraint, the first inequality from
$\boldsymbol{l}_{\dd{p}}
\le
\boldsymbol{s}
\le
\boldsymbol{u}_{\dd{p}}$
and 
$\boldsymbol{\mu}^+,\boldsymbol{\mu}^-
\ge
\boldsymbol{0}$,
and the second inequality from the nonpositive dual objective value.
Define
$H_{\dd{p}}
\coloneqq
\{
    \boldsymbol{x}\in\mathbb{R}^n:
    \sum_{j\in[m]}
    \lambda_j(
        \boldsymbol{a}_j^\top\boldsymbol{x}
        -\beta_j
    )
    \le0
\}$.
Since $\lambda_j\ge0$ for all $j\in[m]$,
we have $\mathcal{C}\subseteq H_{\dd{p}}$,
and \eqref{eq:dual-inequality} gives
$\operatorname{int}(H_{\dd{p}})
\cap
\mathcal{S}_{\dd{p}}
=\emptyset$.
Moreover, since
$\boldsymbol{z}
\in
\operatorname{relint}(F_i)$,
we have
\[
    \sum_{j\in[m]}
    \lambda_j
    (
        \boldsymbol{a}_j^\top\boldsymbol{z}
        -\beta_j
    )
    =
    \sum_{j\in[m]\setminus\{i\}}
    \lambda_j
    (
        \boldsymbol{a}_j^\top\boldsymbol{z}
        -\beta_j
    )
    <0,
\]
where the equality is because
$\boldsymbol{a}_i^\top\boldsymbol{z}
=\beta_i$,
and the strict inequality follows from
$\boldsymbol{a}_j^\top\boldsymbol{z}<\beta_j$
and $\lambda_j\ge0$ for all
$j\in[m]\setminus\{i\}$,
together with
$\sum_{j\in[m]\setminus\{i\}}\lambda_j=1$.
Hence, 
$\boldsymbol{z}
\in
\operatorname{int}(H_{\dd{p}})$.

Let
$\mathcal{Q}
\coloneqq
\bigcap_{\dd{p}\in\dd{P}(\dd{t})}
H_{\dd{p}}$.
Since $\dd{P}(\dd{t})$ is finite,
$\mathcal{Q}$ is DD-free and satisfies
$\mathcal{C}\subseteq\mathcal{Q}$ and
$\boldsymbol{z}
\in
\operatorname{int}(\mathcal{Q})$.
Because 
$\boldsymbol{z}
\in
\operatorname{bd}(\mathcal{C})$,
we have $\mathcal{C}\subsetneq\mathcal{Q}$,
contradicting maximality.

Conversely, suppose that every facet of $\mathcal{C}$ satisfies the stated condition.
Assume, to the contrary, that there exists a DD-free set
$\mathcal{C}'$ with
$\mathcal{C}\subsetneq\mathcal{C}'$.
Choose
$\boldsymbol{y}
\in
\mathcal{C}'
\setminus
\mathcal{C}$,
and let $i\in[m]$ satisfy
$\boldsymbol{a}_i^\top\boldsymbol{y}
>\beta_i$.
By assumption, there exists
$\boldsymbol{z}
\in
\operatorname{relint}(F_i)
\cap
\operatorname{Sol}(\dd{D})$.
Therefore, since
$\boldsymbol{z}
\in
\operatorname{relint}(F_i)$,
for sufficiently small $\varepsilon>0$, the point
$\boldsymbol{q}
\coloneqq
\boldsymbol{z}
+
\varepsilon(\boldsymbol{z}-\boldsymbol{y})$
belongs to $\operatorname{int}(\mathcal{C})$.
Moreover,
$\boldsymbol{z}
=
\frac{1}{1+\varepsilon}\boldsymbol{q}
+
\frac{\varepsilon}{1+\varepsilon}\boldsymbol{y}$.
Since
$\boldsymbol{q}\in
\operatorname{int}(\mathcal{C})
\subseteq\operatorname{int}(\mathcal{C}')$,
$\boldsymbol{y}\in\mathcal{C}'$, and
$\mathcal{C}'$ is convex,
we have
$\boldsymbol{z}\in
\{\boldsymbol{t}:
\boldsymbol{t}
=
\lambda\boldsymbol{y}
+(1-\lambda)\boldsymbol{q},
\lambda\in[0,1)\}
\subseteq\operatorname{int}(\mathcal{C}')$.
This contradicts the fact that $\mathcal{C}'$ is DD-free
because $\boldsymbol{z}\in\operatorname{Sol}(\dd{D})$.
\Halmos
\endproof

Despite the facet-contact characterization established above, the following proposition shows that deciding maximality is computationally difficult even when DD-freeness is guaranteed.

\begin{proposition}
\label{prop:DD-free-maximality-hard}
    Given an IDD $\dd{D}$ with rational interval endpoints and a full-dimensional polyhedron $\mathcal{C}$ described by rational linear inequalities, deciding whether $\mathcal{C}$ is maximal is NP-complete, even under the promise that $\mathcal{C}$ is DD-free with respect to $\dd{D}$.
    This conclusion holds even when $\dd{D}$ has width one, every arc label is either $[0,0]$ or $[1,1]$, and $\mathcal{C}$ has only two facets.
\end{proposition}

\proof{Proof.}

We first show that the problem belongs to NP under the stated promise.
Write
$\mathcal{C}
=
\{
    \boldsymbol{x}\in\mathbb{R}^n:
    \boldsymbol{A}\boldsymbol{x}
    \le
    \boldsymbol{\beta}
\}$.
We may assume that this description is irredundant, since redundant inequalities can be deleted in polynomial time by LP.
Let $\boldsymbol{A}$ have $m$ rows.
For each $i\in[m]$, let 
$\boldsymbol{a}_i^\top$ denote its $i$-th row and define
$F_i
\coloneqq
\{
    \boldsymbol{x}\in\mathcal{C}:
    \boldsymbol{a}_i^\top\boldsymbol{x}
    =\beta_i
\}$.
By Proposition~\ref{prop:polyhedral-maximal-characterization},
$\mathcal{C}$ is maximal if and only if, for every $i\in[m]$, there exists an 
$\dd{r}$--$\dd{t}$ path $\dd{p}_i$ such that
$\mathcal{S}_{\dd{p}_i}
\cap
\operatorname{relint}(F_i)
\neq
\emptyset$.
For a given path $\dd{p}_i$, the condition
$\mathcal{S}_{\dd{p}_i}
\cap
\operatorname{relint}(F_i)
\neq
\emptyset$
holds if and only if the linear program
\[
    \max_{\boldsymbol{x},\varepsilon}
    \left\{
        \varepsilon:
        \boldsymbol{x}
        \in
        \mathcal{S}_{\dd{p}_i},\
        \boldsymbol{a}_i^\top\boldsymbol{x}
        =\beta_i,\
        \boldsymbol{a}_j^\top\boldsymbol{x}
        +\varepsilon
        \le\beta_j
        \quad
        \forall j\in[m]\setminus\{i\},\
        0\le\varepsilon\le1
    \right\}
\]
has a positive optimal value.
Since each path $\dd{p}_i$ consists of $n$ arcs, there are $m$ facets in $\mathcal{C}$, and each linear program has rational data and encoding length polynomial in the input size, the paths $\dd{p}_1,\ldots,\dd{p}_m$ form a polynomial-size certificate that can be verified in polynomial time through solving the above LP.
Hence the problem belongs to NP under the stated promise.

We now establish NP-completeness by a reduction from
\textsc{Partition} \cite{karp1972reducibility}.
Consider an instance of \textsc{Partition} specified by positive integers $a_1,\ldots,a_n$, where the goal is to determine whether they can be divided into two subsets of equal sum.
Construct an IDD with $n$ arc layers, one node in each node layer, and two parallel arcs, labeled $[0,0]$ and $[1,1]$, between each pair of consecutive node layers.
This IDD has $n+1$ nodes, $2n$ arcs, and width one, with
$\operatorname{Sol}(\dd{D})=\{0,1\}^n$.
For $\boldsymbol{x}\in\mathbb{R}^n$, write
$q(\boldsymbol{x})
\coloneqq
\sum_{i=1}^n a_i(2x_i-1)$,
and define
\[
    \mathcal{C}
    \coloneqq
    \left\{
        \boldsymbol{x}\in\mathbb{R}^n:
        0\le 2q(\boldsymbol{x})+x_1\le1
    \right\}.
\]
The set $\mathcal{C}$ is a polyhedron described by rational linear inequalities and is full-dimensional, since
$(1/2,\ldots,1/2)
\in
\operatorname{int}(\mathcal{C})$.
Its two facets are the hyperplanes
$F_k
\coloneqq
\{\boldsymbol{x}\in\mathbb{R}^n:
2q(\boldsymbol{x})+x_1=k\}$
for $k\in\{0,1\}$.
In particular, $\operatorname{relint}(F_k)=F_k$.
For every 
$\boldsymbol{x}
\in
\operatorname{Sol}(\dd{D})$,
the expression $2q(\boldsymbol{x})+x_1$ is an integer.
Thus, no such point belongs to 
$\operatorname{int}(\mathcal{C})$,
and $\mathcal{C}$ is DD-free for every \textsc{Partition} instance.
Moreover, for every
$\boldsymbol{x}
\in
\operatorname{Sol}(\dd{D})$,
we have 
$q(\boldsymbol{x})\in\mathbb{Z}$ and 
$x_1\in\{0,1\}$.
Thus, for each $k\in\{0,1\}$,
$\boldsymbol{x}\in F_k$ if and only if
$q(\boldsymbol{x})=0$ and $x_1=k$.
The condition 
$q(\boldsymbol{x})=0$ is equivalent to
$\sum_{i:x_i=1}a_i=\sum_{i:x_i=0}a_i$.
If the \textsc{Partition} instance has a solution, choose
$\boldsymbol{x}\in\{0,1\}^n$ with 
$q(\boldsymbol{x})=0$.
Then 
$\boldsymbol{1}-\boldsymbol{x}
\in
\{0,1\}^n$ and
$q(\boldsymbol{1}-\boldsymbol{x})
=-q(\boldsymbol{x})
=0$.
Since these two points have different first coordinates,
one belongs to $F_0$ and the other to $F_1$.
Conversely, if $\mathcal{C}$ is maximal, then
Proposition~\ref{prop:polyhedral-maximal-characterization}
implies that 
$\operatorname{Sol}(\dd{D})
\cap 
F_0
\neq
\emptyset$.
Any point 
$\boldsymbol{x}
\in
\operatorname{Sol}(\dd{D})
\cap 
F_0$
satisfies
$q(\boldsymbol{x})=0$ and hence yields a solution to
the \textsc{Partition} instance.
Consequently, $\mathcal{C}$ is maximal if and only if the
\textsc{Partition} instance has a solution.
This construction can be carried out in polynomial time and always satisfies the DD-freeness promise.
Therefore, deciding whether a given full-dimensional polyhedron $\mathcal{C}$ is maximal with respect to an IDD $\dd{D}$ is NP-complete, even under the promise that $\mathcal{C}$ is DD-free with respect to $\dd{D}$.
\Halmos
\endproof

It follows from Proposition~\ref{prop:polyhedral-maximal-characterization} that a key step in determining whether a given polyhedron is maximal is verifying the facet-contact condition for each of its facets.
The reduction in 
Proposition~\ref{prop:DD-free-maximality-hard}, however,
implies that such a verification task is itself computationally difficult for a general polyhedron, as presented next.

\begin{corollary}
\label{cor:DD-free-facet-contact-hard}
    Given an IDD $\dd{D}$ with rational interval endpoints and a full-dimensional DD-free polyhedron $\mathcal{C}$ described by rational linear inequalities, and a specified facet $F$ of $\mathcal{C}$, deciding whether 
    $\operatorname{relint}(F)
    \cap
    \operatorname{Sol}(\dd{D})
    \neq
    \emptyset$ is NP-complete.
    Moreover, deciding whether
    $F
    \cap
    \operatorname{Sol}(\dd{D})
    \neq
    \emptyset$ is NP-complete.
\end{corollary}

\proof{Proof.}

We first show that both problems belong to NP.
By the definition of $\operatorname{Sol}(\dd{D})$, a certificate for either problem is an $\dd{r}$--$\dd{t}$ path $\dd{p}$ whose induced hyperrectangle $\mathcal{S}_{\dd{p}}$ satisfies the corresponding contact condition.
We may assume that this description of $\mathcal{C}$ is irredundant, since redundant inequalities can be deleted in polynomial time by LP while retaining an inequality defining $F$.
Then, the condition
$\mathcal{S}_{\dd{p}}
\cap
\operatorname{relint}(F)
\neq
\emptyset$
can be verified by solving the LP used to establish membership in NP in the proof of
Proposition~\ref{prop:DD-free-maximality-hard}.
The condition
$\mathcal{S}_{\dd{p}}
\cap F\neq\emptyset$
can be verified by solving a linear feasibility problem.
Since each path consists of $n$ arcs and the resulting LPs have rational data and encoding length polynomial in the input size,
these certificates can be verified in polynomial time.
Hence both problems belong to NP.

To establish NP-hardness, consider the reduction from \textsc{Partition} in the proof of
Proposition~\ref{prop:DD-free-maximality-hard},
and let $F=F_0$.
Recall that
$q(\boldsymbol{x})=\sum_{i=1}^n a_i(2x_i-1)$
and
$F_0=\{
    \boldsymbol{x}\in\mathbb{R}^n:
    2q(\boldsymbol{x})+x_1=0
\}$.
Since $F_0$ is a hyperplane, we have
$\operatorname{relint}(F_0)=F_0$.
Moreover,
$\operatorname{Sol}(\dd{D})=\{0,1\}^n$
and $q(\boldsymbol{x})\in\mathbb{Z}$ for every
$\boldsymbol{x}\in\{0,1\}^n$, so
\[
    \operatorname{relint}(F_0)
    \cap
    \operatorname{Sol}(\dd{D})
    =
    F_0
    \cap
    \operatorname{Sol}(\dd{D})
    =
    \left\{
        \boldsymbol{x}\in\{0,1\}^n:
        q(\boldsymbol{x})=0,\ 
        x_1=0
    \right\}.
\]
For $\boldsymbol{x}\in\{0,1\}^n$,
$q(\boldsymbol{x})=0$ is equivalent to
$\sum_{i:x_i=1}a_i=\sum_{i:x_i=0}a_i$.
If the \textsc{Partition} instance has a solution, choose
$\boldsymbol{x}\in\{0,1\}^n$ satisfying
$q(\boldsymbol{x})=0$.
Since
$q(\boldsymbol{1}-\boldsymbol{x})
=-q(\boldsymbol{x})=0$,
we may replace $\boldsymbol{x}$ by its complement, if necessary, so that $x_1=0$.
Hence,
$\boldsymbol{x}
\in
\operatorname{relint}(F_0)
\cap
\operatorname{Sol}(\dd{D})$.
Conversely, any
$\boldsymbol{x}\in F_0\cap\operatorname{Sol}(\dd{D})$
satisfies $q(\boldsymbol{x})=0$ and hence yields a solution to the \textsc{Partition} instance.
Thus, each of the two contact conditions holds if and only if the \textsc{Partition} instance has a solution.
The reduction can be carried out in polynomial time and satisfies the DD-freeness promise for every instance.
Therefore, both problems are NP-complete.
\Halmos
\endproof

Although deciding whether a specified facet of a general DD-free polyhedron contains a contact point is NP-complete, 
boundary contacts of DD-free $\ell_p$-balls admit an efficiently computable implicit representation.
Building on Proposition~\ref{prop:distance-to-path},
the following result presents a polynomial-time procedure for identifying, without enumerating its 
$\dd{r}$--$\dd{t}$ paths, a subgraph of the IDD whose paths induce exactly the hyperrectangles that contain a contact point.

\begin{proposition}
\label{prop:ball-contact-subgraph}
    Let 
    $\dd{D}=(\dd{U},\dd{A},\dd{l}(\cdot))$ be an IDD with 
    $\operatorname{Sol}(\dd{D})
    \neq\emptyset$.
    Fix $p\in[1,\infty]$, and let
    $\mathcal{B}_p(\boldsymbol{c},r)$ be a DD-free ball with $r>0$.
    A subgraph $\dd{D}^{\partial}$ of $\dd{D}$ containing $\dd{r}$ and $\dd{t}$ can be constructed in
    $\mathcal{O}(|\dd{A}|)$ time such that, for every
    $\dd{r}$--$\dd{t}$ path $\dd{q}$ of $\dd{D}$,
    the path $\dd{q}$ is retained in $\dd{D}^{\partial}$ if and only if its induced hyperrectangle
    $\mathcal{S}_{\dd{q}}$ intersects the boundary of $\mathcal{B}_p(\boldsymbol{c},r)$.
    For each retained path, a contact point can be recovered in $\mathcal{O}(n)$ time by projecting $\boldsymbol{c}$ coordinatewise onto its induced hyperrectangle.
    If $p<\infty$, this projection is the unique contact point in the induced hyperrectangle, and every contact point of the ball is obtained from a retained path in this way.
\end{proposition}

\proof{Proof.}

For every $\dd{r}$--$\dd{t}$ path $\dd{q}$,
DD-freeness implies that
$d_p(
    \boldsymbol{c},
    \mathcal{S}_{\dd{q}}
)\ge r$.
Since $\mathcal{S}_{\dd{q}}$ is compact,
it intersects
$\operatorname{bd}\bigl(
    \mathcal{B}_p(\boldsymbol{c},r)
\bigr)$
if and only if
$d_p(
    \boldsymbol{c},
    \mathcal{S}_{\dd{q}}
)=r$.
Using Algorithm~\ref{alg:lp-ball-DD-free}, we compute
$d_p(
    \boldsymbol{c},
    \operatorname{Sol}(\dd{D})
)$.
Since the ball is DD-free, this distance is at least $r$.
If it exceeds $r$, no path-induced hyperrectangle contains a contact point, so we take $\dd{D}^{\partial}$ to consist of $\dd{r}$ and $\dd{t}$ with no arcs.
Otherwise, the distance equals $r$ and the subgraph can be constructed as follows.

For $p<\infty$, assign each arc $\dd{a}\in\dd{A}$ the weight
$\dd{w}(\dd{a})
\coloneqq
\xi(\dd{a};\boldsymbol{c})^p$, where $\xi$ is the arc deviation defined in \eqref{eq:arc-deviation}.
By Proposition~\ref{prop:distance-to-path},
the contact paths are precisely the shortest $\dd{r}$--$\dd{t}$ paths under these weights.
Since the path costs are additive, we augment the forward pass of Algorithm~\ref{alg:lp-ball-DD-free} to record all minimizing predecessor arcs at each node.
Traversing these arcs backward from $\dd{t}$ identifies all arcs that belong to at least one shortest $\dd{r}$--$\dd{t}$ path.
We take $\dd{D}^{\partial}$ to consist of these arcs and their incident nodes.

For $p=\infty$, we first construct a subgraph consisting of $\dd{r}$, $\dd{t}$, and all arcs
$\dd{a}\in\dd{A}$ satisfying
$\xi(\dd{a};\boldsymbol{c})\le r$,
together with their incident nodes.
We then obtain $\dd{D}^{\partial}$ by removing all nodes and arcs that do not lie on any $\dd{r}$--$\dd{t}$ path of this subgraph, using forward and backward reachability passes from $\dd{r}$ and $\dd{t}$, respectively.
An $\dd{r}$--$\dd{t}$ path $\dd{q}$ is retained if and only if its maximum arc deviation is at most $r$.
By the definition of arc deviation and the fact that $\mathcal{S}_{\dd{q}}$ is a Cartesian product of intervals, this condition is equivalent to
$d_\infty(
    \boldsymbol{c},
    \mathcal{S}_{\dd{q}}
)\le r$.
Since DD-freeness implies
$d_\infty(
    \boldsymbol{c},
    \mathcal{S}_{\dd{q}}
)\ge r$,
the retained $\dd{r}$--$\dd{t}$ paths are precisely those whose induced hyperrectangles intersect
$\operatorname{bd}\bigl(
    \mathcal{B}_\infty(\boldsymbol{c},r)
\bigr)$.

The initial distance computation requires
$\mathcal{O}(|\dd{A}|)$ time by
Proposition~\ref{prop:lp-ball-DD-free-time-complexity}.
In either case, the remaining construction requires a constant number of graph traversals and hence takes
$\mathcal{O}(|\dd{U}|+|\dd{A}|)$ time.
By \eqref{eq:DD-size-bounds},
we have $|\dd{U}|\le 2|\dd{A}|$.
Therefore, the total running time is
$\mathcal{O}(|\dd{A}|)$.

Finally, consider a retained path
$\dd{q}=(\dd{a}_1,\dots,\dd{a}_n)$,
and let $\boldsymbol{z}$ be the coordinatewise projection of $\boldsymbol{c}$ onto $\mathcal{S}_{\dd{q}}$.
Since $\boldsymbol{z}$ is the coordinatewise projection,
\eqref{eq:arc-deviation} gives
$|c_i-z_i|=\xi(\dd{a}_i;\boldsymbol{c})$
for all $i\in[n]$.
Because each $z_i$ is the closest point to $c_i$ in $\dd{l}(\dd{a}_i)$, $\boldsymbol{z}$ minimizes the $\ell_p$-distance from $\boldsymbol{c}$ over $\mathcal{S}_{\dd{q}}$.
Hence,
$\lVert
    \boldsymbol{c}-\boldsymbol{z}
\rVert_p
=d_p(
    \boldsymbol{c},
    \mathcal{S}_{\dd{q}}
)=r$,
where the last equality follows from the retention
of $\dd{q}$.
Therefore, $\boldsymbol{z}$ is a contact point of $\mathcal{B}_p(\boldsymbol{c},r)$.
Each coordinate of $\boldsymbol{z}$ is obtained by projecting a scalar onto an interval, which takes constant time, giving
$\mathcal{O}(n)$ time in total.
If $p<\infty$, consider any
$\boldsymbol{y}
\in
\mathcal{S}_{\dd{q}}
\setminus
\{\boldsymbol{z}\}$.
By coordinatewise projection,
$|c_i-y_i|\ge|c_i-z_i|$ for all $i\in[n]$.
Since projection onto a closed interval is unique
and 
$\boldsymbol{y}
\neq
\boldsymbol{z}$,
at least one of these inequalities is strict.
Hence,
$\lVert
    \boldsymbol{c}-\boldsymbol{y}
\rVert_p>r$,
which shows that $\boldsymbol{z}$ is the unique contact point in $\mathcal{S}_{\dd{q}}$.
Moreover, every contact point of the ball belongs to a hyperrectangle induced by some $\dd{r}$--$\dd{t}$ path retained in $\dd{D}^{\partial}$.
For $p<\infty$, the uniqueness established above implies that projecting $\boldsymbol{c}$ coordinatewise onto the hyperrectangles induced by the retained paths recovers all contact points of the ball.
\Halmos
\endproof

Proposition~\ref{prop:ball-contact-subgraph}
provides an efficient representation of boundary contact points for DD-free $\ell_p$-balls through a compact subgraph structure. This representation facilitates the characterization of maximality among DD-free sets. For example, unlike the general case in Corollary~\ref{cor:DD-free-facet-contact-hard}, which shows that determining whether a facet of a DD-free polyhedron contains a contact point is computationally difficult, this verification can be performed in time linear in the size of the IDD when the polyhedron is an $\ell_p$-ball with $p\in\{1,\infty\}$.

\begin{corollary}
\label{cor:l1-linf-ball-facet-contact}
    Let $\dd{D}=(\dd{U},\dd{A},\dd{l}(\cdot))$ be an IDD with 
    $\operatorname{Sol}(\dd{D})
    \neq\emptyset$.
    Fix $p\in\{1,\infty\}$, and let
    $\mathcal{B}_p(\boldsymbol{c},r)$ be a DD-free ball with $r>0$.
    Given a facet $F$ of this ball, it can be determined in $\mathcal{O}(|\dd{A}|)$ time whether $F$ contains a contact point and whether its relative interior contains one.
    Whenever either condition holds, a corresponding $\dd{r}$--$\dd{t}$ path can be identified within the same time bound, from which a contact point satisfying that condition can be recovered in $\mathcal{O}(n)$ time.
\end{corollary}

\proof{Proof.}

Construct the subgraph $\dd{D}^{\partial}$ from 
Proposition~\ref{prop:ball-contact-subgraph}.
For each $i\in[n]$ and $\dd{a}\in\dd{A}_i$, let $z(\dd{a})$ denote the projection of $c_i$ onto $\dd{l}(\dd{a})$.
For each norm, we characterize $F$ and
$\operatorname{relint}(F)$ and then reduce the corresponding contact conditions to
$\dd{r}$--$\dd{t}$ reachability tests in subgraphs
of $\dd{D}^{\partial}$.

Suppose first that $p=1$.
There exists $\boldsymbol{\sigma}\in\{-1,1\}^n$ such that
\[
    F=\left\{
        \boldsymbol{x}
        \in
        \mathcal{B}_1(\boldsymbol{c},r):
        \boldsymbol{\sigma}^{\top}
        (\boldsymbol{x}-\boldsymbol{c})=r
    \right\}.
\]
We first show that a point 
$\boldsymbol{x}\in\mathbb{R}^n$
belongs to $F$ if and only if
$\boldsymbol{\sigma}^{\top}(
    \boldsymbol{x}-\boldsymbol{c}
)=r$
and $\sigma_i(x_i-c_i)\ge0$ for all $i\in[n]$.
For every $\boldsymbol{x}
\in
\mathcal{B}_1(\boldsymbol{c},r)$,
\begin{equation}
\label{eq:l1-facet}
    \boldsymbol{\sigma}^{\top}(
        \boldsymbol{x}-\boldsymbol{c}
    )
    =
    \sum_{i=1}^n \sigma_i(x_i-c_i)
    \le
    \sum_{i=1}^n |x_i-c_i|
    \le r.
\end{equation}
If $\boldsymbol{x}\in F$, then
$\boldsymbol{\sigma}^{\top}(
    \boldsymbol{x}-\boldsymbol{c}
)=r$,
and hence both inequalities in \eqref{eq:l1-facet} hold with equality.
In particular, the first inequality in \eqref{eq:l1-facet} holds at equality only if 
$\sigma_i(x_i-c_i)=|x_i-c_i|$ for all $i\in[n]$,
or equivalently,
$\sigma_i(x_i-c_i)\ge0$ for all $i\in[n]$.
Conversely, suppose that
$\boldsymbol{\sigma}^{\top}(
    \boldsymbol{x}-\boldsymbol{c}
)=r$
and $\sigma_i(x_i-c_i)\ge0$ for all $i\in[n]$.
Since $\sigma_i\in\{-1,1\}$, these inequalities give
$|x_i-c_i|=\sigma_i(x_i-c_i)$ for all $i\in[n]$.
Hence, $\lVert\boldsymbol{x}-\boldsymbol{c}\rVert_1=r$,
which shows that $\boldsymbol{x}\in F$.
We next characterize the relative interior of $F$.
Since $r>0$, the point
$\boldsymbol{y}=
\boldsymbol{c}+(r/n)\boldsymbol{\sigma}$
belongs to $F$ and satisfies
$\sigma_i(y_i-c_i)=r/n>0$ for all $i\in[n]$.
Hence, a point 
$\boldsymbol{x}\in\mathbb{R}^n$
belongs to $\operatorname{relint}(F)$ if and only if
$\boldsymbol{\sigma}^{\top}(
    \boldsymbol{x}-\boldsymbol{c}
)=r$
and $\sigma_i(x_i-c_i)>0$ for all $i\in[n]$.

Let $\dd{D}^{F}$ be obtained from $\dd{D}^{\partial}$ by retaining, for each $i\in[n]$, only those arcs
$\dd{a}\in\dd{A}_i$ satisfying
$\sigma_i(z(\dd{a})-c_i)\ge0$.
Similarly, let $\dd{D}^{\operatorname{relint}(F)}$ be obtained by replacing the non-strict inequality above with a strict inequality.
By Proposition~\ref{prop:ball-contact-subgraph},
each $\dd{r}$--$\dd{t}$ path
$\dd{q}=(\dd{a}_1,\ldots,\dd{a}_n)$
of $\dd{D}^{\partial}$ induces a hyperrectangle whose unique contact point is
$\boldsymbol{z}_{\dd{q}}
\coloneqq
\bigl(z(\dd{a}_1),\ldots,z(\dd{a}_n)\bigr)$.
The facet descriptions above imply that
$\boldsymbol{z}_{\dd{q}}\in F$ if and only if $\dd{q}$ is retained in $\dd{D}^{F}$, and that
$\boldsymbol{z}_{\dd{q}}
\in
\operatorname{relint}(F)$
if and only if $\dd{q}$ is retained in
$\dd{D}^{\operatorname{relint}(F)}$.
Since every contact point of the $\ell_1$-ball equals
$\boldsymbol{z}_{\dd{q}}$ for some
$\dd{r}$--$\dd{t}$ path $\dd{q}$ of
$\dd{D}^{\partial}$, the equivalences above imply
that 
$F\cap
\operatorname{Sol}(\dd{D})\neq
\emptyset$
if and only if $\dd{D}^{F}$ contains an
$\dd{r}$--$\dd{t}$ path.
The same holds for $\operatorname{relint}(F)$
and $\dd{D}^{\operatorname{relint}(F)}$.

Now suppose that $p=\infty$.
Since
$\mathcal{B}_\infty(\boldsymbol{c},r)
=\prod_{i\in[n]}[c_i-r,c_i+r]$
is a full-dimensional hyperrectangle, each of its facets is obtained by fixing one coordinate at an endpoint of its corresponding interval.
Thus, there exist $j\in[n]$ and
$\sigma\in\{-1,1\}$ such that
\[
    F=
    \left\{
        \boldsymbol{x}\in\mathbb{R}^n:
        x_j=c_j+\sigma r,\ 
        |x_i-c_i|\le r
        \quad\forall i\in[n]\setminus\{j\}
    \right\}.
\]
In $F$, the $j$-th coordinate is fixed at
$c_j+\sigma r$, while each remaining coordinate
ranges over $[c_i-r,c_i+r]$.
Since $r>0$, $\operatorname{relint}(F)$ is obtained
by replacing $|x_i-c_i|\le r$ with
$|x_i-c_i|<r$ for every $i\neq j$.

Let $\dd{D}^{F}$ be obtained from $\dd{D}^{\partial}$
by retaining in arc layer $j$ only those arcs
$\dd{a}\in\dd{A}_j$ satisfying
$c_j+\sigma r\in\dd{l}(\dd{a})$,
and leaving all other arc layers unchanged.
Let $\dd{D}^{\operatorname{relint}(F)}$ be obtained from
$\dd{D}^{F}$ by additionally retaining, for each
$i\neq j$, only those arcs $\dd{a}\in\dd{A}_i$ satisfying
$|z(\dd{a})-c_i|<r$.
We next show that the two subgraphs $\dd{D}^{F}$ and 
$\dd{D}^{\operatorname{relint}(F)}$ retain precisely the paths whose induced hyperrectangles intersect
$F$ and $\operatorname{relint}(F)$, respectively.
Consider an arbitrary $\dd{r}$--$\dd{t}$ path
$\dd{q}=(\dd{a}_1,\ldots,\dd{a}_n)$
of $\dd{D}^{\partial}$.
By Proposition~\ref{prop:ball-contact-subgraph},
the coordinatewise projection
$\boldsymbol{z}_{\dd{q}}
\coloneqq
\bigl(z(\dd{a}_1),\ldots,z(\dd{a}_n)\bigr)$
is a contact point of $\mathcal{B}_\infty(\boldsymbol{c},r)$.
Hence, $\max_{i\in[n]}|z(\dd{a}_i)-c_i|=r$,
which implies that
$|z(\dd{a}_i)-c_i|\le r$
for every $i\in[n]$.
Hence,
$\dd{l}(\dd{a}_i)
\cap
[c_i-r,c_i+r]
\neq\emptyset$
for every $i\in[n]$.
Moreover, for every 
$i\in[n]\setminus\{j\}$,
the interval $\dd{l}(\dd{a}_i)$ intersects
$(c_i-r,c_i+r)$ if and only if
$|z(\dd{a}_i)-c_i|<r$, since $z(\dd{a}_i)$
is the projection of $c_i$ onto
$\dd{l}(\dd{a}_i)$.
Since $\mathcal{S}_{\dd{q}}$ is a Cartesian product
of the arc labels, it follows that
$\mathcal{S}_{\dd{q}}
\cap 
F
\neq\emptyset$
if and only if $\dd{q}$ is retained in $\dd{D}^{F}$,
and that
$\mathcal{S}_{\dd{q}}
\cap
\operatorname{relint}(F)
\neq\emptyset$
if and only if $\dd{q}$ is retained in
$\dd{D}^{\operatorname{relint}(F)}$.
Since every contact point of the $\ell_\infty$-ball
belongs to $\mathcal{S}_{\dd{q}}$ for some
$\dd{r}$--$\dd{t}$ path $\dd{q}$ of
$\dd{D}^{\partial}$, the equivalences above imply
that
$F
\cap
\operatorname{Sol}(\dd{D})
\neq\emptyset$
if and only if $\dd{D}^{F}$ contains an
$\dd{r}$--$\dd{t}$ path.
The same holds for $\operatorname{relint}(F)$
and $\dd{D}^{\operatorname{relint}(F)}$.
For any $\dd{r}$--$\dd{t}$ path
$\dd{q}=(\dd{a}_1,\ldots,\dd{a}_n)$
of $\dd{D}^{F}$, set
$x_j=c_j+\sigma r$ and
$x_i=z(\dd{a}_i)$ for all $i\neq j$.
Then 
$\boldsymbol{x}
\in
\mathcal{S}_{\dd{q}}\cap F$,
so $\boldsymbol{x}$ is a contact point.
If $\dd{q}$ is a path of
$\dd{D}^{\operatorname{relint}(F)}$,
the same construction yields
$\boldsymbol{x}\in\operatorname{relint}(F)$.

In either case, the subgraph $\dd{D}^{\partial}$ can be constructed in $\mathcal{O}(|\dd{A}|)$ time by
Proposition~\ref{prop:ball-contact-subgraph}.
For each arc $\dd{a}$ of $\dd{D}^{\partial}$,
computing $z(\dd{a})$ and determining whether
$\dd{a}$ is retained in $\dd{D}^{F}$ and
$\dd{D}^{\operatorname{relint}(F)}$ take constant time.
Thus, both $\dd{D}^{F}$ and
$\dd{D}^{\operatorname{relint}(F)}$ can be constructed in
$\mathcal{O}(|\dd{A}|)$ time.
Testing $\dd{r}$--$\dd{t}$ reachability in either subgraph
takes $\mathcal{O}(|\dd{U}|+|\dd{A}|)$ time.
By retaining predecessor arcs during the reachability tests,
an $\dd{r}$--$\dd{t}$ path can be recovered in
$\mathcal{O}(n)$ time whenever one exists.
Given such a path, a contact point can be recovered in
$\mathcal{O}(n)$ time by the constructions above.
Since $|\dd{U}|\le 2|\dd{A}|$ and $n\le|\dd{A}|$ by
\eqref{eq:DD-size-bounds}, the total running time is
$\mathcal{O}(|\dd{A}|)$.
\Halmos
\endproof

The results of Proposition~\ref{prop:ball-contact-subgraph} and Corollary~\ref{cor:l1-linf-ball-facet-contact} show that verifying maximality of DD-free sets can be facilitated by restricting attention to structured families such as $\ell_p$-balls. This motivates the study of maximality relative to a prescribed family $\mathfrak{F}$, where a DD-free set is considered maximal with respect to other DD-free sets belonging to $\mathfrak{F}$.

\begin{definition}[$\mathfrak{F}$-maximal DD-free set]
\label{def:maximal-DD-free}
    Given an IDD $\dd{D}$ and a family $\mathfrak{F}$ of convex sets
    of $\mathbb{R}^n$, a set $\mathcal{C}\in\mathfrak{F}$ that is
    DD-free with respect to $\dd{D}$ is
    \emph{$\mathfrak{F}$-maximal} if there is no DD-free set
    $\mathcal{C}'\in\mathfrak{F}$ such that
    $\mathcal{C}\subsetneq\mathcal{C}'$.
    When $\mathfrak{F}$ is the family of all convex subsets of $\mathbb{R}^n$, this definition reduces to the general notion of maximality for $\mathcal{C}$ introduced earlier.
\end{definition}

Motivated by the above results for $\ell_p$-balls, we focus on two structured families. The first consists of boxes, which generalize $\ell_{\infty}$-balls, and is denoted by $\mathfrak{F}^{\square}$.
Because $\operatorname{Sol}(\dd{D})$ is a finite union of compact hyperrectangles, boxes are a natural choice for identifying maximality of DD-free sets. 
Here, a box is a Cartesian product of nonempty closed intervals, which are allowed to be unbounded.
A DD-free box that is maximal within $\mathfrak{F}^{\square}$ is called \emph{box-maximal}.
Section~\ref{subsec:box-maximal-DD-free} shows that, for full-dimensional boxes, box-maximality coincides with maximality among all DD-free sets and develops an efficient sequential facet expansion algorithm for constructing box-maximal DD-free sets.

The second family consists of $\ell_p$-balls, for $p \in [1, \infty)$, and is denoted by $\mathfrak{F}^{\ell_p}$.
For a fixed $p\in[1,\infty)$, a DD-free $\ell_p$-ball that is maximal within $\mathfrak{F}^{\ell_p}$ is called \emph{$\ell_p$-ball-maximal}.
The $\ell_p$-ball $\mathcal{C}_p(\bar{\boldsymbol{x}})$ is centered at $\bar{\boldsymbol{x}}$, whereas an $\ell_p$-ball-maximal DD-free ball containing $\bar{\boldsymbol{x}}$ in its interior need not be centered at $\bar{\boldsymbol{x}}$.
Section~\ref{subsec:lp-ball-maximal-DD-free} characterizes $\ell_p$-ball-maximality for every fixed $p\in[1,\infty)$ and develops a radial expansion algorithm for enlarging $\mathcal{C}_p(\bar{\boldsymbol{x}})$.

\subsection{Box-Maximal DD-Free Sets}
\label{subsec:box-maximal-DD-free}

We restrict attention to full-dimensional boxes, since every lower-dimensional box has empty interior in $\mathbb{R}^n$ and is therefore automatically DD-free.
We first characterize DD-freeness through reachability in the IDD and box-maximality through a facet-contact condition.
We then use these characterizations to develop a sequential facet expansion procedure initialized at
$\mathcal{C}_{\infty}(\bar{\boldsymbol{x}})$ 
and show that it constructs a box-maximal DD-free set in time linear in the size of the IDD.
Finally, we derive the intersection cut associated with the expanded box, compare its strength with the initial $\ell_{\infty}$-ball cut, and establish the corresponding cut-closure result.

\subsubsection{Reachability and Box-Maximality}
\label{subsubsec:box-maximality-characterization}

Recall that $\operatorname{Sol}(\dd{D})$ is the union of the hyperrectangles induced by 
all $\dd{r}$--$\dd{t}$ paths of $\dd{D}$.
Testing each path-induced hyperrectangle for intersection with the interior of a box may require examining exponentially many paths.
The following reachability characterization avoids this enumeration.

Let 
$\mathcal{R}
\coloneqq
\prod_{i\in[n]}[l_i,u_i]$ 
be a box, where $l_i\in\mathbb{R}\cup\{-\infty\}$, $u_i\in\mathbb{R}\cup\{+\infty\}$, and $l_i<u_i$ for every $i\in[n]$.
For each $i\in[n]$, define
\[
    \dd{A}_i^{\mathrm{act}}(\mathcal{R})
    \coloneqq
    \left\{
        \dd{a}\in\dd{A}_i:
        \dd{l}(\dd{a})\cap(l_i,u_i)\neq\emptyset
    \right\}
    =
    \left\{
        \dd{a}\in\dd{A}_i:
        \dd{l}(\dd{a})\downarrow \,< u_i,\;
        \dd{l}(\dd{a})\uparrow \, > l_i
    \right\}.
\]
We call these arcs \emph{active with respect to $\mathcal{R}$}.
The open interval $(l_i,u_i)$ reflects the use of 
$\operatorname{int}(\mathcal{R})$ in the definition of DD-free sets; in particular, an arc whose label only touches an endpoint of $[l_i,u_i]$ is not active.
Let
$\dd{A}^{\mathrm{act}}(\mathcal{R})
\coloneqq
\bigcup_{i\in[n]}\dd{A}_i^{\mathrm{act}}(\mathcal{R})$,
and refer to the subgraph obtained from $\dd{D}$ by retaining only the arcs in
$\dd{A}^{\mathrm{act}}(\mathcal{R})$ as the \emph{active subgraph} induced by $\mathcal{R}$. It is clear that this subgraph can be computed in $\mathcal{O}(|\dd{A}|)$ time.

\begin{lemma}
\label{lemma:box-DD-free-reachability}
    Let $\dd{D}=(\dd{U},\dd{A},\dd{l}(\cdot))$ be an IDD, and let 
    $\mathcal{R}
    \coloneqq
    \prod_{i\in[n]}[l_i,u_i]$ 
    be a full-dimensional box.
    Then $\mathcal{R}$ is DD-free with respect to $\dd{D}$ if and only if the active subgraph induced by $\mathcal{R}$ contains no $\dd{r}$--$\dd{t}$ path.
\end{lemma}

\proof{Proof.}

Fix an arbitrary $\dd{r}$--$\dd{t}$ path $\dd{p}=(\dd{a}_1,\ldots,\dd{a}_n)$.
Then 
$\mathcal{S}_{\dd{p}}
\cap
\operatorname{int}(\mathcal{R})
=
\prod_{i\in[n]}
\bigl(
    \dd{l}(\dd{a}_i)\cap(l_i,u_i)
\bigr)$.
This intersection is nonempty if and only if every arc of $\dd{p}$ is active with respect to $\mathcal{R}$. 
Since $\dd{p}$ was arbitrary and
$\operatorname{Sol}(\dd{D})
=
\bigcup_{\dd{p}\in\dd{P}(\dd{t})}
\mathcal{S}_{\dd{p}}$,
it follows that
$\operatorname{int}(\mathcal{R})\cap\operatorname{Sol}(\dd{D})
\neq\emptyset$
if and only if the active subgraph induced by $\mathcal{R}$ contains 
an $\dd{r}$--$\dd{t}$ path.
The result then follows from
Definition~\ref{def:DD-free-sets}.
\Halmos
\endproof

We next show that box-maximality is characterized by the requirement that every finite facet of $\mathcal{R}$ contains a contact point in its relative interior, thereby coinciding with the maximality condition for general polyhedra established in Proposition~\ref{prop:polyhedral-maximal-characterization}.
For each $i\in[n]$ with $l_i\in\mathbb{R}$, define the
\emph{finite lower facet} of $\mathcal{R}$ in coordinate $i$ by
$F_i^-(\mathcal{R})
\coloneqq
\left\{
    \boldsymbol{x}\in\mathcal{R}:
    x_i=l_i
\right\}$.
Similarly, for each $i\in[n]$ with $u_i\in\mathbb{R}$, define its
\emph{finite upper facet} by
$F_i^+(\mathcal{R})
\coloneqq
\left\{
    \boldsymbol{x}\in\mathcal{R}:
    x_i=u_i
\right\}$.

\begin{proposition}[Maximality characterization for DD-free boxes]
\label{prop:box-maximal-characterization}
    Let
    $\dd{D}=(\dd{U},\dd{A},\dd{l}(\cdot))$
    be an IDD, and let
    $\mathcal{R}\coloneqq\prod_{i\in[n]}[l_i,u_i]$
    be a full-dimensional DD-free box with respect to $\dd{D}$.
    The following statements are equivalent:
    \begin{enumerate}[label=\textup{(\roman*)}, ref=\textup{(\roman*)}]
        \item
        $\mathcal{R}$ is box-maximal;
        \item
        every finite facet $F$ of $\mathcal{R}$ satisfies
        $\operatorname{relint}(F)
        \cap
        \operatorname{Sol}(\dd{D})
        \neq
        \emptyset$;
        \item
        $\mathcal{R}$ is maximal among all DD-free sets.
    \end{enumerate}
\end{proposition}

\proof{Proof.}

We first prove that \textup{(i)} implies \textup{(ii)}.
Suppose that $\mathcal{R}$ is box-maximal, and fix an arbitrary 
finite lower facet $F_i^-(\mathcal{R})$.
Since $\dd{A}_i$ contains finitely many arcs, there exists $\varepsilon>0$ such that
$\dd{l}(\dd{a})\uparrow
\notin
(l_i-\varepsilon,l_i)$
for all $\dd{a}\in\dd{A}_i$.
Let $\widehat{\mathcal{R}}$ be obtained from $\mathcal{R}$ by 
replacing $l_i$ with $l_i-\varepsilon$.
Because
$\mathcal{R}\subsetneq\widehat{\mathcal{R}}$,
box-maximality implies that $\widehat{\mathcal{R}}$ is not DD-free.
By Lemma~\ref{lemma:box-DD-free-reachability}, there exists an
$\dd{r}$--$\dd{t}$ path
$\dd{p}=(\dd{a}_1,\ldots,\dd{a}_n)$
whose arcs are all active with respect to
$\widehat{\mathcal{R}}$.

For every $j\neq i$, the interval corresponding to coordinate $j$ of $\widehat{\mathcal{R}}$ is unchanged, so $\dd{a}_j$ is active with respect to $\mathcal{R}$.
Since $\mathcal{R}$ is DD-free, $\dd{a}_i$ cannot be active with respect to $\mathcal{R}$.
Because $\dd{a}_i$ is active with respect to $\widehat{\mathcal{R}}$ but inactive with respect to $\mathcal{R}$, and only the lower endpoint in coordinate $i$ of $\widehat{\mathcal{R}}$ has changed, we have
$l_i-\varepsilon
<
\dd{l}(\dd{a}_i)\uparrow
\le
l_i$.
By the choice of $\varepsilon$, this inequality can hold only if
$\dd{l}(\dd{a}_i)\uparrow=l_i$.
Set $z_i\coloneqq l_i$ and, for each $j\neq i$, choose
$z_j\in
\dd{l}(\dd{a}_j)\cap(l_j,u_j)$, which is nonempty as $\dd{a}_j$ is active with respect to $\mathcal{R}$.
Then, we can write
\[
    \boldsymbol{z}
    \in
    \mathcal{S}_{\dd{p}}
    \cap
    \operatorname{relint}\bigl(F_i^-(\mathcal{R})\bigr)
    \subseteq
    \operatorname{Sol}(\dd{D})
    \cap
    \operatorname{relint}\bigl(F_i^-(\mathcal{R})\bigr),
\]
where the first inclusion follows from the fact that $\boldsymbol{z} \in \prod_{j\in[n]}\dd{l}(\dd{a}_j) = \mathcal{S}_{\dd{p}}$ and the definition of $F_i^-(\mathcal{R})$.
Since the finite lower facet $F_i^-(\mathcal{R})$ was arbitrary, every finite lower facet satisfies the contact condition.
The argument for finite upper facets is symmetric, using the lower endpoints of the arc labels.
This proves \textup{(ii)}.

Because $\mathcal{R}$ is a full-dimensional DD-free polyhedron and every facet of $\mathcal{R}$ is a finite lower or upper facet,
\textup{(ii)} implies \textup{(iii)} by
Proposition~\ref{prop:polyhedral-maximal-characterization}.
Finally, \textup{(iii)} implies \textup{(i)}, since every DD-free box is a DD-free set.
\Halmos
\endproof

We next use the reachability characterization to develop a sequential facet expansion procedure that preserves DD-freeness and yields a box satisfying the facet-contact condition in
Proposition~\ref{prop:box-maximal-characterization}.
As shown below, this procedure also provides an efficient test for box-maximality by checking whether the input box is returned unchanged.

\subsubsection{Sequential Box Expansion}
\label{subsubsec:box-expansion}

Let
$r
\coloneqq
d_\infty\bigl(
    \bar{\boldsymbol{x}},
    \operatorname{Sol}(\dd{D})
\bigr)$,
and set
\[
    \mathcal{R}^0
    \coloneqq
    \mathcal{C}_\infty(\bar{\boldsymbol{x}})
    =
    \prod_{i\in[n]}
    [\bar{x}_i-r,\bar{x}_i+r].
\]
By Proposition~\ref{prop:lp-ball-DD-free}, 
$\mathcal{R}^0$ is DD-free and contains 
$\bar{\boldsymbol{x}}$ in its interior.

Consider a DD-free box
$\mathcal{R}
\coloneqq
\prod_{j\in[n]}[l_j,u_j]$.
We first determine the largest expansion of a single facet of
$\mathcal{R}$ that preserves DD-freeness.
For each $i\in[n+1]$, let
$\dd{U}_i^{\mathrm{fwd}}(\mathcal{R})$
denote the set of nodes in $\dd{U}_i$ that are reachable from $\dd{r}$ using only arcs active with respect to $\mathcal{R}$.
Similarly, let
$\dd{U}_i^{\mathrm{bwd}}(\mathcal{R})$
denote the set of nodes in $\dd{U}_i$ from which $\dd{t}$ is reachable using only arcs active with respect to $\mathcal{R}$.
It is easy to verify that both of these sets can be computed in $\mathcal{O}(|\dd{A}|)$ time.

For each $i\in[n]$, define the set of \emph{critical arcs} in layer $i$ by
\begin{equation}
\label{eq:critical-arcs}
    \dd{A}_i^{\mathrm{crit}}(\mathcal{R})
    \coloneqq
    \left\{
        \dd{a}\in\dd{A}_i:
        \dd{tail}(\dd{a})
        \in
        \dd{U}_i^{\mathrm{fwd}}(\mathcal{R}),
        \;
        \dd{head}(\dd{a})
        \in
        \dd{U}_{i+1}^{\mathrm{bwd}}(\mathcal{R})
    \right\}.
\end{equation}
The next lemma shows that, when only coordinate $i$ is expanded, the resulting box remains DD-free if and only if no critical arc in layer $i$ becomes active.

\begin{lemma}
\label{lemma:single-coordinate-expansion}
    Let
    $\dd{D}=(\dd{U},\dd{A},\dd{l}(\cdot))$
    be an IDD, and let
    $\mathcal{R}=\prod_{j\in[n]}[l_j,u_j]$
    be a full-dimensional DD-free box with respect to $\dd{D}$. 
    Fix $i\in[n]$.
    Let
    $\widehat{\mathcal{R}}
    =
    \prod_{j\in[n]}[\widehat{l}_j,\widehat{u}_j]$
    be a box obtained from $\mathcal{R}$ by enlarging only its 
    $i$-th coordinate interval, i.e., 
    $[l_i,u_i]\subseteq[\widehat{l}_i,\widehat{u}_i]$ 
    and 
    $[\widehat{l}_j,\widehat{u}_j]=[l_j,u_j]$
    for all $j\in[n]\setminus\{i\}$.
    Then $\widehat{\mathcal{R}}$ is DD-free if and only if
    \[
        \dd{A}_i^{\mathrm{crit}}(\mathcal{R})
        \cap
        \dd{A}_i^{\mathrm{act}}(\widehat{\mathcal{R}})
        =
        \emptyset.
    \]
\end{lemma}

\proof{Proof.}

Suppose first that
$\dd{A}_i^{\mathrm{crit}}(\mathcal{R})
\cap
\dd{A}_i^{\mathrm{act}}(\widehat{\mathcal{R}})
\neq
\emptyset$,
and let $\dd{a}$ belong to this intersection.
By the definition of $\dd{A}_i^{\mathrm{crit}}(\mathcal{R})$, there exists an $\dd{r}$--$\dd{tail}(\dd{a})$ path and an $\dd{head}(\dd{a})$--$\dd{t}$ path in the active subgraph induced by $\mathcal{R}$.
Since neither path contains an arc in layer $i$ and $\widehat{\mathcal{R}}$ differs from $\mathcal{R}$ only in its $i$-th coordinate interval, all arcs on these paths remain active with respect to $\widehat{\mathcal{R}}$.
Concatenating the $\dd{r}$--$\dd{tail}(\dd{a})$ path, the arc $\dd{a}$, and the $\dd{head}(\dd{a})$--$\dd{t}$ path yields an $\dd{r}$--$\dd{t}$ path in the active subgraph induced by $\widehat{\mathcal{R}}$.
Lemma~\ref{lemma:box-DD-free-reachability} therefore implies that $\widehat{\mathcal{R}}$ is not DD-free.

Conversely, suppose that $\widehat{\mathcal{R}}$ is not DD-free.
By Lemma~\ref{lemma:box-DD-free-reachability}, the active subgraph induced by $\widehat{\mathcal{R}}$ contains an $\dd{r}$--$\dd{t}$ path $\dd{p}=(\dd{a}_1,\ldots,\dd{a}_n)$.
For every $j\neq i$, the $j$-th coordinate interval is unchanged, so $\dd{a}_j$ is also active with respect to $\mathcal{R}$.
Hence the prefix $(\dd{a}_1,\ldots,\dd{a}_{i-1})$
shows that
$\dd{tail}(\dd{a}_i)
\in
\dd{U}_i^{\mathrm{fwd}}(\mathcal{R})$,
and the suffix $(\dd{a}_{i+1},\ldots,\dd{a}_n)$
shows that
$\dd{head}(\dd{a}_i)
\in
\dd{U}_{i+1}^{\mathrm{bwd}}(\mathcal{R})$.
Therefore, 
$\dd{a}_i\in\dd{A}_i^{\mathrm{crit}}(\mathcal{R})$.
Since $\dd{a}_i$ is active with respect to $\widehat{\mathcal{R}}$, we conclude that
$\dd{A}_i^{\mathrm{crit}}(\mathcal{R})
\cap
\dd{A}_i^{\mathrm{act}}(\widehat{\mathcal{R}})
\neq
\emptyset$.
\Halmos
\endproof

Lemma~\ref{lemma:single-coordinate-expansion} reduces a single coordinate expansion to determining when a critical arc becomes active.
Accordingly, for each $i\in[n]$, define the \emph{lower-facet expansion threshold} by
\begin{equation}
\label{eq:lower-facet-expansion-threshold}
    \theta_i^-(\mathcal{R})
    \coloneqq
    \max\left(
        \{-\infty\}
        \cup
        \left\{
            \dd{l}(\dd{a})\uparrow:
            \dd{a}\in
            \dd{A}_i^{\mathrm{crit}}(\mathcal{R}),
            \;
            \dd{l}(\dd{a})\uparrow \, \le l_i
        \right\}
    \right).
\end{equation}
Similarly, define the \emph{upper-facet expansion threshold} by
\begin{equation}
\label{eq:upper-facet-expansion-threshold}
    \theta_i^+(\mathcal{R})
    \coloneqq
    \min\left(
        \{+\infty\}
        \cup
        \left\{
            \dd{l}(\dd{a})\downarrow:
            \dd{a}\in
            \dd{A}_i^{\mathrm{crit}}(\mathcal{R}),
            \;
            \dd{l}(\dd{a})\downarrow \,\ge u_i
        \right\}
    \right).
\end{equation}
Thus, $\theta_i^-(\mathcal{R})$ is the first upper endpoint of a critical arc encountered as the lower facet is moved toward $-\infty$, whereas $\theta_i^+(\mathcal{R})$ is the first lower endpoint of a critical arc encountered as the upper facet is moved toward $+\infty$.
If no such critical arc exists in the corresponding direction, then
$\theta_i^-(\mathcal{R})=-\infty$ or
$\theta_i^+(\mathcal{R})=+\infty$, respectively.

For notational convenience, let $\sigma\in\{-,+\}$, where
$\sigma=-$ and $\sigma=+$ indicate expansion of the lower and
upper sides, respectively.
Thus, $(i,\sigma)\in[n]\times\{-,+\}$ specifies which side in
coordinate $i$ is expanded.
We write $\theta_i^\sigma(\mathcal{R})$ to encompass both cases, with
$\theta_i^-(\mathcal{R})$ when $\sigma=-$ and
$\theta_i^+(\mathcal{R})$ when $\sigma=+$.
Whenever the corresponding endpoint of $\mathcal{R}$ is finite,
we similarly write $F_i^\sigma(\mathcal{R})$ to encompass both cases, with
$F_i^-(\mathcal{R})$ when $\sigma=-$ and
$F_i^+(\mathcal{R})$ when $\sigma=+$.
For fixed $(i,\sigma)\in[n]\times\{-,+\}$, a
\emph{single-facet expansion} of $\mathcal{R}$ associated with $(i,\sigma)$ is a box obtained by replacing $l_i$ with a value no larger than $l_i$ if $\sigma=-$, or $u_i$ with a value no smaller than $u_i$ if $\sigma=+$, while leaving all other endpoints unchanged. 
The replacement value may be $-\infty$ or $+\infty$, respectively.

\begin{proposition}[Single-facet expansion threshold]
    \label{prop:single-facet-expansion-threshold}
    Let
    $\dd{D}=(\dd{U},\dd{A},\dd{l}(\cdot))$
    be an IDD, and let
    $\mathcal{R}=\prod_{j\in[n]}[l_j,u_j]$
    be a full-dimensional DD-free box with respect to $\dd{D}$.
    Fix $(i,\sigma)\in[n]\times\{-,+\}$.
    Let $\widehat{\mathcal{R}}$ be obtained from $\mathcal{R}$
    by replacing $l_i$ with $\theta_i^-(\mathcal{R})$
    if $\sigma=-$, or $u_i$ with $\theta_i^+(\mathcal{R})$
    if $\sigma=+$, while keeping all other endpoints unchanged.
    Then, among the single-facet expansions of $\mathcal{R}$
    associated with $(i,\sigma)$ that are DD-free,
    $\widehat{\mathcal{R}}$ is the largest with respect to set
    inclusion.
    Moreover, if $\theta_i^\sigma(\mathcal{R})$ is finite, then
    $\operatorname{relint}
    \bigl(F_i^\sigma(\widehat{\mathcal{R}})\bigr)
    \cap
    \operatorname{Sol}(\dd{D})
    \neq
    \emptyset$.
\end{proposition}

\proof{Proof.}

We prove the result for $\sigma=-$; 
the case $\sigma=+$ is symmetric.
It follows from \eqref{eq:lower-facet-expansion-threshold} that
$\theta_i^-(\mathcal{R})\le l_i$.
Thus, $\mathcal{R}\subseteq\widehat{\mathcal{R}}$
and
$\widehat{\mathcal{R}}$ is full-dimensional.
Since $\mathcal{R}$ is DD-free, 
Lemma~\ref{lemma:single-coordinate-expansion} implies that
$\dd{A}_i^{\mathrm{crit}}(\mathcal{R})
\cap
\dd{A}_i^{\mathrm{act}}(\mathcal{R})
=
\emptyset$.
Hence, every
$\dd{a}
\in
\dd{A}_i^{\mathrm{crit}}(\mathcal{R})$
satisfies either
$\dd{l}(\dd{a})\uparrow\le l_i$ or
$\dd{l}(\dd{a})\downarrow\ge u_i$.
For critical arcs on the lower side, the definition of 
$\theta_i^-(\mathcal{R})$ gives
$\dd{l}(\dd{a})\uparrow
\le
\theta_i^-(\mathcal{R})$;
critical arcs on the upper side remain inactive because $u_i$ is unchanged.
Therefore, no critical arc is active with respect to $\widehat{\mathcal{R}}$, and 
Lemma~\ref{lemma:single-coordinate-expansion} shows that 
$\widehat{\mathcal{R}}$ is DD-free.

If $\theta_i^-(\mathcal{R})=-\infty$, 
then no single-facet expansion of $\mathcal{R}$ associated with $(i,-)$ strictly contains $\widehat{\mathcal{R}}$, so maximality follows.
Now suppose that $\theta_i^-(\mathcal{R})>-\infty$.
Since $\dd{A}_i^{\mathrm{crit}}(\mathcal{R})$ is finite,
the maximum in \eqref{eq:lower-facet-expansion-threshold}
is attained by some
$\dd{a}^\ast
\in
\dd{A}_i^{\mathrm{crit}}(\mathcal{R})$
such that
$\dd{l}(\dd{a}^\ast)\uparrow
=
\theta_i^-(\mathcal{R})$.
It follows from \eqref{eq:critical-arcs} that this arc can be concatenated with an active 
$\dd{r}$--$\dd{tail}(\dd{a}^\ast)$ path and an active
$\dd{head}(\dd{a}^\ast)$--$\dd{t}$ path to obtain an
$\dd{r}$--$\dd{t}$ path
$\dd{p}=(\dd{a}_1,\ldots,\dd{a}_n)$,
where $\dd{a}_i=\dd{a}^\ast$ and every $\dd{a}_j$ with $j\neq i$ is active with respect to $\mathcal{R}$.
Set
$z_i\coloneqq\theta_i^-(\mathcal{R})$
and choose
$z_j\in\dd{l}(\dd{a}_j)\cap(l_j,u_j)$
for each $j\neq i$.
Since
$z_i=\dd{l}(\dd{a}^\ast)\uparrow
\in\dd{l}(\dd{a}^\ast)$,
we have
$\boldsymbol{z}\in\mathcal{S}_{\dd{p}}
\subseteq\operatorname{Sol}(\dd{D})$.
Moreover,
$z_i=\theta_i^-(\mathcal{R})$ and
$l_j<z_j<u_j$ for every $j\neq i$, so
$\boldsymbol{z}\in
\operatorname{relint}\bigl(F_i^-(\widehat{\mathcal{R}})\bigr)$.
Hence,
$\operatorname{relint}\bigl(F_i^-(\widehat{\mathcal{R}})\bigr)
\cap
\operatorname{Sol}(\dd{D})
\neq\emptyset$.
Any single-facet expansion associated with $(i,-)$ that strictly
contains $\widehat{\mathcal{R}}$ contains $\boldsymbol{z}$ in its
interior and is therefore not DD-free.
As a result, $\widehat{\mathcal{R}}$ is the largest DD-free single-facet
expansion of $\mathcal{R}$ associated with $(i,-)$.
\Halmos
\endproof

Starting from $\mathcal{R}^0$, we can apply the single-facet expansion method established above sequentially.
Let
$\Pi
=
\bigl(
    (i_1,\sigma_1),\ldots,(i_{2n},\sigma_{2n})
\bigr)$
be an arbitrary ordering of $[n]\times\{-,+\}$.
Each pair is processed once, using the expansion threshold associated with the current box.
An expansion can alter the reachability sets and, consequently, the thresholds used in subsequent iterations. 
Therefore, the resulting box generally depends on the ordering $\Pi$.

\begin{remark} 
\label{remark:reachable}
For a fixed box, all forward and backward reachability sets can be computed in $\mathcal{O}(|\dd{A}|)$ time.
Thus, recomputing these sets from scratch for each of the $2n$
single-facet expansions in the ordering $\Pi$ would require
$\mathcal{O}(n|\dd{A}|)$ time in total, and hence
$\mathcal{O}(|\dd{A}|^2)$ time when expressed solely in terms of
$|\dd{A}|$.
This quadratic-time repeated computation can be avoided and, in fact, reduced to linear time by maintaining the reachability sets incrementally.
Since each single-facet expansion replaces the current box by a
containing box, the set of active arcs and the forward and backward
reachability sets can only grow across iterations.
This enables the following update procedure, which can be carried out much more efficiently.

After each single-facet expansion, we first identify all newly active arcs and mark them active.
For the forward update, each newly active arc
$\dd{a}=(\dd{u},\dd{v})$ is processed.
If $\dd{u}$ is currently forward reachable and $\dd{v}$ is not, then $\dd{v}$ is marked forward reachable.
Whenever a node becomes forward reachable for the first time,
its outgoing arcs are scanned, and the same test is applied to
each active outgoing arc.
This propagation continues until no additional node is marked.
The backward update is symmetric, using incoming arcs and interchanging the roles of tails and heads.
As shown in 
Proposition~\ref{prop:box-expansion-time-complexity},
these incremental reachability updates require only $\mathcal{O}(|\dd{A}|)$ time over the entire expansion sequence.
\end{remark}
\medskip

\begin{algorithm}
\caption{Sequential box expansion}
\label{alg:box-expansion}
\begin{algorithmic}[1]
    \Require an IDD
    $\dd{D}=(\dd{U},\dd{A},\dd{l}(\cdot))$;
    an initial full-dimensional DD-free box $\mathcal{R}^0$;
    an ordering
    $\Pi=((i_1,\sigma_1),\ldots,(i_{2n},\sigma_{2n}))$
    of $[n]\times\{-,+\}$
    \Ensure a box-maximal DD-free box $\mathcal{R}^{2n}$
    containing $\mathcal{R}^0$
    \State initialize the active arcs with respect to $\mathcal{R}^0$
    \State compute
    $\dd{U}_j^{\mathrm{fwd}}(\mathcal{R}^0)$ and
    $\dd{U}_j^{\mathrm{bwd}}(\mathcal{R}^0)$
    for all $j\in[n+1]$
    \For{$q\in[2n]$}
        \State $(i,\sigma)\gets(i_q,\sigma_q)$
        \If{$\sigma=-$}
            \State compute
            $\theta_i^-(\mathcal{R}^{q-1})$
            using equation \eqref{eq:lower-facet-expansion-threshold}
            \State obtain $\mathcal{R}^{q}$ from
            $\mathcal{R}^{q-1}$ by replacing its lower endpoint in
            coordinate $i$ with
            $\theta_i^-(\mathcal{R}^{q-1})$
        \Else
            \State compute
            $\theta_i^+(\mathcal{R}^{q-1})$
            using equation \eqref{eq:upper-facet-expansion-threshold}
            \State obtain $\mathcal{R}^{q}$ from
            $\mathcal{R}^{q-1}$ by replacing its upper endpoint in
            coordinate $i$ with
            $\theta_i^+(\mathcal{R}^{q-1})$
        \EndIf
        \State
        $\Delta\dd{A}^q
        \gets
        \dd{A}_i^{\mathrm{act}}(\mathcal{R}^q)
        \setminus
        \dd{A}_i^{\mathrm{act}}(\mathcal{R}^{q-1})$
        \State mark all arcs in $\Delta\dd{A}^q$ as active
        \State update the forward and backward reachability sets from $\Delta\dd{A}^q$ using the incremental procedure in Remark~\ref{remark:reachable}.
    \EndFor
    \State\Return $\mathcal{R}^{2n}$
\end{algorithmic}
\end{algorithm}

Algorithm~\ref{alg:box-expansion} presents the sequential facet-expansion procedure.
We next establish the correctness and running-time complexity of
Algorithm~\ref{alg:box-expansion}.
In particular, Theorem~\ref{thm:box-expansion-valid} shows that the algorithm returns
a box-maximal DD-free box containing $\mathcal{R}^0$, while Proposition~\ref{prop:box-expansion-time-complexity} established that this box can be computed in time linear in the size of the IDD.

\begin{theorem}
\label{thm:box-expansion-valid}
    Let
    $\dd{D}=(\dd{U},\dd{A},\dd{l}(\cdot))$
    be an IDD, 
    let $\mathcal{R}^0$ be a full-dimensional DD-free box 
    with respect to $\dd{D}$, 
    and let $\Pi$ be an ordering of $[n]\times\{-,+\}$.
    Algorithm~\ref{alg:box-expansion} computes the exact forward
    and backward reachability sets at initialization and after
    each iteration of facet-expansion.
    Moreover, the generated boxes
    $\mathcal{R}^1,\ldots,\mathcal{R}^{2n}$ are full-dimensional
    and DD-free and satisfy
    $\mathcal{R}^0\subseteq\mathcal{R}^1\subseteq\cdots
    \subseteq\mathcal{R}^{2n}$.
    The box $\mathcal{R}^{2n}$ returned in the last iteration is box-maximal and,
    consequently, maximal among all DD-free sets.
\end{theorem}

\proof{Proof.}

For each $q\in\{0,\ldots,2n\}$, write
$\dd{A}^{q,\mathrm{act}}
\coloneqq
\dd{A}^{\mathrm{act}}(\mathcal{R}^q)
=
\bigcup_{i\in[n]}
\dd{A}_i^{\mathrm{act}}(\mathcal{R}^q)$.
Let
$\dd{U}^{q,\mathrm{fwd}}
\coloneqq
\bigcup_{j\in[n+1]}
\dd{U}_j^{\mathrm{fwd}}(\mathcal{R}^q)$
and
$\dd{U}^{q,\mathrm{bwd}}
\coloneqq
\bigcup_{j\in[n+1]}
\dd{U}_j^{\mathrm{bwd}}(\mathcal{R}^q)$.
We prove by induction on $q$ that, for every
$q\in\{0,\ldots,2n\}$, $\mathcal{R}^q$ is a full-dimensional DD-free box and the forward and backward reachability sets maintained by the algorithm after line 14 are exactly
$\dd{U}^{q,\mathrm{fwd}}$ and
$\dd{U}^{q,\mathrm{bwd}}$, respectively.
In each inductive step, we also establish that
$\mathcal{R}^{q-1}\subseteq\mathcal{R}^q$.

\emph{Base case ($q=0$).}
By assumption, $\mathcal{R}^0$ is a full-dimensional DD-free box.
The initialization steps identify exactly the arcs in
$\dd{A}^{0,\mathrm{act}}$ and compute
$\dd{U}^{0,\mathrm{fwd}}$ and
$\dd{U}^{0,\mathrm{bwd}}$.
Hence the claim holds for $q=0$.

\emph{Inductive step.}
Fix $q\in[2n]$ and suppose that the claim holds for $q-1$.
Let $(i,\sigma)=(i_q,\sigma_q)$ be the facet processed in iteration $q$, and write the $i$-th coordinate interval of $\mathcal{R}^{q-1}$ as $[l_i,u_i]$.
By the induction hypothesis, the maintained forward and backward reachability sets for $\mathcal{R}^{q-1}$ are
$\dd{U}^{q-1,\mathrm{fwd}}$ and
$\dd{U}^{q-1,\mathrm{bwd}}$, respectively.
Hence, the set of critical arcs used to compute the threshold in iteration $q$ is $\dd{A}_i^{\mathrm{crit}}(\mathcal{R}^{q-1})$.

If $\sigma=-$, then
$\theta_i^-(\mathcal{R}^{q-1})\le l_i$ by definition, and
$\mathcal{R}^q$ is obtained by replacing $l_i$ with this threshold.
If $\sigma=+$, then
$\theta_i^+(\mathcal{R}^{q-1})\ge u_i$, and
$\mathcal{R}^q$ is obtained by replacing $u_i$ with this threshold.
In either case,
$\mathcal{R}^{q-1}\subseteq\mathcal{R}^q$, and
$\mathcal{R}^q$ remains full-dimensional.
Proposition~\ref{prop:single-facet-expansion-threshold} further implies
that $\mathcal{R}^q$ is DD-free.

Because only the $i$-th coordinate interval is enlarged, every arc active with respect to $\mathcal{R}^{q-1}$ remains active with respect 
to $\mathcal{R}^q$.
Thus,
$\dd{A}^{q-1,\mathrm{act}}
\subseteq
\dd{A}^{q,\mathrm{act}}$.
Let
$\Delta\dd{A}^q
\coloneqq
\dd{A}^{q,\mathrm{act}}
\setminus
\dd{A}^{q-1,\mathrm{act}}$
denote the set of newly active arcs.
Since only the $i$-th coordinate interval changes,
$\Delta\dd{A}^q
=
\dd{A}_i^{\mathrm{act}}(\mathcal{R}^q)
\setminus
\dd{A}_i^{\mathrm{act}}(\mathcal{R}^{q-1})$,
which is precisely the set computed in line 12 of iteration $q$ of Algorithm~\ref{alg:box-expansion}.

Consider the forward reachability update.
By the incremental update rule in Remark~\ref{remark:reachable}, the algorithm retains all nodes in 
$\dd{U}^{q-1,\mathrm{fwd}}$, processes every arc in
$\Delta\dd{A}^q$, and processes every active outgoing arc of a node when that node is first marked forward reachable.
We show that the maintained forward reachability set is exactly
$\dd{U}^{q,\mathrm{fwd}}$ by proving both forward and reverse inclusions.

First, since
$\dd{A}^{q-1,\mathrm{act}}
\subseteq
\dd{A}^{q,\mathrm{act}}$,
every node in $\dd{U}^{q-1,\mathrm{fwd}}$ remains reachable from
$\dd{r}$ under $\mathcal{R}^q$.
Every node marked subsequently is the head of an arc in
$\dd{A}^{q,\mathrm{act}}$ whose tail is already marked.
An induction on the order in which nodes are marked therefore shows that every marked node is reachable from $\dd{r}$ using arcs in
$\dd{A}^{q,\mathrm{act}}$.
Hence, the maintained forward reachability set is contained in
$\dd{U}^{q,\mathrm{fwd}}$.

For the reverse inclusion, fix
$\dd{v}\in\dd{U}^{q,\mathrm{fwd}}$ and choose an
$\dd{r}$--$\dd{v}$ path consisting of arcs in
$\dd{A}^{q,\mathrm{act}}$.
We show, in order along this path, that every node on the path is marked.
The root node $\dd{r}$ belongs to
$\dd{U}^{q-1,\mathrm{fwd}}$ and is therefore marked.
Now consider an arc
$\dd{a}=(\dd{u},\dd{u}')$ on the path and suppose that
$\dd{u}$ is marked.
If
$\dd{u}\notin\dd{U}^{q-1,\mathrm{fwd}}$, then $\dd{u}$ was first
marked during the current update, so all of its active outgoing arcs,
including $\dd{a}$, are processed.
If
$\dd{u}\in\dd{U}^{q-1,\mathrm{fwd}}$, then either
$\dd{a}\in\dd{A}^{q-1,\mathrm{act}}$, in which case
$\dd{u}'\in\dd{U}^{q-1,\mathrm{fwd}}$, or
$\dd{a}\in\Delta\dd{A}^q$, in which case $\dd{a}$ is processed as a newly active arc.
Thus, in every case, $\dd{u}'$ is marked.
Proceeding along the path shows that $\dd{v}$ is marked.
Hence,
$\dd{U}^{q,\mathrm{fwd}}$ is contained in the maintained forward
reachability set, proving equality.

The backward reachability update follows the symmetric incremental rule.
The algorithm retains all nodes in
$\dd{U}^{q-1,\mathrm{bwd}}$, processes every arc in
$\Delta\dd{A}^q$, and processes every active incoming arc of a node when that node is first marked backward reachable.
Applying the above two-inclusion argument in the reverse direction, with $\dd{t}$ in place of $\dd{r}$ and incoming arcs in place of outgoing arcs, shows that the maintained backward reachability set is exactly
$\dd{U}^{q,\mathrm{bwd}}$.

This completes the induction and establishes the claimed reachability sets, DD-freeness, and nestedness.
Each incremental reachability update terminates because it only marks previously unmarked nodes and $\dd{U}$ is finite.
Thus, the algorithm completes its $2n$ iterations.

It remains to prove maximality of the returned box.
For each $q\in\{0,\ldots,2n\}$, write
$\mathcal{R}^q=\prod_{j\in[n]}[l_j^q,u_j^q]$.
Fix a finite lower facet $F_i^-(\mathcal{R}^{2n})$, and let
$q$ be the unique iteration with $(i_q,\sigma_q)=(i,-)$.
Since the lower endpoint in coordinate $i$ is not modified after iteration $q$, we have
$l_i^{2n}
=
l_i^q
=
\theta_i^-(\mathcal{R}^{q-1})
\in\mathbb{R}$.
Proposition~\ref{prop:single-facet-expansion-threshold}
therefore provides a point
$\boldsymbol{z}
\in
\operatorname{relint}\bigl(F_i^-(\mathcal{R}^q)\bigr)
\cap
\operatorname{Sol}(\dd{D})$.
For every $j\neq i$, nestedness gives
\[
    l_j^{2n}
    \le
    l_j^q
    <
    z_j
    <
    u_j^q
    \le
    u_j^{2n}.
\]
Together with
$z_i=l_i^q=l_i^{2n}$, this implies
\[
    \boldsymbol{z}
    \in
    \operatorname{relint}\bigl(F_i^-(\mathcal{R}^{2n})\bigr)
    \cap
    \operatorname{Sol}(\dd{D}).
\]
Thus, every finite lower facet of the returned box satisfies the facet-contact condition (ii) in Proposition~\ref{prop:box-maximal-characterization}.
The argument for finite upper facets is symmetric.
Proposition~\ref{prop:box-maximal-characterization} now implies that
$\mathcal{R}^{2n}$ is box-maximal and maximal among all
DD-free sets.
\Halmos
\endproof

\begin{proposition}
\label{prop:box-expansion-time-complexity}
    Let
    $\dd{D}=(\dd{U},\dd{A},\dd{l}(\cdot))$
    be an IDD with
    $\operatorname{Sol}(\dd{D})\neq\emptyset$,
    let $\mathcal{R}^0$ be a full-dimensional DD-free box
    with respect to $\dd{D}$,
    and let $\Pi$ be an ordering of $[n]\times\{-,+\}$.
    Then Algorithm~\ref{alg:box-expansion} can be implemented in
    $\mathcal{O}(|\dd{A}|)$ time.
\end{proposition}

\proof{Proof.}

Consider the IDD representation by outgoing and incoming adjacency lists, with the arcs also indexed by layer, and maintain Boolean flags indicating arc activity and forward and backward node reachability.
One scan of the arcs identifies the active arcs with respect to  $\mathcal{R}^0$, and forward and backward traversals of the resulting active subgraph compute the initial forward and backward reachability sets.
Together, these initialization steps require
$\mathcal{O}(|\dd{U}|+|\dd{A}|)$ time.

Fix an iteration $q\in[2n]$ and let
$(i,\sigma)=(i_q,\sigma_q)$.
Using the maintained reachability flags, membership in
$\dd{A}_i^{\mathrm{crit}}(\mathcal{R}^{q-1})$ can be tested in constant time for each arc in $\dd{A}_i$.
Hence, one scan of $\dd{A}_i$
computes $\theta_i^\sigma(\mathcal{R}^{q-1})$.
After the corresponding endpoint is updated, a second scan of
$\dd{A}_i$ identifies the arcs in $\Delta\dd{A}^q$.
Since $\Pi$ is an ordering of $[n]\times\{-,+\}$, each pair
$(i,\sigma)$ is processed exactly once.
Hence, each arc participates in at most four such layer scans, and all threshold computations and newly active arc identifications require
$\mathcal{O}(|\dd{A}|)$ time in total.

It remains to bound the incremental reachability updates as described 
in Remark~\ref{remark:reachable}.
By Theorem~\ref{thm:box-expansion-valid}, the generated boxes are nested, so the sets of active arcs and forward and backward reachable nodes are monotone across iterations.
Consequently, each arc becomes newly active at most once, and each node becomes newly reachable at most once in each direction.
When a node first becomes forward reachable, its outgoing adjacency list is scanned once; when it first becomes backward reachable, its incoming adjacency list is scanned once.
Thus, each arc is examined at most once through the forward scan of its tail node and at most once through the backward scan of its head node.
Moreover, since each arc belongs to $\Delta\dd{A}^q$ for at most one iteration $q$, it is tested at most once in each direction upon activation.
Hence, each arc is examined only a constant number of times over all incremental updates.
Since each node is marked at most once in each direction and the sums of the outgoing and incoming degrees are both $|\dd{A}|$, all incremental reachability updates require
$\mathcal{O}(|\dd{U}|+|\dd{A}|)$ time.

With the current box stored by its endpoints and updated in place,
the remaining work over the $2n$ iterations requires
$\mathcal{O}(n)$ time.
Thus, the total running time is
$\mathcal{O}(|\dd{U}|+|\dd{A}|+n)$.
By \eqref{eq:DD-size-bounds},
$|\dd{U}|\le2|\dd{A}|$ and $n\le|\dd{A}|$.
Therefore, Algorithm~\ref{alg:box-expansion} can be implemented in
$\mathcal{O}(|\dd{A}|)$ time.
\Halmos
\endproof

Algorithm~\ref{alg:box-expansion} also provides an efficient test for box-maximality as presented in Proposition~\ref{prop:box-maximal-characterization}.
Indeed, Theorem~\ref{thm:box-expansion-valid} guarantees that, for any full-dimensional DD-free box $\mathcal{R}^0$, the algorithm returns a box-maximal DD-free set $\mathcal{R}^{2n}$ containing $\mathcal{R}^0$.
Hence, $\mathcal{R}^0$ is box-maximal if and only if 
$\mathcal{R}^{2n}=\mathcal{R}^0$.
By Proposition~\ref{prop:box-expansion-time-complexity},
this test can be performed in $\mathcal{O}(|\dd{A}|)$ time.

\begin{figure}[tb]
    \centering
    \includegraphics[width=\linewidth]{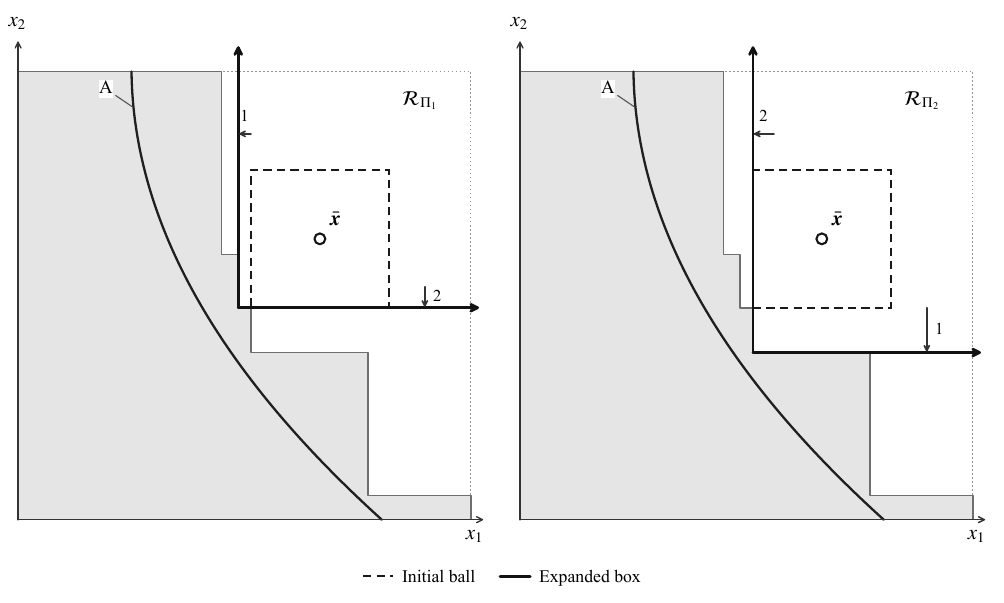}
    \caption{Box-maximal DD-free sets obtained by sequential
        facet expansion.
        The solid curve labeled $A$ depicts the nonlinear
        boundary of $\mathcal{G}^A$, with the feasible region
        lying on the side containing the origin.
        The shaded region represents
        $\operatorname{Sol}(\dd{D}^A)$.
        In each panel, the dashed square bounds the initial
        DD-free $\ell_\infty$-ball
        $\mathcal{R}^0
        =\mathcal{C}_\infty(\bar{\boldsymbol{x}})$,
        and the thick solid rays bound the expanded boxes
        $\mathcal{R}_{\Pi_1}$ (left) and
        $\mathcal{R}_{\Pi_2}$ (right), which are unbounded
        in both positive coordinate directions.
        The numbered arrows indicate the order in which
        the two lower facets are processed.
        }
    \label{fig:box-expansion-ordering}
\end{figure}

To illustrate the box expansion procedure introduced above, we return to the example in
Section~\ref{subsec:illustrative-example},
focusing on constraint $A$ with the same IDD relaxation
$\dd{D}^A$ and point $\bar{\boldsymbol{x}}$.
Figure~\ref{fig:box-expansion-ordering} shows the
expansion of the initial DD-free $\ell_\infty$-ball
$\mathcal{R}^0=\mathcal{C}_\infty(\bar{\boldsymbol{x}})$
by Algorithm~\ref{alg:box-expansion}.
The left and right panels depict two resulting 
box-maximal DD-free boxes,
$\mathcal{R}_{\Pi_1}$ and $\mathcal{R}_{\Pi_2}$,
obtained under different facet orderings.
In both cases, the expanded box contains
$\mathcal{R}^0$ and is unbounded in both positive
coordinate directions, while each finite facet
contains a contact point with
$\operatorname{Sol}(\dd{D}^A)$ in its relative interior.

\subsubsection{Intersection Cuts from Expanded Boxes}
\label{subsubsec:box-IC}

We now return to the corner relaxation setting of
Section~\ref{subsec:IC-derivation}.
Let 
$\bar{\boldsymbol{x}}
\in
\mathcal{P}
\setminus
\operatorname{Sol}(\dd{D})$ 
be a basic feasible solution of $\mathcal{P}$,
and fix an associated corner relaxation 
$\mathcal{P}'
=\bar{\boldsymbol{x}}
+\operatorname{cone}\{
    \boldsymbol{v}^j:j\in R
\}$.
Recall that $R$ indexes the extreme rays of $\mathcal{P}'$, and every
$\boldsymbol{x}\in\mathcal{P}'$ admits a unique representation
$\boldsymbol{x}
=\bar{\boldsymbol{x}}
+\sum_{j\in R}
\lambda_j(\boldsymbol{x})\boldsymbol{v}^j$,
where $\lambda_j(\boldsymbol{x})\ge0$ for all $j\in R$.

\begin{corollary}
\label{cor:box-IC}
    Let
    $\mathcal{R}^{2n}=\prod_{i\in[n]}[l_i,u_i]$
    be the box returned by Algorithm~\ref{alg:box-expansion}  under an arbitrary
    facet ordering, initialized at
    $\mathcal{R}^0
    =
    \mathcal{C}_\infty(\bar{\boldsymbol{x}})$.
    For each $i\in[n]$, define
    \[
        \rho_i^-
        \coloneqq
        \begin{cases}
            \bar{x}_i-l_i,
            & l_i\in\mathbb{R},\\
            +\infty,
            & l_i=-\infty,
        \end{cases}
        \qquad
        \rho_i^+
        \coloneqq
        \begin{cases}
            u_i-\bar{x}_i,
            & u_i\in\mathbb{R},\\
            +\infty,
            & u_i=+\infty.
        \end{cases}
    \]
    Then the intersection cut obtained from $\mathcal{R}^{2n}$ on $\mathcal{P}'$ is
    \begin{equation}
    \label{eq:box-IC}
        \sum_{j\in R}
        \max_{i\in[n]}
        \left\{
            \frac{(v_i^j)_+}{\rho_i^+},
            \frac{(-v_i^j)_+}{\rho_i^-}
        \right\}
        \lambda_j(\boldsymbol{x})
        \ge
        1.
    \end{equation}
    This inequality cuts off $\bar{\boldsymbol{x}}$ and is valid for 
    $\mathcal{P}
    \cap
    \operatorname{Sol}(\dd{D})$ and hence for 
    $\mathcal{P}\cap\mathcal{G}$.
\end{corollary}

\proof{Proof.}

Since
$\bar{\boldsymbol{x}}
\notin
\operatorname{Sol}(\dd{D})$,
the initial DD-free box
$\mathcal{R}^0
=\mathcal{C}_\infty(\bar{\boldsymbol{x}})$
contains $\bar{\boldsymbol{x}}$ in its interior.
By Theorem~\ref{thm:box-expansion-valid},
$\mathcal{R}^{2n}$ is DD-free and contains $\mathcal{R}^0$.
Hence,
$\bar{\boldsymbol{x}}
\in
\operatorname{int}(\mathcal{R}^{2n})$
and
$\rho_i^-,\rho_i^+
\in(0,+\infty]$ 
for all $i\in[n]$.
By definition,
$\mathcal{R}^{2n}-\bar{\boldsymbol{x}}
=\prod_{i\in[n]}[-\rho_i^-,\rho_i^+]$.
For every
$\boldsymbol{w}\in\mathbb{R}^n$ and $t>0$, the condition
$\boldsymbol{w}
\in
t(\mathcal{R}^{2n}-\bar{\boldsymbol{x}})$
is equivalent to the coordinatewise inequalities
\[
    t
    \ge
    \frac{(w_i)_+}{\rho_i^+}
    \qquad\text{and}\qquad
    t
    \ge
    \frac{(-w_i)_+}{\rho_i^-},
    \qquad
    \forall i\in[n].
\]
Taking the infimum over $t>0$ gives
\begin{equation}
\label{eq:box-gauge}
    \gamma_{\mathcal{R}^{2n}-\bar{\boldsymbol{x}}}
    (\boldsymbol{w})
    =
    \max_{i\in[n]}
    \left\{
        \frac{(w_i)_+}{\rho_i^+},
        \frac{(-w_i)_+}{\rho_i^-}
    \right\}.
\end{equation}
The terms corresponding to unbounded sides of
$\mathcal{R}^{2n}$ are zero.
Substituting \eqref{eq:box-gauge} into
\eqref{eq:IC-gauge} yields \eqref{eq:box-IC}.
By the validity argument in
Section~\ref{subsec:IC-derivation},
inequality \eqref{eq:box-IC} is valid for
$\mathcal{P}
\cap
\operatorname{Sol}(\dd{D})$
and hence for 
$\mathcal{P}\cap\mathcal{G}$.
Finally,
$\lambda_j(\bar{\boldsymbol{x}})=0$ for every $j\in R$,
so $\bar{\boldsymbol{x}}$ violates \eqref{eq:box-IC}.
\Halmos
\endproof

The inclusion
$\mathcal{R}^0\subseteq\mathcal{R}^{2n}$
yields the following comparison with the initial
$\ell_\infty$-ball intersection cut on the same corner relaxation.

\begin{corollary}
\label{cor:box-IC-strength}
    Let
    $r\coloneqq
    d_\infty\bigl(
        \bar{\boldsymbol{x}},
        \operatorname{Sol}(\dd{D})
    \bigr)$,
    and let $\mathcal{R}^{2n}$ be the box returned by
    Algorithm~\ref{alg:box-expansion} from
    $\mathcal{R}^0=\mathcal{C}_\infty(\bar{\boldsymbol{x}})$
    under any sequential expansion ordering $\Pi$ of
    $[n]\times\{-,+\}$.
    Then, for every $j\in R$,
    $\gamma_{\mathcal{R}^{2n}-\bar{\boldsymbol{x}}}
    (\boldsymbol{v}^j)
    \le
    \lVert\boldsymbol{v}^j\rVert_\infty/r$.
    Consequently, \eqref{eq:box-IC} is at least as strong on the corner relaxation $\mathcal{P}'$ as the $\ell_\infty$-ball intersection cut
    obtained from \eqref{eq:lp-ball-IC} by setting $p=\infty$.
\end{corollary}

\proof{Proof.}

By Theorem~\ref{thm:box-expansion-valid},
$\mathcal{B}_\infty(\bar{\boldsymbol{x}},r)
=
\mathcal{R}^0
\subseteq
\mathcal{R}^{2n}$.
Translating both sets by $-\bar{\boldsymbol{x}}$ and using
gauge monotonicity under set inclusion gives
$\gamma_{\mathcal{R}^{2n}-\bar{\boldsymbol{x}}}
(\boldsymbol{v}^j)
\le
\gamma_{\mathcal{B}_\infty(\boldsymbol{0},r)}
(\boldsymbol{v}^j)
=
\lVert\boldsymbol{v}^j\rVert_\infty/r$
for every $j\in R$.
Since $\lambda_j(\boldsymbol{x})\ge0$ for every $j\in R$ and
$\boldsymbol{x}\in\mathcal{P}'$, every point of $\mathcal{P}'$
satisfying \eqref{eq:box-IC} also satisfies the initial
$\ell_\infty$-ball intersection cut.
\Halmos
\endproof

We now consider the cut closures defined in
\eqref{eq:rank-t-cut-closure} and reinstate the constraint index.
For each $k\in K$, fix any sequential expansion ordering $\Pi^k$ of
$[n]\times\{-,+\}$.
For
$\boldsymbol{x}\notin\operatorname{Sol}(\dd{D}^k)$,
let
$\mathcal{C}_{\mathrm{box}}^k(\boldsymbol{x})$
denote the box returned by Algorithm~\ref{alg:box-expansion}
applied to $\dd{D}^k$ with initial box
$\mathcal{C}_\infty^k(\boldsymbol{x})$ and expansion ordering $\Pi^k$.
For
$\boldsymbol{x}\in\operatorname{Sol}(\dd{D}^k)$,
set
$\mathcal{C}_{\mathrm{box}}^k(\boldsymbol{x})\coloneqq\emptyset$.

\begin{corollary}
\label{cor:box-cut-closure-convergence}
    For every $k\in K$, the mapping
    $\mathcal{C}_{\mathrm{box}}^k(\cdot)$ defined above
    is a DD-free-set oracle
    for $\dd{D}^k$, and the resulting family
    $\{\mathcal{C}_{\mathrm{box}}^k(\cdot)\}_{k\in K}$
    satisfies the local covering condition
    \eqref{eq:local-covering-condition}.
    Consequently, if
    $\{\mathcal{P}_t\}_{t\in\mathbb{W}}$
    denotes the cut closures in
    \eqref{eq:rank-t-cut-closure} induced by this oracle family,
    then
    \[
        \mathcal{P}_\infty
        =
        \operatorname{conv}\left(
            \mathcal{P}
            \cap
            \bigcap_{k\in K}
            \operatorname{Sol}(\dd{D}^k)
        \right).
    \]
\end{corollary}

\proof{Proof.}

Fix $k\in K$ and
$\boldsymbol{x}\notin\operatorname{Sol}(\dd{D}^k)$.
By Proposition~\ref{prop:lp-ball-DD-free},
$\mathcal{C}_\infty^k(\boldsymbol{x})$ is a DD-free box
containing $\boldsymbol{x}$ in its interior and hence is
full-dimensional.
Applying Theorem~\ref{thm:box-expansion-valid} shows that
$\mathcal{C}_{\mathrm{box}}^k(\boldsymbol{x})$ is DD-free and
contains $\mathcal{C}_\infty^k(\boldsymbol{x})$.
Therefore,
$\boldsymbol{x}\in
\operatorname{int}\bigl(
    \mathcal{C}_{\mathrm{box}}^k(\boldsymbol{x})
\bigr)$.
Together with the definition on
$\operatorname{Sol}(\dd{D}^k)$, this proves that
$\mathcal{C}_{\mathrm{box}}^k(\cdot)$ is a DD-free-set oracle
for $\dd{D}^k$.

Moreover,
$\mathcal{C}_\infty^k(\boldsymbol{x})
\subseteq
\mathcal{C}_{\mathrm{box}}^k(\boldsymbol{x})$
for every $k\in K$ and $\boldsymbol{x}\in\mathbb{R}^n$.
Corollary~\ref{cor:lp-ball-enlargement-cut-closure-convergence}
therefore applies with $p=\infty$, yielding the local covering
condition and the stated cut-closure equality.
\Halmos
\endproof

\subsection{$\ell_p$-Ball-Maximal DD-Free Sets}
\label{subsec:lp-ball-maximal-DD-free}

In this section, we focus on DD-free families defined by $\ell_p$-balls for $p\in[1,\infty)$, characterize $\ell_p$-ball maximality, and study procedures for expanding the initial DD-free ball $\mathcal{C}_p(\bar{\boldsymbol{x}})$.
Fix a contact point of $\mathcal{C}_p(\bar{\boldsymbol{x}})$.
Beginning with $\mathcal{C}_p(\bar{\boldsymbol{x}})$, we construct a nested family of $\ell_p$-balls, parameterized by their radius, whose centers move along the ray from the contact point through $\bar{\boldsymbol{x}}$, while the contact point remains on the boundary.
We call this construction the \emph{radial expansion} induced by the selected contact point.

For $p\in(1,\infty)$, we show that every DD-free $\ell_p$-ball
containing $\mathcal{C}_p(\bar{\boldsymbol{x}})$ belongs to this family.
Consequently, if there is a largest radius for which a ball in this family is DD-free, the corresponding ball is the unique largest
DD-free $\ell_p$-ball containing
$\mathcal{C}_p(\bar{\boldsymbol{x}})$ and is therefore
$\ell_p$-ball-maximal.
If this family contains DD-free balls of arbitrarily large radius, no such $\ell_p$-ball-maximal DD-free ball exists; instead, the unbounded case yields a DD-free halfspace containing the entire family, which is itself a maximal DD-free set.

For $p=1$, however, the resulting family need not include every DD-free $\ell_1$-ball containing
$\mathcal{C}_1(\bar{\boldsymbol{x}})$.
Consequently, a largest DD-free member of this family, when it exists, need not be $\ell_1$-ball-maximal.
We instead characterize $\ell_1$-ball-maximality by examining how the path-induced hyperrectangles touching the ball boundary constrain expansion along the coordinate directions.

The remainder of this section is organized as follows.
We first study boundary contacts and, for $p\in[1,\infty)$, characterize $\ell_p$-ball-maximality.
We then analyze the radial expansion induced by a selected contact point and develop an IDD-based procedure for constructing the resulting expanded DD-free set.
Finally, we derive intersection cuts from these expanded sets, compare their strength with the initial $\ell_p$-ball cut, and establish the corresponding cut-closure result.

\subsubsection{Boundary Contacts and $\ell_p$-Ball-Maximality}
\label{subsubsec:ball-contact-geometry}

For $p\in[1,\infty)$ and a DD-free $\ell_p$-ball 
$\mathcal{B}_p(\boldsymbol{c},r)$ with $r>0$,
define its \emph{contact set} by
\[
    \mathcal{Z}_p(\boldsymbol{c},r)
    \coloneqq
    \operatorname{Sol}(\dd{D})
    \cap
    \operatorname{bd}
    \bigl(\mathcal{B}_p(\boldsymbol{c},r)\bigr).
\]
Whenever this set is nonempty, we have
$r=d_p(\boldsymbol{c},\operatorname{Sol}(\dd{D}))$, 
and the contact points are precisely the nearest points in 
$\operatorname{Sol}(\dd{D})$ to
$\boldsymbol{c}$ under the $\ell_p$-norm.
In this case, 
Proposition~\ref{prop:ball-contact-subgraph} provides
an efficiently computable DD-based representation of this contact set through coordinatewise projections onto the hyperrectangles induced by shortest
$\dd{r}$--$\dd{t}$ paths.

The following lemma shows a basic inclusion property for $\ell_p$-balls.

\begin{lemma}
\label{lemma:lp-ball-containment}
    Let $p\in[1,\infty)$,
    $\boldsymbol{c}_1,\boldsymbol{c}_2\in\mathbb{R}^n$, and
    $r_1,r_2\ge0$.
    Then
    $\mathcal{B}_p(\boldsymbol{c}_1,r_1)
    \subseteq
    \mathcal{B}_p(\boldsymbol{c}_2,r_2)$
    if and only if
    $\lVert
        \boldsymbol{c}_2-\boldsymbol{c}_1
    \rVert_p+r_1
    \le r_2$.
\end{lemma}

\proof{Proof.}

Suppose first that
$\mathcal{B}_p(\boldsymbol{c}_1,r_1)
\subseteq
\mathcal{B}_p(\boldsymbol{c}_2,r_2)$.
If $\boldsymbol{c}_1=\boldsymbol{c}_2$, 
then the desired inequality reduces to $r_1\le r_2$, which follows immediately from the containment.
Otherwise, let 
$\boldsymbol{w}
\coloneqq 
(\boldsymbol{c}_1-\boldsymbol{c}_2)
/
\lVert\boldsymbol{c}_1-\boldsymbol{c}_2\rVert_p$
so that $\lVert\boldsymbol{w}\rVert_p=1$,
and set
$\boldsymbol{x}\coloneqq\boldsymbol{c}_1+r_1\boldsymbol{w}$.
Then
$\boldsymbol{x}\in\mathcal{B}_p(\boldsymbol{c}_1,r_1)$
and hence
$\boldsymbol{x}\in\mathcal{B}_p(\boldsymbol{c}_2,r_2)$.
Since
$\boldsymbol{c}_1-\boldsymbol{c}_2
=
\lVert\boldsymbol{c}_1-\boldsymbol{c}_2\rVert_p
\boldsymbol{w}$,
we have
$\boldsymbol{x}-\boldsymbol{c}_2
=
(\lVert\boldsymbol{c}_1-\boldsymbol{c}_2\rVert_p+r_1)
\boldsymbol{w}$.
Thus,
$\lVert\boldsymbol{x}-\boldsymbol{c}_2\rVert_p
=
\lVert
    \boldsymbol{c}_1-\boldsymbol{c}_2
\rVert_p+r_1
\le r_2$.

Conversely, suppose that
$\lVert
    \boldsymbol{c}_2-\boldsymbol{c}_1
\rVert_p+r_1
\le r_2$, and let
$\boldsymbol{x}\in\mathcal{B}_p(\boldsymbol{c}_1,r_1)$.
By the triangle inequality,
$\lVert
    \boldsymbol{x}-\boldsymbol{c}_2
\rVert_p
\le
\lVert
    \boldsymbol{x}-\boldsymbol{c}_1
\rVert_p
+
\lVert
    \boldsymbol{c}_1-\boldsymbol{c}_2
\rVert_p
\le
r_1 + 
\lVert
    \boldsymbol{c}_1-\boldsymbol{c}_2
\rVert_p
\le r_2$.
Thus,
$\boldsymbol{x}\in\mathcal{B}_p(\boldsymbol{c_2},r_2)$.
Since $\boldsymbol{x}$ was arbitrary, it follows that
$\mathcal{B}_p(\boldsymbol{c}_1,r_1)
\subseteq
\mathcal{B}_p(\boldsymbol{c}_2,r_2)$.
\Halmos
\endproof

The following lemma characterizes that, for $p\in(1,\infty)$, once a
contact point $\boldsymbol{z}$ of a DD-free $\ell_p$-ball
$\mathcal{B}_p(\boldsymbol{c},r)$ is fixed, every DD-free
$\ell_p$-ball containing $\mathcal{B}_p(\boldsymbol{c},r)$ must
have its center on the ray from $\boldsymbol{z}$ through
$\boldsymbol{c}$.

\begin{lemma}
\label{lemma:ball-containment-rigidity}
    Let $\dd{D}$ be an IDD and let $p\in(1,\infty)$.
    Let $\mathcal{B}_p(\boldsymbol{c},r)$ be a DD-free
    $\ell_p$-ball with respect to $\dd{D}$, where $r>0$,
    and let
    $\boldsymbol{z}\in\mathcal{Z}_p(\boldsymbol{c},r)$.
    If another DD-free $\ell_p$-ball
    $\mathcal{B}_p(\widehat{\boldsymbol{c}},\widehat{r})$
    contains $\mathcal{B}_p(\boldsymbol{c},r)$, then
    $\widehat{r}\ge r$ and
    \begin{equation}
    \label{eq:containing-ball-center}
        \widehat{\boldsymbol{c}}
        =
        \boldsymbol{z}
        +
        \frac{\widehat{r}}{r}
        (\boldsymbol{c}-\boldsymbol{z}).
    \end{equation}
\end{lemma}

\proof{Proof.}

By Lemma~\ref{lemma:lp-ball-containment}, 
$\mathcal{B}_p(\boldsymbol{c},r)
\subseteq
\mathcal{B}_p(\widehat{\boldsymbol{c}},\widehat{r})$
gives
$\lVert
    \widehat{\boldsymbol{c}}-\boldsymbol{c}
\rVert_p+r
\le
\widehat{r}$,
and hence $\widehat{r}\ge r$.
Since
$\boldsymbol{z}\in\mathcal{Z}_p(\boldsymbol{c},r)$,
we have
$\boldsymbol{z}\in\operatorname{Sol}(\dd{D})$
and
$\lVert\boldsymbol{c}-\boldsymbol{z}\rVert_p=r$.
The inclusion
$\mathcal{B}_p(\boldsymbol{c},r)
\subseteq
\mathcal{B}_p(\widehat{\boldsymbol{c}},\widehat{r})$
then implies that
$\boldsymbol{z}
\in
\mathcal{B}_p(\widehat{\boldsymbol{c}},\widehat{r})$.
Since 
$\mathcal{B}_p(\widehat{\boldsymbol{c}},\widehat{r})$ is DD-free,
$\boldsymbol{z}$ cannot lie in its interior.
Thus,
$\boldsymbol{z}
\in
\operatorname{bd}\bigl(
    \mathcal{B}_p(\widehat{\boldsymbol{c}},\widehat{r})
\bigr)$,
which implies that
$\lVert
    \widehat{\boldsymbol{c}}-\boldsymbol{z}
\rVert_p
=
\widehat{r}$.
Combining these relations with the triangle inequality yields
\[
    \widehat{r}
    =
    \lVert
        \widehat{\boldsymbol{c}}-\boldsymbol{z}
    \rVert_p
    \le
    \lVert
        \widehat{\boldsymbol{c}}-\boldsymbol{c}
    \rVert_p
    +
    \lVert
        \boldsymbol{c}-\boldsymbol{z}
    \rVert_p
    =
    \lVert
        \widehat{\boldsymbol{c}}-\boldsymbol{c}
    \rVert_p+r
    \le
    \widehat{r}.
\]
Since the chain begins and ends at $\widehat{r}$, both inequalities must hold with equality.
In particular,
$\lVert
    \widehat{\boldsymbol{c}}-\boldsymbol{z}
\rVert_p
=
\lVert
    \widehat{\boldsymbol{c}}-\boldsymbol{c}
\rVert_p
+
\lVert
    \boldsymbol{c}-\boldsymbol{z}
\rVert_p$
and
$\lVert
    \widehat{\boldsymbol{c}}-\boldsymbol{c}
\rVert_p
=
\widehat{r}-r$.
Because the $\ell_p$-norm is strictly convex for $p\in(1,\infty)$,
equality in the triangle inequality implies that
there exists $\lambda\ge0$ such that
$\widehat{\boldsymbol{c}}-\boldsymbol{c}
=\lambda(\boldsymbol{c}-\boldsymbol{z})$.
Taking norms on both sides and using
$\lVert\boldsymbol{c}-\boldsymbol{z}\rVert_p=r$
gives
$\lambda=(\widehat{r}-r)/r$.
Hence,
$\widehat{\boldsymbol{c}}
=
\boldsymbol{c}
+
\frac{\widehat{r}-r}{r}
(\boldsymbol{c}-\boldsymbol{z})
=
\boldsymbol{z}
+
\frac{\widehat{r}}{r}
(\boldsymbol{c}-\boldsymbol{z})$,
which proves \eqref{eq:containing-ball-center}.
\Halmos
\endproof

We next establish that, for every $p\in[1,\infty)$,
a DD-free $\ell_p$-ball with a unique contact point can be strictly enlarged while preserving DD-freeness.

\begin{lemma}
\label{lemma:singleton-contact-expansion}
    Let $\dd{D}$ be an IDD and consider $p\in[1,\infty)$.
    Let $\mathcal{B}_p(\boldsymbol{c},r)$
    be a DD-free $\ell_p$-ball with respect to
    $\dd{D}$, where $r>0$,
    and suppose that
    $\mathcal{Z}_p(\boldsymbol{c},r)
    =\{\boldsymbol{z}\}$.
    Then there exists $\varepsilon_0>0$ such that,
    for every $\varepsilon\in(0,\varepsilon_0]$,
    $\mathcal{B}_p\bigl(
        \boldsymbol{c}
        +\varepsilon
        (\boldsymbol{c}-\boldsymbol{z})/r,
        r+\varepsilon
    \bigr)$
    is DD-free and strictly contains
    $\mathcal{B}_p(\boldsymbol{c},r)$.
\end{lemma}

\proof{Proof.}

For each $\varepsilon>0$, define
\[
    \boldsymbol{u}
    \coloneqq
    \frac{\boldsymbol{c}-\boldsymbol{z}}{r},
    \qquad
    \boldsymbol{c}(\varepsilon)
    \coloneqq
    \boldsymbol{c}+\varepsilon\boldsymbol{u}
    =
    \boldsymbol{z}+(r+\varepsilon)\boldsymbol{u}.
\]
Since $\lVert\boldsymbol{u}\rVert_p=1$, we have
$\lVert
    \boldsymbol{c}(\varepsilon)-\boldsymbol{c}
\rVert_p
=
\varepsilon$.
Hence, Lemma~\ref{lemma:lp-ball-containment} gives
$\mathcal{B}_p(\boldsymbol{c},r)
\subseteq
\mathcal{B}_p\bigl(
    \boldsymbol{c}(\varepsilon),r+\varepsilon
\bigr)$.
Moreover, the inclusion is strict, since
$\boldsymbol{c}(\varepsilon)+(r+\varepsilon)\boldsymbol{u}$
belongs to
$\mathcal{B}_p\bigl(
    \boldsymbol{c}(\varepsilon),r+\varepsilon
\bigr)$
but not to
$\mathcal{B}_p(\boldsymbol{c},r)$.
It remains to show that
$\mathcal{B}_p\bigl(
    \boldsymbol{c}(\varepsilon),r+\varepsilon
\bigr)$
is DD-free for sufficiently small $\varepsilon>0$.

Consider first any path
$\dd{q}=(\dd{a}_1,\ldots,\dd{a}_n)\in\dd{P}(\dd{t})$
such that
$\boldsymbol{z}\in\mathcal{S}_{\dd{q}}$.
Since
$\boldsymbol{z}\in\mathcal{S}_{\dd{q}}$ and
$\mathcal{B}_p(\boldsymbol{c},r)$ is DD-free, we have
$d_p(\boldsymbol{c},\mathcal{S}_{\dd{q}})=r$.
By Proposition~\ref{prop:ball-contact-subgraph},
$\boldsymbol{z}$ is the coordinatewise projection of $\boldsymbol{c}$ onto $\mathcal{S}_{\dd{q}}$, which implies that each $z_i$ is a closest point in $\dd{l}(\dd{a}_i)$ to $c_i$.
If $c_i\in\dd{l}(\dd{a}_i)$,
then we must have $z_i=c_i$ as $z_i$ is a closest point 
in $\dd{l}(\dd{a}_i)$ to $c_i$.
Hence,
$\bigl(\boldsymbol{c}(\varepsilon)\bigr)_i
=c_i+(\varepsilon/r)(c_i-z_i)
=c_i
=z_i$.
If $c_i\notin\dd{l}(\dd{a}_i)$,
then $z_i$ is the endpoint of $\dd{l}(\dd{a}_i)$ closest to $c_i$.
Since
$\bigl(\boldsymbol{c}(\varepsilon)\bigr)_i
=
z_i+(1+\varepsilon/r)(c_i-z_i)$
with $1+\varepsilon/r>1$,
the point
$\bigl(\boldsymbol{c}(\varepsilon)\bigr)_i$
lies beyond $c_i$ on the ray emanating from $z_i$ through $c_i$.
Thus, in either case, $z_i$ remains a closest point in
$\dd{l}(\dd{a}_i)$ to
$\bigl(\boldsymbol{c}(\varepsilon)\bigr)_i$.
Consequently,
$d_p\bigl(
    \boldsymbol{c}(\varepsilon),\mathcal{S}_{\dd{q}}
\bigr)
=
\lVert
    \boldsymbol{c}(\varepsilon)-\boldsymbol{z}
\rVert_p
=
r+\varepsilon$.
Thus, this equality holds for every path satisfying
$\boldsymbol{z}\in\mathcal{S}_{\dd{q}}$.

Now consider the paths
$\dd{q}\in\dd{P}(\dd{t})$ satisfying
$\boldsymbol{z}\notin\mathcal{S}_{\dd{q}}$.
For every such path,
$d_p(\boldsymbol{c},\mathcal{S}_{\dd{q}})>r$.
Indeed, because 
$\mathcal{B}_p(\boldsymbol{c},r)$ is DD-free,
$d_p(\boldsymbol{c},\mathcal{S}_{\dd{q}})\ge r$.
If equality held, compactness of
$\mathcal{S}_{\dd{q}}$ would yield a point
$\boldsymbol{y}\in\mathcal{S}_{\dd{q}}$ such that
$\lVert\boldsymbol{c}-\boldsymbol{y}\rVert_p=r$.
It would follow that
$\boldsymbol{y}\in
\mathcal{Z}_p(\boldsymbol{c},r)
=
\{\boldsymbol{z}\}$,
contradicting
$\boldsymbol{z}\notin\mathcal{S}_{\dd{q}}$.
Since there are finitely many such paths, we may choose
$\varepsilon_0>0$ sufficiently small that
$d_p(\boldsymbol{c},\mathcal{S}_{\dd{q}})
>
r+2\varepsilon_0$
for all of them.
Fix any $\varepsilon\in(0,\varepsilon_0]$.
For every such path, it follows from the triangle inequality that
$d_p\bigl(
    \boldsymbol{c}(\varepsilon),\mathcal{S}_{\dd{q}}
\bigr)
\ge
d_p(\boldsymbol{c},\mathcal{S}_{\dd{q}})
-
\lVert
    \boldsymbol{c}(\varepsilon)-\boldsymbol{c}
\rVert_p
>
r+\varepsilon$.

Since
$\boldsymbol{z}\in\operatorname{Sol}(\dd{D})$,
at least one path $\dd{q}$ satisfies
$\boldsymbol{z}\in\mathcal{S}_{\dd{q}}$.
Combining the two cases and using
$\operatorname{Sol}(\dd{D})
=
\bigcup_{\dd{q}\in\dd{P}(\dd{t})}
\mathcal{S}_{\dd{q}}$,
we obtain
$d_p\bigl(
    \boldsymbol{c}(\varepsilon),\operatorname{Sol}(\dd{D})
\bigr)
=
r+\varepsilon$.
Therefore, 
for every $\varepsilon\in(0,\varepsilon_0]$, 
the $\ell_p$-ball
$\mathcal{B}_p\bigl(
    \boldsymbol{c}(\varepsilon),r+\varepsilon
\bigr)$
is DD-free and strictly contains
$\mathcal{B}_p(\boldsymbol{c},r)$.
\Halmos
\endproof

Lemmas~\ref{lemma:ball-containment-rigidity} 
and~\ref{lemma:singleton-contact-expansion}, 
together with the fact that
$\operatorname{Sol}(\dd{D})$ is a finite union of hyperrectangles
induced by $\dd{r}$--$\dd{t}$ paths, yields the following
characterization of $\ell_p$-ball-maximality for $p\in(1,\infty)$.

\begin{theorem}[Maximality characterization for DD-free $\ell_p$-balls for $p\in(1,\infty)$]
\label{thm:ball-maximal-characterization}
    Let
    $\dd{D}=(\dd{U},\dd{A},\dd{l}(\cdot))$
    be an IDD with 
    $\operatorname{Sol}(\dd{D})\neq\emptyset$.
    Consider $p\in(1,\infty)$, 
    and let
    $\mathcal{B}_p(\boldsymbol{c},r)$ be a 
    DD-free $\ell_p$-ball 
    with respect to $\dd{D}$, where $r>0$.
    Then $\mathcal{B}_p(\boldsymbol{c},r)$ is 
    $\ell_p$-ball-maximal
    if and only if
    $\mathcal{Z}_p(\boldsymbol{c},r)$ contains 
    at least two distinct points.
\end{theorem}

\proof{Proof.}

First, suppose that
$\boldsymbol{z}_1,\boldsymbol{z}_2
\in\mathcal{Z}_p(\boldsymbol{c},r)$
with $\boldsymbol{z}_1\neq\boldsymbol{z}_2$.
Let
$\mathcal{B}_p(\widehat{\boldsymbol{c}},\widehat{r})$
be any DD-free $\ell_p$-ball containing
$\mathcal{B}_p(\boldsymbol{c},r)$.
Applying Lemma~\ref{lemma:ball-containment-rigidity}
to each contact point gives
$\widehat{\boldsymbol{c}}
=
\boldsymbol{z}_i
+
\frac{\widehat{r}}{r}
(\boldsymbol{c}-\boldsymbol{z}_i)
$
for $i\in\{1,2\}$.
Subtracting these identities yields
$\left(1-\frac{\widehat{r}}{r}\right)
(\boldsymbol{z}_1-\boldsymbol{z}_2)
=
\boldsymbol{0}$.
Since $\boldsymbol{z}_1\neq\boldsymbol{z}_2$,
we obtain $\widehat{r}=r$ and consequently
$\widehat{\boldsymbol{c}}=\boldsymbol{c}$.
Thus, no DD-free $\ell_p$-ball strictly contains
$\mathcal{B}_p(\boldsymbol{c},r)$, proving $\ell_p$-ball-maximality.

Conversely, suppose that
$\mathcal{B}_p(\boldsymbol{c},r)$ is $\ell_p$-ball-maximal.
If $\mathcal{Z}_p(\boldsymbol{c},r)=\emptyset$, then we have
$\mathcal{B}_p(\boldsymbol{c},r)
\cap\operatorname{Sol}(\dd{D})=\emptyset$
because
$\mathcal{B}_p(\boldsymbol{c},r)$ is DD-free.
Since $\operatorname{Sol}(\dd{D})$ is nonempty and compact,
$d_p(\boldsymbol{c},\operatorname{Sol}(\dd{D}))>r$.
Then
$\mathcal{B}_p\bigl(
    \boldsymbol{c},
    d_p(\boldsymbol{c},\operatorname{Sol}(\dd{D}))
\bigr)$
is a DD-free $\ell_p$-ball strictly containing
$\mathcal{B}_p(\boldsymbol{c},r)$, a contradiction.
If $\mathcal{Z}_p(\boldsymbol{c},r)=\{\boldsymbol{z}\}$, then
Lemma~\ref{lemma:singleton-contact-expansion} yields a DD-free
$\ell_p$-ball strictly containing $\mathcal{B}_p(\boldsymbol{c},r)$,
again a contradiction.
Therefore, $\mathcal{Z}_p(\boldsymbol{c},r)$ contains at least two
distinct points.
\Halmos
\endproof

The next remark shows that the conclusion of Lemma~\ref{lemma:ball-containment-rigidity} and the necessary and sufficient conditions for $\ell_p$-ball-maximality established in Theorem~\ref{thm:ball-maximal-characterization} do not extend to the case $p=1$. 
In particular, the existence of two distinct contact points
is no longer sufficient for $\ell_1$-ball-maximality.
As a result, in Theorem~\ref{thm:l1-ball-maximal-characterization}, 
we characterize $\ell_1$-ball-maximality from a different perspective by examining the coordinate bounds of the path-induced hyperrectangles that intersect the boundary of the ball.

\begin{remark}
\label{remark:p-1-limitations}    
    Consider an IDD $\dd{D}$ satisfying
    $\operatorname{Sol}(\dd{D})=\{(-1,0),(1,0)\}$.  
    Let
    $\widehat{\mathcal{B}}
    \coloneqq
    \mathcal{B}_1((0,1),2)$.
    Since $\lVert(0,1)\rVert_1+1=2$,
    Lemma~\ref{lemma:lp-ball-containment} gives
    $\mathcal{B}_1(\boldsymbol{0},1)
    \subseteq\widehat{\mathcal{B}}$.
    The inclusion is strict because
    $(0,3)\in
    \widehat{\mathcal{B}}
    \setminus\mathcal{B}_1(\boldsymbol{0},1)$.
    Moreover,
    $\lVert(\pm1,0)-(0,1)\rVert_1=2$,
    so both points of $\operatorname{Sol}(\dd{D})$ lie on
    $\operatorname{bd}(\widehat{\mathcal{B}})$, and
    $\widehat{\mathcal{B}}$ is DD-free.
    Thus, $\mathcal{B}_1(\boldsymbol{0},1)$ is not
    $\ell_1$-ball-maximal despite having two distinct contact points.
    The center $(0,1)$ of $\widehat{\mathcal{B}}$ does not lie on either of the rays emanating from the two contact points of $\mathcal{B}_1(\boldsymbol{0},1)$ through $\boldsymbol{0}$, and hence does not satisfy \eqref{eq:containing-ball-center}.
    Thus,
    Lemma~\ref{lemma:ball-containment-rigidity} does not extend to
    $p=1$.
    For the unit $\ell_1$-ball
    $\mathcal{B}_1(\boldsymbol{0},1)$ in this example, the radial
    expansion induced by $\boldsymbol{z}=(1,0)$ consists of the $\ell_1$-balls
    $\mathcal{B}_1((1-\rho,0),\rho)$ with $\rho\ge1$.
    These balls contain $\mathcal{B}_1(\boldsymbol{0},1)$ by
    Lemma~\ref{lemma:lp-ball-containment}.
    For every $\rho>1$, however,
    $\lVert(-1,0)-(1-\rho,0)\rVert_1=|\rho-2|<\rho$,
    so the other contact point $(-1,0)$ lies in the interior of the expanded ball.
    The argument for the contact point $(-1,0)$ is symmetric.
    Thus, for either choice of contact point, the largest DD-free member of the corresponding radial family is $\mathcal{B}_1(\boldsymbol{0},1)$, which is not $\ell_1$-ball-maximal.
\end{remark}


\begin{theorem}[Maximality characterization for DD-free $\ell_1$-balls]
\label{thm:l1-ball-maximal-characterization}
    Let $\dd{D}=(\dd{U},\dd{A},\dd{l}(\cdot))$ be an IDD with
    $\operatorname{Sol}(\dd{D})\neq\emptyset$.
    Let $\mathcal{B}_1(\boldsymbol{c},r)$ be a DD-free
    $\ell_1$-ball with respect to $\dd{D}$, where $r>0$.
    For each $\dd{p}\in\dd{P}(\dd{t})$, write
    $\mathcal{S}_{\dd{p}}
    =\prod_{i\in[n]}[l_i^{\dd{p}},u_i^{\dd{p}}]$.
    Let 
    $\dd{P}^{\ast}
    \coloneqq
    \{
        \dd{p}\in\dd{P}(\dd{t}):
        d_1(\boldsymbol{c},\mathcal{S}_{\dd{p}})=r
    \}$
    denote the set of \emph{active paths}.
    Then, $\mathcal{B}_1(\boldsymbol{c},r)$ is
    $\ell_1$-ball-maximal if and only if
    $\dd{P}^{\ast}\neq\emptyset$ and
    \begin{equation}
    \label{eq:l1-ball-maximal-coordinate-bounds}
        \min_{\dd{p}\in\dd{P}^{\ast}} l_i^{\dd{p}}
        < c_i <
        \max_{\dd{p}\in\dd{P}^{\ast}} u_i^{\dd{p}},
        \qquad \forall i\in[n].
    \end{equation}
\end{theorem}

\proof{Proof.}

First, suppose that $\dd{P}^{\ast}\neq\emptyset$ and
\eqref{eq:l1-ball-maximal-coordinate-bounds} holds.
Let $\mathcal{B}_1(\widehat{\boldsymbol{c}},\widehat{r})$
be any DD-free $\ell_1$-ball containing
$\mathcal{B}_1(\boldsymbol{c},r)$, and set
$\boldsymbol{h}
\coloneqq
\widehat{\boldsymbol{c}}-\boldsymbol{c}$.
By Lemma~\ref{lemma:lp-ball-containment}, we have
$\widehat{r}\ge r+\lVert\boldsymbol{h}\rVert_1$.
For every $\dd{p}\in\dd{P}^{\ast}$, we have
\[
    \widehat{r}
    \le
    d_1(\widehat{\boldsymbol{c}},\mathcal{S}_{\dd{p}})
    =
    \min_{\boldsymbol{z}\in\mathcal{S}_{\dd{p}}}
    \lVert\boldsymbol{c}+\boldsymbol{h}-\boldsymbol{z}\rVert_1
    \le
    \min_{\boldsymbol{z}\in\mathcal{S}_{\dd{p}}}
    \bigl(
        \lVert\boldsymbol{c}-\boldsymbol{z}\rVert_1
        +\lVert\boldsymbol{h}\rVert_1
    \bigr)
    =
    r+\lVert\boldsymbol{h}\rVert_1
    \le
    \widehat{r},
\]
where the first inequality holds because
$\mathcal{B}_1(\widehat{\boldsymbol{c}},\widehat{r})$ is DD-free,
the second inequality follows from the triangle inequality, the subsequent equality is implied by the definition of
$\dd{P}^{\ast}$, and the last inequality is given by
Lemma~\ref{lemma:lp-ball-containment}.
Since the chain begins and ends at $\widehat{r}$, all inequalities must hold with equality.
In particular, we have
$d_1(\widehat{\boldsymbol{c}},\mathcal{S}_{\dd{p}})
=r+\lVert\boldsymbol{h}\rVert_1$
and
$\widehat{r}=r+\lVert\boldsymbol{h}\rVert_1$.
We next show that $\boldsymbol{h}=\boldsymbol{0}$.
Fix $i\in[n]$ and suppose, to the contrary, that $h_i>0$.
By \eqref{eq:l1-ball-maximal-coordinate-bounds}, there exists
an active path $\dd{p}\in\dd{P}^{\ast}$
such that $c_i<u^\dd{p}_i$.
It follows that
\[
    d_1(c_i+h_i,[l^\dd{p}_i,u^\dd{p}_i])
    <
    d_1(c_i,[l^\dd{p}_i,u^\dd{p}_i])+h_i.
\]
Indeed, if $c_i<l^\dd{p}_i$, the inequality follows from
$|c_i+h_i-l^\dd{p}_i|<l^\dd{p}_i-c_i+h_i$.
Otherwise, $l^\dd{p}_i\le c_i<u^\dd{p}_i$, so
$d_1(c_i,[l^\dd{p}_i,u^\dd{p}_i])=0$ and
$d_1(c_i+h_i,[l^\dd{p}_i,u^\dd{p}_i])
=\max\{c_i+h_i-u^\dd{p}_i,0\}<h_i$.
For every $j\in[n]\setminus\{i\}$, the triangle inequality gives
$d_1(c_j+h_j,[l^\dd{p}_j,u^\dd{p}_j])
\le d_1(c_j,[l^\dd{p}_j,u^\dd{p}_j])+|h_j|$.
Summing these inequalities and using the coordinatewise distance decomposition in the proof of
Proposition~\ref{prop:distance-to-path} yields
$d_1(\widehat{\boldsymbol{c}},\mathcal{S}_{\dd{p}})
<r+\lVert\boldsymbol{h}\rVert_1$,
contradicting the equality
$d_1(\widehat{\boldsymbol{c}},\mathcal{S}_{\dd{p}})
=r+\lVert\boldsymbol{h}\rVert_1$
established above.
Hence, $h_i\le0$.
Similarly, using an active path
$\dd{p}\in\dd{P}^{\ast}$ such that $l^\dd{p}_i<c_i$
rules out $h_i<0$.
Thus, $h_i=0$.
Since $i$ was arbitrary, $\boldsymbol{h}=\boldsymbol{0}$.
Substituting this into
$\widehat{r}=r+\lVert\boldsymbol{h}\rVert_1$
gives $\widehat{r}=r$.
Thus, no DD-free $\ell_1$-ball strictly contains
$\mathcal{B}_1(\boldsymbol{c},r)$, proving
$\ell_1$-ball-maximality.

Conversely, suppose that $\mathcal{B}_1(\boldsymbol{c},r)$
is $\ell_1$-ball-maximal.
If $\dd{P}^{\ast}=\emptyset$, then DD-freeness and compactness
of $\operatorname{Sol}(\dd{D})$ imply
$d_1(\boldsymbol{c},\operatorname{Sol}(\dd{D}))>r$.
Consequently,
$\mathcal{B}_1\bigl(
    \boldsymbol{c},
    d_1(\boldsymbol{c},\operatorname{Sol}(\dd{D}))
\bigr)$
is a DD-free $\ell_1$-ball strictly containing
$\mathcal{B}_1(\boldsymbol{c},r)$, a contradiction.
Therefore, $\dd{P}^{\ast}\neq\emptyset$.

It remains to establish
\eqref{eq:l1-ball-maximal-coordinate-bounds}.
We consider two cases. For the first case, suppose by contradiction that
$u^\dd{p}_i\le c_i$ for some $i\in[n]$ and every
$\dd{p}\in\dd{P}^{\ast}$.
For $\varepsilon>0$ to be chosen below, set
$\boldsymbol{c}(\varepsilon)
\coloneqq
\boldsymbol{c}+\varepsilon\boldsymbol{e}_i$,
where $\boldsymbol{e}_i$ is the $i$-th unit vector.
For every $\dd{p}\in\dd{P}^{\ast}$, since
$u^\dd{p}_i\le c_i$, shifting the center 
from $\boldsymbol{c}$ to $\boldsymbol{c}(\varepsilon)$
increases the distance 
to $[l^\dd{p}_i,u^\dd{p}_i]$ by exactly $\varepsilon$, while leaving the distances to the other coordinate intervals unchanged.
Thus,
$d_1(
    \boldsymbol{c}(\varepsilon),
    \mathcal{S}_{\dd{p}}
)
=r+\varepsilon$.
For every 
$\dd{q}
\in
\dd{P}(\dd{t})\setminus\dd{P}^{\ast}$,
DD-freeness and the definition of $\dd{P}^{\ast}$ give
$d_1(\boldsymbol{c},\mathcal{S}_{\dd{q}})>r$.
Since there are finitely many such paths, we may choose
$\varepsilon>0$ sufficiently small that
$d_1(\boldsymbol{c},\mathcal{S}_{\dd{q}})>r+2\varepsilon$
for all of them.
Thus, for every such path, we have
\[
d_1(\boldsymbol{c}(\varepsilon),\mathcal{S}_{\dd{q}})
\ge
d_1(\boldsymbol{c},\mathcal{S}_{\dd{q}})-
\lVert\boldsymbol{c}-\boldsymbol{c}(\varepsilon)\rVert_1
=
d_1(\boldsymbol{c},\mathcal{S}_{\dd{q}})-
\varepsilon
>
r+\varepsilon,
\]
where the first inequality follows from the triangle inequality, and the last inequality is given by $d_1(\boldsymbol{c},\mathcal{S}_{\dd{q}})>r+2\varepsilon$.
Combining the above two results and using
$\dd{P}^{\ast}\neq\emptyset$, we obtain
$d_1(\boldsymbol{c}(\varepsilon),
\operatorname{Sol}(\dd{D}))=r+\varepsilon$.
For this choice of $\varepsilon$, the ball
$\mathcal{B}_1(\boldsymbol{c}(\varepsilon),r+\varepsilon)$
is DD-free and contains $\mathcal{B}_1(\boldsymbol{c},r)$
by Lemma~\ref{lemma:lp-ball-containment}.
The inclusion is strict, since
$\boldsymbol{c}+(r+2\varepsilon)\boldsymbol{e}_i$
belongs to 
$\mathcal{B}_1(\boldsymbol{c}(\varepsilon),r+\varepsilon)$ 
but not to $\mathcal{B}_1(\boldsymbol{c},r)$.
This contradicts the assumed $\ell_1$-ball-maximality of
$\mathcal{B}_1(\boldsymbol{c},r)$.
The second case in which $l^\dd{p}_i\ge c_i$ for every
$\dd{p}\in\dd{P}^{\ast}$ is analogous, using the direction
$-\boldsymbol{e}_i$.
Therefore, \eqref{eq:l1-ball-maximal-coordinate-bounds}
holds for every $i\in[n]$.
\Halmos
\endproof

Geometrically, the active path-induced hyperrectangles must collectively extend to both sides of the $\ell_1$-ball center in each coordinate, thereby precluding any strictly containing DD-free $\ell_1$-ball.
We show in Corollary~\ref{cor:radial-ball-maximality-test} of the next section that verifying $\ell_p$-ball maximality for $p\in(0,\infty)$ as characterized in Theorem~\ref{thm:ball-maximal-characterization}, as well as $\ell_1$-ball maximality in Theorem~\ref{thm:l1-ball-maximal-characterization}, can be performed efficiently in time linear in the size of the IDD by exploiting the compact subgraph representation of contact points established in Proposition~\ref{prop:ball-contact-subgraph}.

The following remark discusses the difference between $\ell_p$-ball-maximality, as characterized in Theorems~\ref{thm:ball-maximal-characterization} and~\ref{thm:l1-ball-maximal-characterization}, and box-maximality, as described in Proposition~\ref{prop:box-maximal-characterization}, with respect to obtaining maximal DD-free sets.

\begin{remark}
\label{remark:lp-ball-maximality-limitations}
    While Proposition~\ref{prop:box-maximal-characterization} shows that box-maximality coincides with maximality for full-dimensional DD-free boxes, $\ell_p$-ball-maximality need not imply maximality among all DD-free sets. In particular, Proposition~\ref{prop:polyhedral-maximal-DD-free} implies that $\ell_p$-ball-maximal sets for $p\in(0,\infty)$ are not maximal among all DD-free sets. Furthermore, since $\ell_1$-ball-maximal sets are polyhedra, Proposition~\ref{prop:polyhedral-maximal-characterization} provides a necessary and sufficient condition for determining their maximality among all DD-free sets.
    For example, consider an IDD $\dd{D}$ satisfying
    $\operatorname{Sol}(\dd{D})=\{(-1,0),(1,0),(0,-1),(0,1)\}$.
    For every $p\in(1,\infty)$, the unit ball
    $\mathcal{B}_p(\boldsymbol{0},1)$ is $\ell_p$-ball-maximal by Theorem~\ref{thm:ball-maximal-characterization}, since all four points are contact points.
    For $p=1$, all four points are contact points, and in each coordinate, both $-1$ and $1$ are attained among these points.
    Thus, $\mathcal{B}_1(\boldsymbol{0},1)$ is $\ell_1$-ball-maximal by Theorem~\ref{thm:l1-ball-maximal-characterization}.
    However, for every $p\in[1,\infty)$, the unit ball is strictly contained in the DD-free box $[-1,1]^2$.
\end{remark}

For every $p\in[1,\infty)$, the initial $\ell_p$-ball
$\mathcal{C}_p(\bar{\boldsymbol{x}})$ has at least one contact
point because $\operatorname{Sol}(\dd{D})$ is nonempty and compact.
We next use a selected contact point to define the radial expansion and characterize the radius range over which the resulting $\ell_p$-balls remain DD-free.

\subsubsection{Radial $\ell_p$-Ball Expansion}
\label{subsubsec:ball-radial-expansion}


Consider an IDD $\dd{D}$ with
$\operatorname{Sol}(\dd{D})\neq\emptyset$
and a point 
$\bar{\boldsymbol{x}} 
\notin 
\operatorname{Sol}(\dd{D})$. 
Fix $p\in[1,\infty)$, and set
$r_0\coloneqq
d_p(\bar{\boldsymbol{x}},\operatorname{Sol}(\dd{D}))>0$.
Choose a contact point
$\boldsymbol{z}\in
\mathcal{Z}_p(\bar{\boldsymbol{x}},r_0)$
using the procedure in 
Proposition~\ref{prop:ball-contact-subgraph}
and define the unit vector
$\boldsymbol{u}
\coloneqq
(\bar{\boldsymbol{x}}-\boldsymbol{z})/r_0$.
For every $r\in[r_0,\infty)$, define
\begin{equation}
\label{eq:radial-center}
    \boldsymbol{c}_{\boldsymbol{z}}^{\boldsymbol{u}}(r)
    \coloneqq
    \boldsymbol{z}+r\boldsymbol{u}.
\end{equation}
Thus, as $r$ increases from $r_0$,
$\boldsymbol{c}_{\boldsymbol{z}}^{\boldsymbol{u}}(r)$ moves along
the ray from $\boldsymbol{z}$ through $\bar{\boldsymbol{x}}$, with
$\boldsymbol{c}_{\boldsymbol{z}}^{\boldsymbol{u}}(r_0)
=
\bar{\boldsymbol{x}}$.
We refer to the family of balls
$\mathcal{B}_p\bigl(
    \boldsymbol{c}_{\boldsymbol{z}}^{\boldsymbol{u}}(r),r
\bigr)$,
for all $r\in[r_0,\infty)$, as the
\emph{radial family induced by $\boldsymbol{z}$}, which is denoted by $\bigl<\mathcal{B}_p\bigl(
    \boldsymbol{c}_{\boldsymbol{z}}^{\boldsymbol{u}}(r),r
\bigr)\bigr>_{r=r_0}^{\infty}$.

\begin{lemma}
\label{lemma:radial-family-geometry}

Consider the radial family $\bigl<\mathcal{B}_p\bigl(
    \boldsymbol{c}_{\boldsymbol{z}}^{\boldsymbol{u}}(r),r
\bigr)\bigr>_{r=r_0}^{\infty}$ induced by $\boldsymbol{z}$.
Then, $\boldsymbol{z}$ lies on the boundary of every member of this family and, whenever that member is DD-free,
$\boldsymbol{z}
\in
\mathcal{Z}_p\bigl(
    \boldsymbol{c}_{\boldsymbol{z}}^{\boldsymbol{u}}(r),r
\bigr)$.
Moreover, if $r_0\le r\le r'$, then
$\mathcal{B}_p\bigl(
    \boldsymbol{c}_{\boldsymbol{z}}^{\boldsymbol{u}}(r),r
\bigr)
\subseteq
\mathcal{B}_p\bigl(
    \boldsymbol{c}_{\boldsymbol{z}}^{\boldsymbol{u}}(r'),r'
\bigr)$,
with strict inclusion whenever $r<r'$.
\end{lemma}

\proof{Proof.}

Since
$\boldsymbol{z}
\in
\mathcal{Z}_p(\bar{\boldsymbol{x}},r_0)$,
we have
$\boldsymbol{z}\in\operatorname{Sol}(\dd{D})$ and
$\lVert\bar{\boldsymbol{x}}-\boldsymbol{z}\rVert_p=r_0$.
Hence, by the definition of $\boldsymbol{u}$, we have
$\lVert\boldsymbol{u}\rVert_p=1$.
If follows that
$\lVert
    \boldsymbol{c}_{\boldsymbol{z}}^{\boldsymbol{u}}(r)
    -\boldsymbol{z}
\rVert_p
=r\lVert\boldsymbol{u}\rVert_p
=r$
for every $r\in[r_0,\infty)$.
Thus,
$\boldsymbol{z}
\in
\operatorname{bd}\bigl(
    \mathcal{B}_p(
        \boldsymbol{c}_{\boldsymbol{z}}^{\boldsymbol{u}}(r),r
    )
\bigr)$.
Whenever this ball is DD-free, the definition of the contact set therefore gives
$\boldsymbol{z}
\in
\mathcal{Z}_p\bigl(
    \boldsymbol{c}_{\boldsymbol{z}}^{\boldsymbol{u}}(r),r
\bigr)$.

Now let $r_0\le r\le r'$.
By construction, we have
$\lVert
    \boldsymbol{c}_{\boldsymbol{z}}^{\boldsymbol{u}}(r')
    -
    \boldsymbol{c}_{\boldsymbol{z}}^{\boldsymbol{u}}(r)
\rVert_p
=(r'-r)\lVert\boldsymbol{u}\rVert_p
=r'-r$.
Lemma~\ref{lemma:lp-ball-containment} then yields
$\mathcal{B}_p\bigl(
    \boldsymbol{c}_{\boldsymbol{z}}^{\boldsymbol{u}}(r),r
\bigr)
\subseteq
\mathcal{B}_p\bigl(
    \boldsymbol{c}_{\boldsymbol{z}}^{\boldsymbol{u}}(r'),r'
\bigr)$.
If $r<r'$, the point
$\boldsymbol{y}
\coloneqq
\boldsymbol{c}_{\boldsymbol{z}}^{\boldsymbol{u}}(r')
+r'\boldsymbol{u}$
satisfies
$\lVert
    \boldsymbol{y}
    -
    \boldsymbol{c}_{\boldsymbol{z}}^{\boldsymbol{u}}(r')
\rVert_p=r'$
and
$\lVert
    \boldsymbol{y}
    -
    \boldsymbol{c}_{\boldsymbol{z}}^{\boldsymbol{u}}(r)
\rVert_p=2r'-r>r$.
Hence,
$\boldsymbol{y}
\in
\mathcal{B}_p\bigl(
    \boldsymbol{c}_{\boldsymbol{z}}^{\boldsymbol{u}}(r'),r'
\bigr)
\setminus
\mathcal{B}_p\bigl(
    \boldsymbol{c}_{\boldsymbol{z}}^{\boldsymbol{u}}(r),r
\bigr)$,
which proves strict inclusion.
\Halmos
\endproof

Due to the nestedness property of the radial family, as established in Lemma~\ref{lemma:radial-family-geometry}, 
we can define the \emph{radial expansion threshold} by
\begin{equation}
\label{eq:radial-expansion-threshold}
    r^\ast
    \coloneqq
    \sup\left\{
        r\in[r_0,\infty):
        \mathcal{B}_p\bigl(
            \boldsymbol{c}_{\boldsymbol{z}}^{\boldsymbol{u}}(r),r
        \bigr)
        \text{ is DD-free}
    \right\}.
\end{equation}

\begin{theorem}[Radial expansion]
\label{thm:normal-ray-ball-expansion}
    Let
    $\dd{D}=(\dd{U},\dd{A},\dd{l}(\cdot))$
    be an IDD with
    $\operatorname{Sol}(\dd{D})\neq\emptyset$,
    let
    $\bar{\boldsymbol{x}}
    \notin
    \operatorname{Sol}(\dd{D})$,
    and fix $p\in[1,\infty)$.
    Set
    $r_0
    \coloneqq
    d_p\bigl(
        \bar{\boldsymbol{x}},
        \operatorname{Sol}(\dd{D})
    \bigr)>0$,
    choose
    $\boldsymbol{z}
    \in
    \mathcal{Z}_p(\bar{\boldsymbol{x}},r_0)$,
    and let
    $\boldsymbol{u}
    \coloneqq
    (\bar{\boldsymbol{x}}-\boldsymbol{z})/r_0$.
    Let
    $\boldsymbol{c}_{\boldsymbol{z}}^{\boldsymbol{u}}(\cdot)$
    and $r^\ast$ be defined by
    \eqref{eq:radial-center} and
    \eqref{eq:radial-expansion-threshold}, respectively.
    Then $r^\ast\in[r_0,+\infty]$, and the following statements hold.
    \begin{enumerate}[label=\textup{(\roman*)}, ref=\textup{(\roman*)}]
        \item
        For every $r\in[r_0,\infty)$,
        $\mathcal{B}_p\bigl(
            \boldsymbol{c}_{\boldsymbol{z}}^{\boldsymbol{u}}(r),r
        \bigr)$
        is DD-free if and only if $r\le r^\ast$.
        \item
        If $r^\ast<+\infty$ and $p\in(1,\infty)$, then
        $\mathcal{B}_p\bigl(
            \boldsymbol{c}_{\boldsymbol{z}}^{\boldsymbol{u}}(r^\ast),
            r^\ast
        \bigr)$ is the unique
        $\ell_p$-ball-maximal DD-free ball containing
        $\mathcal{B}_p(\bar{\boldsymbol{x}},r_0)$.
        Moreover,
        $\mathcal{Z}_p\bigl(
            \boldsymbol{c}_{\boldsymbol{z}}^{\boldsymbol{u}}(r^\ast),
            r^\ast
        \bigr)$
        contains a point distinct from $\boldsymbol{z}$.
        \item If $r^\ast<+\infty$ and $p=1$, then
        $\mathcal{B}_1\bigl(
            \boldsymbol{c}_{\boldsymbol{z}}^{\boldsymbol{u}}(r^\ast),
            r^\ast
        \bigr)$ is the largest DD-free member of the radial family
        $\bigl<\mathcal{B}_p\bigl(
            \boldsymbol{c}_{\boldsymbol{z}}^{\boldsymbol{u}}(r),r
        \bigr)\bigr>_{r=r_0}^{\infty}$. 
        Furthermore, it is $\ell_1$-ball-maximal if and only if $\dd{P}^{\ast}\neq\emptyset$ and \eqref{eq:l1-ball-maximal-coordinate-bounds} holds, where $\dd{P}^{\ast}$ is defined as in Theorem~\ref{thm:l1-ball-maximal-characterization}.
        In addition,
        $\mathcal{Z}_1\bigl(
            \boldsymbol{c}_{\boldsymbol{z}}^{\boldsymbol{u}}(r^\ast),
            r^\ast
        \bigr)$
        contains a point distinct from $\boldsymbol{z}$.
        
    \end{enumerate}
\end{theorem}

\proof{Proof.}

We first prove~\textup{(i)}.
For $r\in[r_0,\infty)$, define
\[
    \psi(r)
    \coloneqq
    d_p(
        \boldsymbol{c}_{\boldsymbol{z}}^{\boldsymbol{u}}(r),
        \operatorname{Sol}(\dd{D})
    )-r.
\]
Since
$\boldsymbol{z}\in\operatorname{Sol}(\dd{D})$
and, by Lemma~\ref{lemma:radial-family-geometry},
$\lVert
\boldsymbol{c}_{\boldsymbol{z}}^{\boldsymbol{u}}(r)
-\boldsymbol{z}
\rVert_p=r$,
we must have $d_p(\boldsymbol{c}_{\boldsymbol{z}}^{\boldsymbol{u}}(r), \operatorname{Sol}(\dd{D})) \leq r$, which yields $\psi(r)\le0$.
By Definition~\ref{def:DD-free-sets},
the $\ell_p$-ball
$\mathcal{B}_p\bigl(
    \boldsymbol{c}_{\boldsymbol{z}}^{\boldsymbol{u}}(r),r
\bigr)$ is DD-free if and only if
$d_p\bigl(
    \boldsymbol{c}_{\boldsymbol{z}}^{\boldsymbol{u}}(r),
    \operatorname{Sol}(\dd{D})
\bigr)\ge r$.
Together with $\psi(r)\le0$, this condition is equivalent to $\psi(r)=0$.
Because $d_p(\cdot,\operatorname{Sol}(\dd{D}))$ is a distance function, it is $1$-Lipschitz.
Therefore, $\psi$ is continuous as $\boldsymbol{c}_{\boldsymbol{z}}^{\boldsymbol{u}}(r)$
depends affinely on $r$.
Moreover, for $r_0\le r\le r'$, it follows from the triangle inequality that
\[
    d_p\bigl(
        \boldsymbol{c}_{\boldsymbol{z}}^{\boldsymbol{u}}(r'),
        \operatorname{Sol}(\dd{D})
    \bigr) -
    d_p\bigl(
        \boldsymbol{c}_{\boldsymbol{z}}^{\boldsymbol{u}}(r),
        \operatorname{Sol}(\dd{D})
    \bigr)
    \le    
    r'-r,
\]
and hence $\psi(r')\le \psi(r)$.
Since
$\boldsymbol{c}_{\boldsymbol{z}}^{\boldsymbol{u}}(r_0)
=
\bar{\boldsymbol{x}}$,
the definition of $r_0$ gives $\psi(r_0)=0$.
Since $\psi(r_0)=0$ and $\psi$ is nonincreasing, the set $\mathcal{I} \coloneqq \{r\in[r_0,\infty):\psi(r)=0\}$
is a nonempty, possibly unbounded interval with left endpoint $r_0$.
Continuity of $\psi$ implies that $\mathcal{I}$ is closed.
Since 
$\mathcal{B}_p\bigl(
    \boldsymbol{c}_{\boldsymbol{z}}^{\boldsymbol{u}}(r),r
\bigr)$ is DD-free if and only if $\psi(r)=0$, it follows from
\eqref{eq:radial-expansion-threshold} that $\sup_{r \in \mathcal{I}} \{r \} = r^\ast$.
Consequently, $\mathcal{I} = [r_0,r^\ast]$ when
$r^\ast<+\infty$ and $\mathcal{I} = [r_0,\infty)$ when
$r^\ast=+\infty$.
Thus,
$\mathcal{B}_p\bigl(
    \boldsymbol{c}_{\boldsymbol{z}}^{\boldsymbol{u}}(r),r
\bigr)$ is DD-free if and only if
$r\le r^\ast$, proving~\textup{(i)}.

We now prove parts~\textup{(ii)} and~\textup{(iii)}.
Suppose that $r^\ast<+\infty$.
By part~\textup{(i)},
$\mathcal{B}_p\bigl(
    \boldsymbol{c}_{\boldsymbol{z}}^{\boldsymbol{u}}(r^\ast),
    r^\ast
\bigr)$
is DD-free, and every DD-free member of the radial family
$\bigl<\mathcal{B}_p\bigl(
    \boldsymbol{c}_{\boldsymbol{z}}^{\boldsymbol{u}}(r),r
\bigr)\bigr>_{r=r_0}^{\infty}$
has radius at most $r^\ast$.
By Lemma~\ref{lemma:radial-family-geometry},
the radial family is nested with respect to the radius.
Hence,
$\mathcal{B}_p\bigl(
    \boldsymbol{c}_{\boldsymbol{z}}^{\boldsymbol{u}}(r^\ast),
    r^\ast
\bigr)$
is the largest DD-free member of the radial family.
Lemma~\ref{lemma:radial-family-geometry} also gives
$\boldsymbol{z}\in
\mathcal{Z}_p\bigl(
    \boldsymbol{c}_{\boldsymbol{z}}^{\boldsymbol{u}}(r^\ast),
    r^\ast
\bigr)$.

We first show that, for $p\in[1,\infty)$,
$\mathcal{Z}_p\bigl(
    \boldsymbol{c}_{\boldsymbol{z}}^{\boldsymbol{u}}(r^\ast),
    r^\ast
\bigr)$
contains a point distinct from $\boldsymbol{z}$.
Suppose, to the contrary, that
$\mathcal{Z}_p\bigl(
    \boldsymbol{c}_{\boldsymbol{z}}^{\boldsymbol{u}}(r^\ast),
    r^\ast
\bigr)=\{\boldsymbol{z}\}$.
For every $\varepsilon>0$, we have
$\boldsymbol{c}_{\boldsymbol{z}}^{\boldsymbol{u}}(r^\ast)
+\varepsilon
\bigl(
    \boldsymbol{c}_{\boldsymbol{z}}^{\boldsymbol{u}}(r^\ast)
    -\boldsymbol{z}
\bigr)/r^\ast
=
\boldsymbol{c}_{\boldsymbol{z}}^{\boldsymbol{u}}(r^\ast+\varepsilon)$.
Thus, Lemma~\ref{lemma:singleton-contact-expansion}
implies that
$\mathcal{B}_p\bigl(
    \boldsymbol{c}_{\boldsymbol{z}}^{\boldsymbol{u}}(r^\ast+\varepsilon),
    r^\ast+\varepsilon
\bigr)$
is DD-free for all sufficiently small $\varepsilon>0$,
contradicting part~\textup{(i)}.
Thus, the contact set contains a point distinct from
$\boldsymbol{z}$.

We next prove the maximaility property in part \textup{(ii)}.
Suppose that $p\in(1,\infty)$.
Let
$\widehat{\mathcal{B}}
\coloneqq
\mathcal{B}_p(
    \widehat{\boldsymbol{c}},
    \widehat{r}
)$
be any DD-free $\ell_p$-ball containing
$\mathcal{B}_p(\bar{\boldsymbol{x}},r_0)$.
Applying Lemma~\ref{lemma:ball-containment-rigidity} to
$\mathcal{B}_p(\bar{\boldsymbol{x}},r_0)$ and its contact point
$\boldsymbol{z}$ gives
$\widehat{r}\ge r_0$ and
$\widehat{\boldsymbol{c}}
=
\boldsymbol{z}
+
\frac{\widehat{r}}{r_0}
(\bar{\boldsymbol{x}}-\boldsymbol{z})
=
\boldsymbol{c}_{\boldsymbol{z}}^{\boldsymbol{u}}
(\widehat{r})$.
Thus, $\widehat{\mathcal{B}}$ belongs to the radial family
$\bigl<\mathcal{B}_p\bigl(
    \boldsymbol{c}_{\boldsymbol{z}}^{\boldsymbol{u}}(r),r
\bigr)\bigr>_{r=r_0}^{\infty}$.
Since $\widehat{\mathcal{B}}$ is DD-free,
part~\textup{(i)} gives
$\widehat{r}\le r^\ast$, and
Lemma~\ref{lemma:radial-family-geometry} yields
$\widehat{\mathcal{B}}
\subseteq
\mathcal{B}_p\bigl(
    \boldsymbol{c}_{\boldsymbol{z}}^{\boldsymbol{u}}(r^\ast),
    r^\ast
\bigr)$.
Therefore,
$\mathcal{B}_p\bigl(
    \boldsymbol{c}_{\boldsymbol{z}}^{\boldsymbol{u}}(r^\ast),
    r^\ast
\bigr)$
is the largest, with respect to set inclusion, among all
DD-free $\ell_p$-balls containing
$\mathcal{B}_p(\bar{\boldsymbol{x}},r_0)$.
Consequently, It is the unique $\ell_p$-ball-maximal DD-free ball
containing $\mathcal{B}_p(\bar{\boldsymbol{x}},r_0)$.
This proves~\textup{(ii)}.

Finally, we prove the maximality property in part \textup{(iii)}. Suppose that $p=1$.
Applying Theorem~\ref{thm:l1-ball-maximal-characterization}
to
$\mathcal{B}_1\bigl(
    \boldsymbol{c}_{\boldsymbol{z}}^{\boldsymbol{u}}(r^\ast),
    r^\ast
\bigr)$
implies the stated necessary and sufficient conditions
for $\ell_1$-ball-maximality.
Together with the conclusions established above,
this proves~\textup{(iii)}.
\Halmos
\endproof

When $r^\ast<+\infty$, Theorem~\ref{thm:normal-ray-ball-expansion} identifies the largest DD-free member of the radial family.
For $p\in(1,\infty)$, this ball is also $\ell_p$-ball-maximal.
For $p=1$, whether this ball is $\ell_1$-ball-maximal can be determined by examining the positions of the path-induced hyperrectangles that intersect the boundary of the ball.
We next characterize the case $r^\ast=+\infty$ in Propositions~\ref{prop:radial-unbounded-smooth} 
and~\ref{prop:radial-unbounded-l1}.
For $p\in(1,\infty)$, we show that the entire radial family is contained in a DD-free halfspace that is maximal among all DD-free sets.
For $p=1$, we identify a DD-free translated polyhedral cone containing the radial family and characterize when this cone is maximal.

\begin{proposition}
\label{prop:radial-unbounded-smooth}
Consider the setting of
    Theorem~\ref{thm:normal-ray-ball-expansion}. Suppose additionally
    that $p\in(1,\infty)$.
    Let
    $\boldsymbol{a}
    \coloneqq
    \bigl(
        \operatorname{sgn}(u_i)|u_i|^{p-1}
    \bigr)_{i\in[n]}$,
    where $\operatorname{sgn}(0)=0$.
    Define
    \begin{equation}
    \label{eq:normal-ray-limit-halfspace}
        \mathcal{H}_p
        \coloneqq
        \left\{
            \boldsymbol{x}\in\mathbb{R}^n:
            \boldsymbol{a}^{\top}
            (\boldsymbol{x}-\boldsymbol{z})\ge0
        \right\},
    \end{equation}
    and set
    \[
        \mu
        \coloneqq
        \max_{\boldsymbol{s}\in\operatorname{Sol}(\dd{D})}
        \boldsymbol{a}^{\top}
        (\boldsymbol{s}-\boldsymbol{z}).
    \]
    Then $\mu\ge0$, and the following statements hold.
    \begin{enumerate}[label=\textup{(\roman*)}, ref=\textup{(\roman*)}]
        \item For every $r\in[r_0,\infty)$,
        $\mathcal{B}_p\bigl(
            \boldsymbol{c}_{\boldsymbol{z}}^{\boldsymbol{u}}(r),r
        \bigr)
        \subseteq
        \mathcal{H}_p$.
        \item The following statements are equivalent:
        $r^\ast=+\infty$; 
        $\mu=0$; 
        and $\mathcal{H}_p$ is DD-free.
        \item If the equivalent statements in part~\textup{(ii)} hold, then $\mathcal{H}_p$ is maximal among all DD-free sets.
    \end{enumerate}
\end{proposition}

\proof{Proof.}

Since $\operatorname{Sol}(\dd{D})$ is nonempty and compact, the
continuous function
$\boldsymbol{s}\mapsto
\boldsymbol{a}^{\top}(\boldsymbol{s}-\boldsymbol{z})$
attains its maximum over $\operatorname{Sol}(\dd{D})$.
Moreover, since
$\boldsymbol{z}\in\operatorname{Sol}(\dd{D})$, we have
$\mu\ge
\boldsymbol{a}^{\top}(\boldsymbol{z}-\boldsymbol{z})=0$.
By construction,
$\lVert\boldsymbol{u}\rVert_p=1$.
Since $p\in(1,\infty)$, the $\ell_p$-norm is differentiable at
$\boldsymbol{u}\neq\boldsymbol{0}$, and its gradient at
$\boldsymbol{u}$ is
\[
    \frac{
        \left(
            \operatorname{sgn}(u_i)|u_i|^{p-1}
        \right)_{i\in[n]}
    }{
        \lVert\boldsymbol{u}\rVert_p^{p-1}
    }
    =
    \boldsymbol{a}.
\]
In particular,
$\boldsymbol{a}^{\top}\boldsymbol{u}
=
\sum_{i\in[n]}|u_i|^p
=
\lVert\boldsymbol{u}\rVert_p^p
=
1$,
and hence $\boldsymbol{a}\neq\boldsymbol{0}$.

We first prove \textup{(i)}.
Fix $r\in[r_0,\infty)$ and let
$\boldsymbol{x}
\in
\mathcal{B}_p\bigl(
    \boldsymbol{c}_{\boldsymbol{z}}^{\boldsymbol{u}}(r),
    r
\bigr)$.
Since $r>0$, we have
\[
    1
    \ge
    \frac{1}{r}
    \left\|
        \boldsymbol{x}
        -
        \boldsymbol{c}_{\boldsymbol{z}}^{\boldsymbol{u}}(r)
    \right\|_p
    =
    \left\|
        \boldsymbol{u}
        -
        \frac{\boldsymbol{x}-\boldsymbol{z}}{r}
    \right\|_p
    \ge
    \lVert\boldsymbol{u}\rVert_p
    -
    \frac{
        \boldsymbol{a}^{\top}
        (\boldsymbol{x}-\boldsymbol{z})
    }{r}
    =
    1
    -
    \frac{
        \boldsymbol{a}^{\top}
        (\boldsymbol{x}-\boldsymbol{z})
    }{r},
\]
where the first inequality is implied by the definition of the selected ball, the first equality is due to the definition of $\boldsymbol{c}_{\boldsymbol{z}}^{\boldsymbol{u}}(r)$, and the second inequality follows from the first-order convexity
inequality for the $\ell_p$-norm at $\boldsymbol{u}$.
Thus,
$\boldsymbol{a}^{\top}(\boldsymbol{x}-\boldsymbol{z})\ge0$,
and hence $\boldsymbol{x}\in\mathcal{H}_p$.
Since $r$ and $\boldsymbol{x}$ were arbitrary,
part~\textup{(i)} follows.
Moreover, by Lemma~\ref{lemma:radial-family-geometry},
$\boldsymbol{z}$ lies on the boundary of every member of the radial
family $\bigl<\mathcal{B}_p\bigl(
    \boldsymbol{c}_{\boldsymbol{z}}^{\boldsymbol{u}}(r),r
\bigr)\bigr>_{r=r_0}^{\infty}$.
Since $\boldsymbol{z}\in\operatorname{bd}(\mathcal{H}_p)$,
the boundary hyperplane of $\mathcal{H}_p$ supports every member of the radial family at $\boldsymbol{z}$.

To prove the necessary and sufficient conditions in part~\textup{(ii)}, suppose first that $\mu=0$.
Then
$\boldsymbol{a}^{\top}(\boldsymbol{s}-\boldsymbol{z})\le0$
for every $\boldsymbol{s}\in\operatorname{Sol}(\dd{D})$.
Since $\boldsymbol{a}\neq\boldsymbol{0}$, we have 
$\operatorname{int}(\mathcal{H}_p)
=
\{
    \boldsymbol{x}\in\mathbb{R}^n:
    \boldsymbol{a}^{\top}(\boldsymbol{x}-\boldsymbol{z})>0
\}$.
It follows that
$\operatorname{int}(\mathcal{H}_p)
\cap
\operatorname{Sol}(\dd{D})
=
\emptyset$ and hence $\mathcal{H}_p$ is DD-free.
The containment established above then implies that every member of the radial family is DD-free.
Theorem~\ref{thm:normal-ray-ball-expansion}\textup{(i)}
therefore gives $r^\ast=+\infty$.

To prove the reverse implication, we argue contrapositively.
Suppose that $\mu>0$.
Choose a maximizer
$\boldsymbol{s}^{+}\in\operatorname{Sol}(\dd{D})$
satisfying
$\boldsymbol{a}^{\top}
(\boldsymbol{s}^{+}-\boldsymbol{z})
=
\mu$.
Since
$\boldsymbol{a}^{\top}
(\boldsymbol{s}^{+}-\boldsymbol{z})
=
\mu
>
0$,
we have
$\boldsymbol{s}^{+}
\in
\operatorname{int}(\mathcal{H}_p)
\cap
\operatorname{Sol}(\dd{D})$.
Hence, $\mathcal{H}_p$ is not DD-free.
Using
$\boldsymbol{c}_{\boldsymbol{z}}^{\boldsymbol{u}}(r)
=
\boldsymbol{z}+r\boldsymbol{u}$
and
$\lVert\boldsymbol{u}\rVert_p=1$, we have
\[
\begin{aligned}
    &\lim_{r\to+\infty}
    \left(
        \left\|
            \boldsymbol{s}^{+}
            -
            \boldsymbol{c}_{\boldsymbol{z}}^{\boldsymbol{u}}(r)
        \right\|_p
        -
        r
    \right)\\
    &\quad=
    \lim_{r\to+\infty}
    r
    \left(
        \left\|
            \boldsymbol{u}
            -
            \frac{\boldsymbol{s}^{+}-\boldsymbol{z}}{r}
        \right\|_p
        -
        \lVert\boldsymbol{u}\rVert_p
    \right)\\
    &\quad=
    \lim_{t\downarrow0}
    \frac{
        \lVert
            \boldsymbol{u}
            -
            t(\boldsymbol{s}^{+}-\boldsymbol{z})
        \rVert_p
        -
        \lVert\boldsymbol{u}\rVert_p
    }{t}\\
    &\quad=
    -\boldsymbol{a}^{\top}(\boldsymbol{s}^{+}-\boldsymbol{z})\\
    &\quad=
    -\mu
    <
    0,
\end{aligned}
\]
where the second equality follows by setting $t=1/r$, and the
third follows from the differentiability of the $\ell_p$-norm at
$\boldsymbol{u}$ and its gradient $\boldsymbol{a}$.
Consequently, there exists $r_1\ge r_0$ such that
$\left\|
    \boldsymbol{s}^{+}
    -
    \boldsymbol{c}_{\boldsymbol{z}}^{\boldsymbol{u}}(r)
\right\|_p
<
r$
for every $r\ge r_1$.
Since $\boldsymbol{s}^{+}\in\operatorname{Sol}(\dd{D})$, the $\ell_p$-ball
$\mathcal{B}_p\bigl(
    \boldsymbol{c}_{\boldsymbol{z}}^{\boldsymbol{u}}(r),
    r
\bigr)$
is not DD-free for every $r\ge r_1$.
By Theorem~\ref{thm:normal-ray-ball-expansion}\textup{(i)},
it follows that $r^\ast<+\infty$.

The preceding arguments show that if $\mu=0$, then
$r^\ast=+\infty$ and $\mathcal{H}_p$ is DD-free, whereas if
$\mu>0$, then $r^\ast<+\infty$ and $\mathcal{H}_p$ is not DD-free.
Since $\mu\ge0$, the three statements in part~\textup{(ii)} are equivalent, proving~\textup{(ii)}.

Finally, we prove \textup{(iii)}.
Suppose that the equivalent statements in 
part~\textup{(ii)} hold.
Then $\mathcal{H}_p$ is a full-dimensional DD-free halfspace, and its unique facet 
$\operatorname{bd}(\mathcal{H}_p)$ contains
$\boldsymbol{z}
\in
\operatorname{Sol}(\dd{D})$ in its relative interior.
By Proposition~\ref{prop:polyhedral-maximal-characterization},
$\mathcal{H}_p$ is maximal among all DD-free sets,
proving~\textup{(iii)}.
\Halmos
\endproof

\begin{proposition}
\label{prop:radial-unbounded-l1}
    Consider the setting of
    Theorem~\ref{thm:normal-ray-ball-expansion}. 
    Suppose additionally that $p=1$.
    Let
    $\operatorname{supp}(\boldsymbol{u})
    \coloneqq
    \{i\in[n]:u_i\neq0\}$,
    and define $\Phi:\mathbb{R}^n\to\mathbb{R}$ by
    \begin{equation}
    \label{eq:radial-l1-function}
        \Phi(\boldsymbol{w})
        \coloneqq
        \sum_{i\in\operatorname{supp}(\boldsymbol{u})}
            \operatorname{sgn}(u_i)w_i
        -
        \sum_{i\in[n]\setminus\operatorname{supp}(\boldsymbol{u})}
            |w_i|.
    \end{equation}
    Define
    \begin{equation}
    \label{eq:radial-l1-limit-set}
        \mathcal{C}_{\Phi}
        \coloneqq
        \left\{
            \boldsymbol{x}\in\mathbb{R}^n:
            \Phi(\boldsymbol{x}-\boldsymbol{z})\ge0
        \right\},
    \end{equation}
    and set
    \[
        \mu
        \coloneqq
        \max_{\boldsymbol{s}\in\operatorname{Sol}(\dd{D})}
        \Phi(\boldsymbol{s}-\boldsymbol{z}).
    \]
    Then $\mu\ge0$, and the following statements hold.
    \begin{enumerate}[label=\textup{(\roman*)}, ref=\textup{(\roman*)}]
        \item The set $\mathcal{C}_{\Phi}$ is a 
        full-dimensional translated polyhedral cone.
        \item For every $r\in[r_0,\infty)$,
        $\mathcal{B}_p\bigl(
            \boldsymbol{c}_{\boldsymbol{z}}^{\boldsymbol{u}}(r),r
        \bigr)
        \subseteq
        \mathcal{C}_{\Phi}$.
        \item The following statements are equivalent:
        $r^\ast=+\infty$; 
        $\mu=0$; 
        and $\mathcal{C}_{\Phi}$ is DD-free.
        \item If the equivalent statements in part~\textup{(iii)} hold, then $\mathcal{C}_{\Phi}$ is maximal among all DD-free sets if and only if every facet $F$ of
        $\mathcal{C}_{\Phi}$ satisfies
        $\operatorname{relint}(F)
        \cap\operatorname{Sol}(\dd{D})\neq\emptyset$.
    \end{enumerate}
\end{proposition}

\proof{Proof.}

Since $\operatorname{Sol}(\dd{D})$ is nonempty and compact, the
continuous function
$\boldsymbol{s}\mapsto\Phi(\boldsymbol{s}-\boldsymbol{z})$
attains its maximum over $\operatorname{Sol}(\dd{D})$.
Moreover, since
$\boldsymbol{z}\in\operatorname{Sol}(\dd{D})$, we have
$\mu\ge\Phi(\boldsymbol{z}-\boldsymbol{z})=0$.

We first prove~\textup{(i)}.
Using $|\delta|=\max\{\delta,-\delta\}$ 
for all $\delta\in\mathbb{R}$,
\eqref{eq:radial-l1-function} and
\eqref{eq:radial-l1-limit-set} imply that
$\boldsymbol{x}\in\mathcal{C}_{\Phi}$ if and only if
\begin{equation}
\label{eq:radial-l1-signed-inequalities}
    \sum_{i\in\operatorname{supp}(\boldsymbol{u})}
        \operatorname{sgn}(u_i)(x_i-z_i)
    -
    \sum_{i\in[n]\setminus\operatorname{supp}(\boldsymbol{u})}
        \sigma_i(x_i-z_i)
    \ge0
\end{equation}
for every choice of signs
$\sigma_i\in\{-1,1\}$, for
$i\in[n]\setminus\operatorname{supp}(\boldsymbol{u})$.
These inequalities define finitely many closed halfspaces
whose boundary hyperplanes pass through $\boldsymbol{z}$.
Thus, $\mathcal{C}_{\Phi}$ is a translated polyhedral cone.
Moreover, since $\lVert\boldsymbol{u}\rVert_1=1$,
we have $\Phi(\boldsymbol{u})=\lVert\boldsymbol{u}\rVert_1=1$.
By continuity of $\Phi$,
$\boldsymbol{z}+\boldsymbol{u}
\in\operatorname{int}(\mathcal{C}_{\Phi})$,
so $\mathcal{C}_{\Phi}$ is full-dimensional.
Consequently, its interior is obtained by requiring all the defining inequalities above to hold strictly, and hence
$\operatorname{int}(\mathcal{C}_{\Phi})
=
\{\boldsymbol{x}\in\mathbb{R}^n:
\Phi(\boldsymbol{x}-\boldsymbol{z})>0\}$.
This proves~\textup{(i)}.

We next prove \textup{(ii)}.
For every $\boldsymbol{w}\in\mathbb{R}^n$,
\eqref{eq:radial-l1-function} gives
$\Phi(\boldsymbol{w})\ge-\lVert\boldsymbol{w}\rVert_1$,
because
$\operatorname{sgn}(u_i)w_i\ge-|w_i|$
for every $i\in\operatorname{supp}(\boldsymbol{u})$.
Furthermore, the definition of $\Phi$ and
$\lVert\boldsymbol{u}\rVert_1=1$ imply that
$\Phi(\boldsymbol{w}+r\boldsymbol{u})
=
\Phi(\boldsymbol{w})+r$
for every $\boldsymbol{w}\in\mathbb{R}^n$.
Fix $r\in[r_0,\infty)$ and let
$\boldsymbol{x}\in
\mathcal{B}_1\bigl(
    \boldsymbol{c}_{\boldsymbol{z}}^{\boldsymbol{u}}(r),r
\bigr)$.
Using
$\boldsymbol{c}_{\boldsymbol{z}}^{\boldsymbol{u}}(r)
=
\boldsymbol{z}+r\boldsymbol{u}$,
we obtain
\[
    \Phi(\boldsymbol{x}-\boldsymbol{z})
    =
    \Phi\left(
        \boldsymbol{x}
        -
        \boldsymbol{c}_{\boldsymbol{z}}^{\boldsymbol{u}}(r)
        +
        r\boldsymbol{u}
    \right) 
    =
    r
    +
    \Phi\left(
        \boldsymbol{x}
        -
        \boldsymbol{c}_{\boldsymbol{z}}^{\boldsymbol{u}}(r)
    \right) 
    \ge
    r
    -
    \left\|
        \boldsymbol{x}
        -
        \boldsymbol{c}_{\boldsymbol{z}}^{\boldsymbol{u}}(r)
    \right\|_1
    \ge0.
\]
Hence, $\boldsymbol{x}\in\mathcal{C}_{\Phi}$.
Since $r$ and $\boldsymbol{x}$ were arbitrary,
$\mathcal{B}_1\bigl(
    \boldsymbol{c}_{\boldsymbol{z}}^{\boldsymbol{u}}(r),r
\bigr)
\subseteq\mathcal{C}_{\Phi}$
for every $r\in[r_0,\infty)$.
This proves~\textup{(ii)}.

To prove the necessary and sufficient conditions in part~\textup{(iii)}, suppose first that $\mu=0$.
Then
$\Phi(\boldsymbol{s}-\boldsymbol{z})\le0$
for every $\boldsymbol{s}\in\operatorname{Sol}(\dd{D})$,
whereas every
$\boldsymbol{x}\in\operatorname{int}(\mathcal{C}_{\Phi})$
satisfies $\Phi(\boldsymbol{x}-\boldsymbol{z})>0$, as established above.
Hence,
$\operatorname{int}(\mathcal{C}_{\Phi})
\cap
\operatorname{Sol}(\dd{D})
=
\emptyset$.
By part~\textup{(i)}, $\mathcal{C}_{\Phi}$ is convex.
It follows that
$\mathcal{C}_{\Phi}$ is DD-free.
The containment established above then implies that every member of
the radial family $\bigl<\mathcal{B}_1\bigl(
    \boldsymbol{c}_{\boldsymbol{z}}^{\boldsymbol{u}}(r),r
\bigr)\bigr>_{r=r_0}^{\infty}$ is DD-free.
Theorem~\ref{thm:normal-ray-ball-expansion}\textup{(i)}
therefore gives $r^\ast=+\infty$.

To prove the reverse implication, we argue contrapositively.
Suppose that $\mu>0$ and choose a maximizer
$\boldsymbol{s}^{+}\in\operatorname{Sol}(\dd{D})$
satisfying
$\Phi(\boldsymbol{s}^{+}-\boldsymbol{z})=\mu$.
Since
$\Phi(\boldsymbol{s}^{+}-\boldsymbol{z})=\mu>0$,
we have $
\boldsymbol{s}^{+}
\in
\operatorname{int}(\mathcal{C}_{\Phi})
\cap
\operatorname{Sol}(\dd{D})$.
Hence, $\mathcal{C}_{\Phi}$ is not DD-free.
For any fixed $\boldsymbol{w}\in\mathbb{R}^n$ and each
$i\in\operatorname{supp}(\boldsymbol{u})$, the quantity
$ru_i-w_i$ has the same sign as $u_i$ for all sufficiently large
$r$.
Hence,
$|w_i-ru_i|
=
r|u_i|-\operatorname{sgn}(u_i)w_i$
for all sufficiently large $r$.
For each
$i\in[n]\setminus\operatorname{supp}(\boldsymbol{u})$,
we have $|w_i-ru_i|=|w_i|$ for every $r$.
Summing these coordinatewise identities and using
$\lVert\boldsymbol{u}\rVert_1=1$ yields
\begin{equation*}
\lVert\boldsymbol{w}-r\boldsymbol{u}\rVert_1-r
=
-\sum_{i\in\operatorname{supp}(\boldsymbol{u})}
            \operatorname{sgn}(u_i)w_i + \sum_{i\in[n]\setminus\operatorname{supp}(\boldsymbol{u})}
            |w_i| + r\lVert\boldsymbol{u}\rVert_1 - r
=
-\Phi(\boldsymbol{w}),
\end{equation*}
for all sufficiently large $r$.
Taking
$\boldsymbol{w}=\boldsymbol{s}^{+}-\boldsymbol{z}$
and using
$\boldsymbol{c}_{\boldsymbol{z}}^{\boldsymbol{u}}(r)
=
\boldsymbol{z}+r\boldsymbol{u}$,
we obtain, for all sufficiently large $r$,
\[
    \left\|
        \boldsymbol{s}^{+}
        -
        \boldsymbol{c}_{\boldsymbol{z}}^{\boldsymbol{u}}(r)
    \right\|_1
    -
    r
    =
    -\Phi(\boldsymbol{s}^{+}-\boldsymbol{z})
    =
    -\mu
    <
    0.
\]
Consequently, there exists $r_1\ge r_0$ such that
$\boldsymbol{s}^{+}\in
\operatorname{int}\bigl(
    \mathcal{B}_1(
        \boldsymbol{c}_{\boldsymbol{z}}^{\boldsymbol{u}}(r),r
    )
\bigr)$
for every $r\ge r_1$.
Since
$\boldsymbol{s}^{+}\in\operatorname{Sol}(\dd{D})$,
these balls are not DD-free.
Theorem~\ref{thm:normal-ray-ball-expansion}\textup{(i)}
therefore implies that $r^\ast<+\infty$.

The preceding arguments show that if $\mu=0$, then
$r^\ast=+\infty$ and $\mathcal{C}_{\Phi}$ is DD-free, whereas if
$\mu>0$, then $r^\ast<+\infty$ and $\mathcal{C}_{\Phi}$ is not
DD-free. 
Since $\mu\ge0$, the three statements in part~\textup{(iii)} are equivalent, proving~\textup{(iii)}.

Finally, we prove~\textup{(iv)}.
Suppose that the equivalent statements in 
part~\textup{(iii)} hold. 
Then $\mathcal{C}_{\Phi}$ is DD-free and, by part~\textup{(i)},
is a full-dimensional polyhedron.
\textup{(iv)} therefore follows from
Proposition~\ref{prop:polyhedral-maximal-characterization}.
\Halmos
\endproof

We will show later in Proposition~\ref{prop:normal-ray-ball-expansion-time-complexity} that condition~(ii) for verifying the DD-freeness of the halfspace in Proposition~\ref{prop:radial-unbounded-smooth} and condition~(iii) for verifying the DD-freeness of the polyhedral cone in Proposition~\ref{prop:radial-unbounded-l1} can both be checked in time linear in the size of the underlying IDD. Furthermore, the maximality condition~(iv) for the polyhedral cone in Proposition~\ref{prop:radial-unbounded-l1} can be checked using an approach analogous to that of Proposition~\ref{prop:ball-contact-subgraph} by identifying the relevant contact points through an induced subgraph of the IDD.

It follows from Theorem~\ref{thm:normal-ray-ball-expansion} that, when $r^\ast<+\infty$, the largest DD-free $\ell_p$-ball in a radial family, which could also be $\ell_p$-ball-maximal depending on $p$, can be identified by determining $r^\ast$. Our goal is therefore to design an algorithm for finding this value. To this end, it is natural to first determine a range containing $r^\ast$ by deriving lower and upper bounds that can guide the search. A simple lower bound is given by $r_0$, while Proposition~\ref{prop:radial-expansion-upper-bound} provides an explicit upper bound on $r^\ast$.

\begin{proposition}
\label{prop:radial-expansion-upper-bound}
    Consider the setting of
    Theorem~\ref{thm:normal-ray-ball-expansion}. Let
    $\mu\coloneqq
    \max_{\boldsymbol{s}\in\operatorname{Sol}(\dd{D})}
    \boldsymbol{a}^{\top}(\boldsymbol{s}-\boldsymbol{z})$
    if $p\in(1,\infty)$, and let
    $\mu\coloneqq
    \max_{\boldsymbol{s}\in\operatorname{Sol}(\dd{D})}
    \Phi(\boldsymbol{s}-\boldsymbol{z})$
    if $p=1$, where $\boldsymbol{a}$ and $\Phi$ are
    defined in
    Propositions~\ref{prop:radial-unbounded-smooth}
    and~\ref{prop:radial-unbounded-l1}, respectively.
    Suppose that $\mu>0$, and fix any
    $\boldsymbol{s}^{+}\in\operatorname{Sol}(\dd{D})$
    attaining the maximum $\mu$.
    Define
    \begin{equation}
    \label{eq:radial-expansion-upper-bound}
        r^{\mathrm{ub}}
        \coloneqq
        \begin{cases}
            \displaystyle
            \left(
                \max_{i\in\operatorname{supp}(\boldsymbol{u})}
                \frac{
                    \operatorname{sgn}(u_i)
                    (\boldsymbol{s}^{+}_i-\boldsymbol{z}_i)
                }{
                    |u_i|
                }
                -
                \frac{\mu}{2}
            \right)_+,
            & p=1,\\[4mm]
            \displaystyle
            \left(
                \frac{
                    2^{2-p}
                    \lVert
                        \boldsymbol{s}^{+}-\boldsymbol{z}
                    \rVert_p^p
                }{
                    p\mu
                }
            \right)^{1/(p-1)},
            & 1<p\le2,\\[4mm]
            \displaystyle
            \frac{
                (p-1)
                \lVert
                    \boldsymbol{s}^{+}-\boldsymbol{z}
                \rVert_p^2
            }{
                2\mu
            },
            & 2<p<\infty.
        \end{cases}
    \end{equation}
    Then
    $r^\ast\le r^{\mathrm{ub}}<+\infty$.
\end{proposition}

\proof{Proof.}

Let
$\boldsymbol{v}\coloneqq
\boldsymbol{s}^{+}-\boldsymbol{z}$.
Since
$\boldsymbol{c}_{\boldsymbol{z}}^{\boldsymbol{u}}(r)
=
\boldsymbol{z}+r\boldsymbol{u}$,
we have
$\lVert
    \boldsymbol{s}^{+}
    -
    \boldsymbol{c}_{\boldsymbol{z}}^{\boldsymbol{u}}(r)
\rVert_p
=
\lVert r\boldsymbol{u}-\boldsymbol{v}\rVert_p$.
We first show that
$\lVert r\boldsymbol{u}-\boldsymbol{v}\rVert_p<r$
for every $r>r^{\mathrm{ub}}$.
Due to their structural differences, we consider separately the cases
$p=1$, $1<p\le2$, and $p>2$.

Suppose first that $p=1$.
For each $i\in\operatorname{supp}(\boldsymbol{u})$, let
$\tau_i\coloneqq
\operatorname{sgn}(u_i)v_i/|u_i|$.
Then,
\[
\begin{aligned}
    \lVert 
        r\boldsymbol{u}-\boldsymbol{v}
    \rVert_1-r
    &=
    \sum_{i\in \operatorname{supp}(\boldsymbol{u})}
    |ru_i-v_i|
    +
    \sum_{i\in[n]\setminus \operatorname{supp}(\boldsymbol{u})}
    |v_i|
    -
    r
    \\
    &=
    \sum_{i\in \operatorname{supp}(\boldsymbol{u})}
    |u_i||r-\tau_i|
    +
    \sum_{i\in[n]\setminus \operatorname{supp}(\boldsymbol{u})}
    |v_i|
    -
    r
    \\
    &=
    \sum_{i\in \operatorname{supp}(\boldsymbol{u})}
    |u_i|\bigl(r-\tau_i+2(\tau_i-r)_+\bigr)
    +
    \sum_{i\in[n]\setminus \operatorname{supp}(\boldsymbol{u})}
    |v_i|
    -
    r
    \\
    &=
    r\left(
        \sum_{i\in \operatorname{supp}(\boldsymbol{u})}|u_i|-1
    \right)
    -
    \sum_{i\in \operatorname{supp}(\boldsymbol{u})}|u_i|\tau_i
    +
    \sum_{i\in[n]\setminus \operatorname{supp}(\boldsymbol{u})}|v_i|
    +
    2\sum_{i\in \operatorname{supp}(\boldsymbol{u})}
    |u_i|(\tau_i-r)_+
    \\
    &=
    -
    \sum_{i\in \operatorname{supp}(\boldsymbol{u})}
        \operatorname{sgn}(u_i)v_i
    +
    \sum_{i\in[n]\setminus \operatorname{supp}(\boldsymbol{u})}|v_i|
    +
    2\sum_{i\in \operatorname{supp}(\boldsymbol{u})}
    |u_i|(\tau_i-r)_+
    \\
    &=
    -\Phi(\boldsymbol{v})
    +
    2\sum_{i\in \operatorname{supp}(\boldsymbol{u})}
    |u_i|(\tau_i-r)_+
    \\
    &=
    -\mu
    +
    2\sum_{i\in \operatorname{supp}(\boldsymbol{u})}
    |u_i|(\tau_i-r)_+.
\end{aligned}
\]
In the above equations, the first equality follows from $u_i=0$ for every
$i\in[n]\setminus\operatorname{supp}(\boldsymbol{u})$.
For the second equality, since
$i\in\operatorname{supp}(\boldsymbol{u})$, the definition of
$\tau_i$ gives $v_i=u_i\tau_i$, and hence
$|ru_i-v_i|=|u_i||r-\tau_i|$.
The third equality follows from the identity
$|\delta|=\delta+2(-\delta)_+$ for $\delta\in\mathbb{R}$.
The fifth equality is due to the facts that
$\sum_{i\in\operatorname{supp}(\boldsymbol{u})}|u_i|
=\lVert\boldsymbol{u}\rVert_1=1$
and
$|u_i|\tau_i=\operatorname{sgn}(u_i)v_i$
for every
$i\in\operatorname{supp}(\boldsymbol{u})$.
Finally, the sixth equality follows from the definition of $\Phi$,
and the last equality follows from
$\Phi(\boldsymbol{v})
=\Phi(\boldsymbol{s}^{+}-\boldsymbol{z})
=\mu$.
For every $r>r^{\mathrm{ub}}$ and
$i\in\operatorname{supp}(\boldsymbol{u})$, the definition of
$r^{\mathrm{ub}}$ and the inequality $(\delta)_+\ge\delta$ give
$r>r^{\mathrm{ub}}
\ge
\max_{j\in\operatorname{supp}(\boldsymbol{u})}\tau_j-\mu/2
\ge
\tau_i-\mu/2$.
Thus, $\tau_i-r<\mu/2$.
Since $\mu>0$, it follows that
$(\tau_i-r)_+<\mu/2$
for every $i\in\operatorname{supp}(\boldsymbol{u})$.
Consequently,
\[
    \lVert r\boldsymbol{u}-\boldsymbol{v}\rVert_1-r
    =
    -\mu
    +
    2\sum_{i\in\operatorname{supp}(\boldsymbol{u})}
    |u_i|(\tau_i-r)_+
    <
    -\mu
    +
    2\sum_{i\in\operatorname{supp}(\boldsymbol{u})}
    |u_i|\frac{\mu}{2}
    =0,
\]
Hence,
$\lVert r\boldsymbol{u}-\boldsymbol{v}\rVert_1<r$
for every $r>r^{\mathrm{ub}}$, 
which proves the desired inequality for $p=1$.

Now suppose that $1<p\le2$, 
so that $p-1\in(0,1]$.
Let $\phi:\mathbb{R}\to\mathbb{R}$ be defined by
$\phi(\delta)\coloneqq
\operatorname{sgn}(\delta)|\delta|^{p-1}$.
If $\delta_1\delta_2\ge0$, assume without loss of generality
that $|\delta_1|\ge|\delta_2|$.
By the definition of $\phi$ and the subadditivity of
$t\mapsto t^{p-1}$ on $\mathbb{R}_+$,
\[
    |\phi(\delta_1)-\phi(\delta_2)|
    =
    |\delta_1|^{p-1}-|\delta_2|^{p-1}
    \le
    (|\delta_1|-|\delta_2|)^{p-1}
    =
    |\delta_1-\delta_2|^{p-1}.
\]
If $\delta_1\delta_2<0$, then, by the definition of $\phi$ and
the concavity of $t\mapsto t^{p-1}$ on $\mathbb{R}_+$,
\[
    |\phi(\delta_1)-\phi(\delta_2)|
    =
    |\delta_1|^{p-1}+|\delta_2|^{p-1}
    \le
    2\left(
        \frac{|\delta_1|+|\delta_2|}{2}
    \right)^{p-1}
    =
    2^{2-p}|\delta_1-\delta_2|^{p-1}.
\]
Since $2^{2-p}\ge1$, the two cases imply
$|\phi(\delta_1)-\phi(\delta_2)|
\le
2^{2-p}|\delta_1-\delta_2|^{p-1}$
for all $\delta_1,\delta_2\in\mathbb{R}$.
For each $i\in[n]$, the fundamental theorem of calculus gives
\[
\begin{aligned}
    \left|u_i-\frac{v_i}{r}\right|^p
    &=
    |u_i|^p
    -
    p\frac{v_i}{r}
    \int_0^1
        \phi\left(
            u_i-\kappa\frac{v_i}{r}
        \right)
    \,d\kappa
    \\
    &=
    |u_i|^p
    -
    p\phi(u_i)\frac{v_i}{r}
    -
    p\frac{v_i}{r}
    \int_0^1
    \left[
        \phi\left(
            u_i-\kappa\frac{v_i}{r}
        \right)
        -
        \phi(u_i)
    \right]
    \,d\kappa
    \\
    &\le
    |u_i|^p
    -
    p\phi(u_i)\frac{v_i}{r}
    +
    p\left|\frac{v_i}{r}\right|
    \int_0^1
    \left|
        \phi\left(
            u_i-\kappa\frac{v_i}{r}
        \right)
        -
        \phi(u_i)
    \right|
    \,d\kappa
    \\
    &\le
    |u_i|^p
    -
    p\phi(u_i)\frac{v_i}{r}
    +
    p\,2^{2-p}
    \left|\frac{v_i}{r}\right|^p
    \int_0^1
        \kappa^{p-1}
    \,d\kappa
    \\
    &=
    |u_i|^p
    -
    p\phi(u_i)\frac{v_i}{r}
    +
    2^{2-p}
    \left|\frac{v_i}{r}\right|^p.
\end{aligned}
\]
In the above equations, the first equality follows from
$\frac{d}{d\delta}|\delta|^p=p\phi(\delta)$
for $\delta\in\mathbb{R}$ and the chain rule, and the second equality is obtained by adding and subtracting $\phi(u_i)$ inside the integral.
The first inequality follows by taking absolute values and applying the triangle inequality to the integral.
For the second inequality, we apply
$|\phi(\delta_1)-\phi(\delta_2)|
\le
2^{2-p}|\delta_1-\delta_2|^{p-1}$
with
$\delta_1=u_i-\kappa v_i/r$
and
$\delta_2=u_i$
for every $\kappa\in[0,1]$.
For every $r>r^{\mathrm{ub}}$, summing over $i\in[n]$ gives
\[
\begin{aligned}
    \left\|
        \boldsymbol{u}
        -
        \frac{\boldsymbol{v}}{r}
    \right\|_p^p
    &=
    \sum_{i\in[n]}
    \left|u_i-\frac{v_i}{r}\right|^p
    \\
    &\le
    \sum_{i\in[n]}|u_i|^p
    -
    \frac{p}{r}
    \sum_{i\in[n]}\phi(u_i)v_i
    +
    \frac{2^{2-p}}{r^p}
    \sum_{i\in[n]}|v_i|^p
    \\
    &=
    1-\frac{p\mu}{r}
    +
    \frac{
        2^{2-p}\lVert\boldsymbol{v}\rVert_p^p
    }{
        r^p
    }
    \\
    &<
    1-\frac{p\mu}{r}
    +
    \frac{p\mu}{r}
    =
    1.
\end{aligned}
\]
In the above equations, the second equality uses
$\lVert\boldsymbol{u}\rVert_p=1$ and
$\sum_{i\in[n]}\phi(u_i)v_i
=
\boldsymbol{a}^{\top}\boldsymbol{v}
=
\mu$.
For the strict inequality, since $p-1>0$, the definition of
$r^{\mathrm{ub}}$ and $r>r^{\mathrm{ub}}$ imply
$r^{p-1}>
(r^{\mathrm{ub}})^{p-1}
=
2^{2-p}\lVert\boldsymbol{v}\rVert_p^p/(p\mu)$,
and hence
$2^{2-p}\lVert\boldsymbol{v}\rVert_p^p/r^p
<
p\mu/r$.
Therefore,
$\lVert r\boldsymbol{u}-\boldsymbol{v}\rVert_p<r$
for every $r>r^{\mathrm{ub}}$, 
which proves the desired inequality for $1<p\le2$.

Finally, suppose that $2<p<\infty$.
Let
$F:\mathbb{R}^n\to\mathbb{R}$ be defined by
$F(\boldsymbol{w})
\coloneqq
\frac12\lVert\boldsymbol{w}\rVert_p^2$.
The function $F$ is continuously differentiable on
$\mathbb{R}^n$ and twice continuously differentiable on
$\mathbb{R}^n\setminus\{\boldsymbol{0}\}$.
For every
$\boldsymbol{w}\in\mathbb{R}^n\setminus\{\boldsymbol{0}\}$,
direct differentiation gives
$\nabla F(\boldsymbol{w})
=
\lVert\boldsymbol{w}\rVert_p^{2-p}
\bigl(
    |w_i|^{p-2}w_i
\bigr)_{i\in[n]}$.
Moreover, for every $\boldsymbol{h}\in\mathbb{R}^n$,
\[
\begin{aligned}
    \boldsymbol{h}^{\top}
    \nabla^2F(\boldsymbol{w})
    \boldsymbol{h}
    &=
    (p-1)\lVert\boldsymbol{w}\rVert_p^{2-p}
    \sum_{i\in[n]}|w_i|^{p-2}h_i^2
    -
    (p-2)\lVert\boldsymbol{w}\rVert_p^{2-2p}
    \left(
        \sum_{i\in[n]}
        |w_i|^{p-2}w_i h_i
    \right)^2
    \\
    &\le
    (p-1)\lVert\boldsymbol{w}\rVert_p^{2-p}
    \sum_{i\in[n]}|w_i|^{p-2}h_i^2
    \\
    &\le
    (p-1)\lVert\boldsymbol{w}\rVert_p^{2-p}
    \left(
        \sum_{i\in[n]}|w_i|^p
    \right)^{\frac{p-2}{p}}
    \left(
        \sum_{i\in[n]}|h_i|^p
    \right)^{\frac{2}{p}}
    \\
    &=
    (p-1)\lVert\boldsymbol{h}\rVert_p^2.
\end{aligned}
\]
where the first inequality follows by dropping the nonpositive
second term, and the second inequality follows from H\"older's inequality with conjugate exponents $p/(p-2)$ and $p/2$.
For arbitrary
$\boldsymbol{w}_1,\boldsymbol{w}_2\in\mathbb{R}^n$,
suppose first that the line segment joining them does not contain
the origin.
A first-order Taylor expansion with integral remainder gives
\[
\begin{aligned}
    F(\boldsymbol{w}_1)
    &=
    F(\boldsymbol{w}_2)
    +
    \nabla F(\boldsymbol{w}_2)^\top
    (\boldsymbol{w}_1-\boldsymbol{w}_2)
    \\
    &\quad
    +
    \int_0^1
        (1-\kappa)
        (\boldsymbol{w}_1-\boldsymbol{w}_2)^\top
        \nabla^2F\bigl(
            \boldsymbol{w}_2
            +
            \kappa(\boldsymbol{w}_1-\boldsymbol{w}_2)
        \bigr)
        (\boldsymbol{w}_1-\boldsymbol{w}_2)
    \,d\kappa
    \\
    &\le
    F(\boldsymbol{w}_2)
    +
    \nabla F(\boldsymbol{w}_2)^\top
    (\boldsymbol{w}_1-\boldsymbol{w}_2)
    \\
    &\quad
    +
    (p-1)
    \lVert
        \boldsymbol{w}_1-\boldsymbol{w}_2
    \rVert_p^2
    \int_0^1(1-\kappa)\,d\kappa
    \\
    &=
    F(\boldsymbol{w}_2)
    +
    \nabla F(\boldsymbol{w}_2)^\top
    (\boldsymbol{w}_1-\boldsymbol{w}_2)
    +
    \frac{p-1}{2}
    \lVert
        \boldsymbol{w}_1-\boldsymbol{w}_2
    \rVert_p^2.
\end{aligned}
\]
where the inequality follows by applying
$\boldsymbol{h}^{\top}\nabla^2F(\boldsymbol{w})\boldsymbol{h}
\le
(p-1)\lVert\boldsymbol{h}\rVert_p^2$
with
$\boldsymbol{h}=\boldsymbol{w}_1-\boldsymbol{w}_2$
and
$\boldsymbol{w}
=
\boldsymbol{w}_2
+
\kappa(\boldsymbol{w}_1-\boldsymbol{w}_2)$
for every $\kappa\in[0,1]$.
Suppose next that the segment contains the origin.
Choose $\kappa_0\in[0,1]$ such that
$\boldsymbol{w}_2
+
\kappa_0(\boldsymbol{w}_1-\boldsymbol{w}_2)
=
\boldsymbol{0}$.
Then,
\[
\begin{aligned}
    &F(\boldsymbol{w}_1)
    -
    F(\boldsymbol{w}_2)
    -
    \nabla F(\boldsymbol{w}_2)^\top
    (\boldsymbol{w}_1-\boldsymbol{w}_2)
    \\
    &\quad=
    F\bigl(
        (1-\kappa_0)
        (\boldsymbol{w}_1-\boldsymbol{w}_2)
    \bigr)
    -
    F\bigl(
        -\kappa_0
        (\boldsymbol{w}_1-\boldsymbol{w}_2)
    \bigr)
    -
    \nabla F\bigl(
        -\kappa_0
        (\boldsymbol{w}_1-\boldsymbol{w}_2)
    \bigr)^\top
    (\boldsymbol{w}_1-\boldsymbol{w}_2)
    \\
    &\quad=
    (1-\kappa_0)^2
    F(\boldsymbol{w}_1-\boldsymbol{w}_2)
    -
    \kappa_0^2
    F(\boldsymbol{w}_1-\boldsymbol{w}_2)
    +
    2\kappa_0
    F(\boldsymbol{w}_1-\boldsymbol{w}_2)
    \\
    &\quad=
    F(\boldsymbol{w}_1-\boldsymbol{w}_2)
    \\
    &\quad=
    \frac12
    \lVert
        \boldsymbol{w}_1-\boldsymbol{w}_2
    \rVert_p^2
    \\
    &\quad\le
    \frac{p-1}{2}
    \lVert
        \boldsymbol{w}_1-\boldsymbol{w}_2
    \rVert_p^2.
\end{aligned}
\]
In the above equations, the first equality follows from
$\boldsymbol{w}_1
=
(1-\kappa_0)(\boldsymbol{w}_1-\boldsymbol{w}_2)$
and
$\boldsymbol{w}_2
=
-\kappa_0(\boldsymbol{w}_1-\boldsymbol{w}_2)$.
The first two terms of the second equality are due to the fact that
$F(\delta\boldsymbol{w})
=
\delta^2F(\boldsymbol{w})$
for every $\delta\in\mathbb{R}$ and $\boldsymbol{w}\in\mathbb{R}^n$.
For the third term, differentiating 
$F(\delta\boldsymbol{w})
=
\delta^2F(\boldsymbol{w})$ 
with respect to $\delta$ and applying the chain rule gives
$\nabla F(\delta\boldsymbol{w})^\top\boldsymbol{w}
=
2\delta F(\boldsymbol{w})$.
Setting
$\delta=-\kappa_0$
and
$\boldsymbol{w}=\boldsymbol{w}_1-\boldsymbol{w}_2$
therefore yields
$-\nabla F(\boldsymbol{w}_2)^\top
(\boldsymbol{w}_1-\boldsymbol{w}_2)
=
2\kappa_0
F(\boldsymbol{w}_1-\boldsymbol{w}_2)$.
The last inequality follows from $p>2$.
Thus, in either case,
\[
    F(\boldsymbol{w}_1)
    \le
    F(\boldsymbol{w}_2)
    +
    \nabla F(\boldsymbol{w}_2)^\top
    (\boldsymbol{w}_1-\boldsymbol{w}_2)
    +
    \frac{p-1}{2}
    \lVert
        \boldsymbol{w}_1-\boldsymbol{w}_2
    \rVert_p^2.
\]
For every $r>r^{\mathrm{ub}}$, we have
\[
\begin{aligned}
    \left\|
        \boldsymbol{u}
        -
        \frac{\boldsymbol{v}}{r}
    \right\|_p^2
    &=
    2F\left(
        \boldsymbol{u}
        -
        \frac{\boldsymbol{v}}{r}
    \right)
    \\
    &\le
    2F(\boldsymbol{u})
    -
    \frac{2}{r}
    \nabla F(\boldsymbol{u})^\top\boldsymbol{v}
    +
    \frac{p-1}{r^2}
    \lVert\boldsymbol{v}\rVert_p^2
    \\
    &=
    1-\frac{2\mu}{r}
    +
    \frac{
        (p-1)\lVert\boldsymbol{v}\rVert_p^2
    }{
        r^2
    }
    \\
    &<
    1-\frac{2\mu}{r}
    +
    \frac{2\mu}{r}
    =
    1.
\end{aligned}
\]
In the above equations, the first inequality follows by applying
$F(\boldsymbol{w}_1)
\le
F(\boldsymbol{w}_2)
+
\nabla F(\boldsymbol{w}_2)^\top
(\boldsymbol{w}_1-\boldsymbol{w}_2)
+
\frac{p-1}{2}
\lVert\boldsymbol{w}_1-\boldsymbol{w}_2\rVert_p^2$
with
$\boldsymbol{w}_1
=
\boldsymbol{u}-\boldsymbol{v}/r$
and
$\boldsymbol{w}_2=\boldsymbol{u}$.
The second equality uses
$F(\boldsymbol{u})=1/2$,
$\nabla F(\boldsymbol{u})=\boldsymbol{a}$,
and
$\boldsymbol{a}^{\top}\boldsymbol{v}=\mu$.
For the strict inequality, the definition of
$r^{\mathrm{ub}}$ and $r>r^{\mathrm{ub}}$ imply
$(p-1)\lVert\boldsymbol{v}\rVert_p^2/r^2<2\mu/r$.
Therefore,
$\lVert r\boldsymbol{u}-\boldsymbol{v}\rVert_p<r$
for every $r>r^{\mathrm{ub}}$, 
which proves the desired inequality for $2<p<\infty$.

Combining the three cases, 
we have
$\lVert
    \boldsymbol{s}^{+}
    -
    \boldsymbol{c}_{\boldsymbol{z}}^{\boldsymbol{u}}(r)
\rVert_p
=
\lVert r\boldsymbol{u}-\boldsymbol{v}\rVert_p
<r$
for every $r>r^{\mathrm{ub}}$.
We now prove $r^\ast\le r^{\mathrm{ub}}$.
Suppose, to the contrary, that
$r^\ast>r^{\mathrm{ub}}$.
Since $r^\ast>r^{\mathrm{ub}}$ and $r^\ast\ge r_0$, there exists a
finite $\tilde{r}\in[r_0,\infty)$ such that
$r^{\mathrm{ub}}<\tilde{r}\le r^\ast$.
By
Theorem~\ref{thm:normal-ray-ball-expansion}\textup{(i)},
$\mathcal{B}_p\bigl(
    \boldsymbol{c}_{\boldsymbol{z}}^{\boldsymbol{u}}(\tilde{r}),
    \tilde{r}
\bigr)$
is DD-free, whereas 
$\lVert
    \boldsymbol{s}^{+}
    -
    \boldsymbol{c}_{\boldsymbol{z}}^{\boldsymbol{u}}(\tilde{r})
\rVert_p
<\tilde{r}$
implies that
$\boldsymbol{s}^{+}
\in
\operatorname{int}\bigl(
    \mathcal{B}_p(
        \boldsymbol{c}_{\boldsymbol{z}}^{\boldsymbol{u}}(\tilde{r}),
        \tilde{r}
    )
\bigr)$.
This contradicts the DD-free property because
$\boldsymbol{s}^{+}\in\operatorname{Sol}(\dd{D})$.
Therefore,
$r^\ast\le r^{\mathrm{ub}}<+\infty$.
\Halmos
\endproof

Next, we present an IDD-based procedure, outlined in Algorithm~\ref{alg:normal-ray-ball-expansion}, that determines whether the radial expansion threshold is infinite and, when it is finite, computes lower and upper bounds whose gap is within a prescribed tolerance.

To specify the radial family used by the procedure, we recover a contact point
$\boldsymbol{z}
\in
\mathcal{Z}_p(\bar{\boldsymbol{x}},r_0)$
from the initial distance computation, as described in 
Proposition~\ref{prop:ball-contact-subgraph}.
The notation introduced above is then applied to the radial family 
$\bigl<\mathcal{B}_p\bigl(
    \boldsymbol{c}_{\boldsymbol{z}}^{\boldsymbol{u}}(r),r
\bigr)\bigr>_{r=r_0}^{\infty}$ induced by this choice of $\boldsymbol{z}$.

For $p\in(1,\infty)$, let $\boldsymbol{a}$ be as defined in
Proposition~\ref{prop:radial-unbounded-smooth}; for $p=1$, let
$\Phi$ be as defined in
Proposition~\ref{prop:radial-unbounded-l1}.
In either case, let $\mu$ denote the corresponding maximum value computed in those propositions.
For each $i\in[n]$ and $\dd{a}\in\dd{A}_i$, define the arc weight
\begin{equation}
\label{eq:radial-expansion-arc-weight}
    \dd{w}^{+}(\dd{a})
    \coloneqq
    \begin{cases}
        \displaystyle
        \max\left\{
            a_i\bigl(\dd{l}(\dd{a})\downarrow-z_i\bigr),
            a_i\bigl(\dd{l}(\dd{a})\uparrow-z_i\bigr)
        \right\},
        & p>1,\\[2mm]
        \displaystyle
        \max\left\{
            \operatorname{sgn}(u_i)
            \bigl(\dd{l}(\dd{a})\downarrow-z_i\bigr),
            \operatorname{sgn}(u_i)
            \bigl(\dd{l}(\dd{a})\uparrow-z_i\bigr)
        \right\},
        & p=1,\ i\in\operatorname{supp}(\boldsymbol{u}),\\[2mm]
        -\xi(\dd{a};\boldsymbol{z}),
        & p=1,\ i\notin\operatorname{supp}(\boldsymbol{u}),
    \end{cases}
\end{equation}
where $\xi(\dd{a};\boldsymbol{z})$ is computed according to \eqref{eq:arc-deviation}.
In \eqref{eq:radial-expansion-arc-weight}, each weight is the maximum, over the corresponding arc label, of the coordinate contribution to the objective defining $\mu$.
A corresponding maximizer can be chosen at an endpoint of the arc label, except when $p=1$ and
$i\notin\operatorname{supp}(\boldsymbol{u})$, in which case it is the projection of $z_i$ onto the arc label.

Using these weights,
Algorithm~\ref{alg:normal-ray-ball-expansion} computes a longest
$\dd{r}$--$\dd{t}$ path and sets $\mu$ equal to its total weight.
If $\mu=0$, the algorithm returns $\mathcal{C}_{\Phi}$ when $p=1$
and $\mathcal{H}_p$ when $p\in(1,\infty)$.
If $\mu>0$, it constructs $\boldsymbol{s}^{+}$ by selecting a
maximizing coordinate value along each arc of the longest path and
computes $r^{\mathrm{ub}}$ from
Proposition~\ref{prop:radial-expansion-upper-bound}.
The resulting upper bound initializes a bisection search for $r^\ast$.
At each trial radius $r$, the distance
$d_p\bigl(
    \boldsymbol{c}_{\boldsymbol{z}}^{\boldsymbol{u}}(r),
    \operatorname{Sol}(\dd{D})
\bigr)$
is computed by the forward pass in
Algorithm~\ref{alg:lp-ball-DD-free}.

\begin{algorithm}
\caption{Radial $\ell_p$-ball expansion}
\label{alg:normal-ray-ball-expansion}
\begin{algorithmic}[1]
    \Require an IDD
    $\dd{D}=(\dd{U},\dd{A},\dd{l}(\cdot))$ with
    $\operatorname{Sol}(\dd{D})\neq\emptyset$;
    a point
    $\bar{\boldsymbol{x}}\notin\operatorname{Sol}(\dd{D})$;
    a norm order $p\in[1,\infty)$;
    a tolerance $\varepsilon>0$
    \Ensure a DD-free set $\widehat{\mathcal{C}}$ containing
    $\mathcal{C}_p(\bar{\boldsymbol{x}})$; if the radial expansion
    threshold $r^\ast$ associated with the recovered contact point
    is finite, also an interval
    $[\underline{r},\overline{r}]$ satisfying
    $r_0\le\underline{r}\le r^\ast\le\overline{r}$ and
    $\overline{r}-\underline{r}\le\varepsilon$, with
    $\widehat{\mathcal{C}}
    =
    \mathcal{B}_p\bigl(
        \boldsymbol{c}_{\boldsymbol{z}}^{\boldsymbol{u}}
        (\underline{r}),
        \underline{r}
    \bigr)$
    \State compute
    $r_0\gets
    d_p(\bar{\boldsymbol{x}},\operatorname{Sol}(\dd{D}))$
    by the forward pass in
    Algorithm~\ref{alg:lp-ball-DD-free}, retaining predecessor
    information to recover a shortest $\dd{r}$--$\dd{t}$ path
    $\dd{q}^0$
    \State obtain $\boldsymbol{z}$ by projecting
    $\bar{\boldsymbol{x}}$ coordinatewise onto
    $\mathcal{S}_{\dd{q}^0}$, and set
    $\boldsymbol{u}\gets
    (\bar{\boldsymbol{x}}-\boldsymbol{z})/r_0$
    \State compute a longest $\dd{r}$--$\dd{t}$ path
    $\dd{q}^{+}
    =
    (\dd{a}_1^{+},\ldots,\dd{a}_n^{+})$
    with arc weights defined by
    \eqref{eq:radial-expansion-arc-weight}
    \State
    $\mu\gets
    \sum_{i\in[n]}
    \dd{w}^{+}(\dd{a}_i^{+})$
    \If{$\mu=0$}
        \If{$p=1$}
            \State\Return $\mathcal{C}_{\Phi}$ as defined in
            \eqref{eq:radial-l1-limit-set}
        \Else
            \State\Return $\mathcal{H}_p$ as defined in
            \eqref{eq:normal-ray-limit-halfspace}
        \EndIf
    \EndIf
    \For{$i\in[n]$}
        \State choose
        $s_i^{+}\in\dd{l}(\dd{a}_i^{+})$
        whose coordinate contribution to the objective defining
        $\mu$ equals
        $\dd{w}^{+}(\dd{a}_i^{+})$
    \EndFor
    \State
    $\boldsymbol{s}^{+}\gets
    (s_1^{+},\ldots,s_n^{+})$
    \State compute $r^{\mathrm{ub}}$ from
    \eqref{eq:radial-expansion-upper-bound}
    \State
    $\underline{r}\gets r_0$;
    $\overline{r}\gets r^{\mathrm{ub}}$
    \While{$\overline{r}-\underline{r}>\varepsilon$}
        \State
        $r\gets(\underline{r}+\overline{r})/2$
        \State compute
        $d_r
        \gets
        d_p\bigl(
            \boldsymbol{c}_{\boldsymbol{z}}^{\boldsymbol{u}}(r),
            \operatorname{Sol}(\dd{D})
        \bigr)$
        by the forward pass in
        Algorithm~\ref{alg:lp-ball-DD-free}
        \If{$d_r=r$}
            \State
            $\underline{r}\gets r$
        \Else
            \State
            $\overline{r}\gets r$
        \EndIf
    \EndWhile
    \State\Return
    $\left(
        \mathcal{B}_p\bigl(
            \boldsymbol{c}_{\boldsymbol{z}}^{\boldsymbol{u}}
            (\underline{r}),
            \underline{r}
        \bigr),
        [\underline{r},\overline{r}]
    \right)$
\end{algorithmic}
\end{algorithm}

We next establish the correctness and running-time complexity of
Algorithm~\ref{alg:normal-ray-ball-expansion}.
Theorem~\ref{thm:normal-ray-ball-expansion-valid} shows that the algorithm returns a DD-free set containing
$\mathcal{C}_p(\bar{\boldsymbol{x}})$.
If $r^\ast<+\infty$, the algorithm returns the DD-free $\ell_p$-ball
$\mathcal{B}_p\bigl(
    \boldsymbol{c}_{\boldsymbol{z}}^{\boldsymbol{u}}(\underline{r}),
    \underline{r}
\bigr)$
satisfying
$0\le r^\ast-\underline{r}\le\varepsilon$, so that its radius is within $\varepsilon$ of that of the largest DD-free member of the radial family
$\bigl<\mathcal{B}_p\bigl(
    \boldsymbol{c}_{\boldsymbol{z}}^{\boldsymbol{u}}(r),r
\bigr)\bigr>_{r=r_0}^{\infty}$, characterized in
Theorem~\ref{thm:normal-ray-ball-expansion}.
Proposition~\ref{prop:normal-ray-ball-expansion-time-complexity}
gives the corresponding running-time bound.

\begin{theorem}
\label{thm:normal-ray-ball-expansion-valid}
    Consider the setting of
    Theorem~\ref{thm:normal-ray-ball-expansion}.
    Let $\boldsymbol{z}$ be the contact point computed in line 2 of
    Algorithm~\ref{alg:normal-ray-ball-expansion}, and let $\mu$ be the value calculated in line 4. 
    Fix $\varepsilon>0$.
    Then Algorithm~\ref{alg:normal-ray-ball-expansion} terminates and
    returns a DD-free set $\widehat{\mathcal{C}}$ containing
    $\mathcal{C}_p(\bar{\boldsymbol{x}})$.
    Moreover, the following statements hold.
    \begin{enumerate}[label=\textup{(\roman*)}, ref=\textup{(\roman*)}]
        \item
        If $\mu = 0$, then $r^\ast=+\infty$, and the algorithm returns
        $\widehat{\mathcal{C}}=\mathcal{C}_{\Phi}$ when $p=1$ and
        $\widehat{\mathcal{C}}=\mathcal{H}_p$ when
        $p\in(1,\infty)$.
        \item
        If $\mu > 0$, then $r^\ast<+\infty$, and the algorithm returns
        $\widehat{\mathcal{C}}
        =
        \mathcal{B}_p\bigl(
            \boldsymbol{c}_{\boldsymbol{z}}^{\boldsymbol{u}}
            (\underline{r}),
            \underline{r}
        \bigr)$
        together with an interval
        $[\underline{r},\overline{r}]$ satisfying
        $r_0\le\underline{r}\le r^\ast\le\overline{r}$ and
        $\overline{r}-\underline{r}\le\varepsilon$.
        In particular,
        $0\le r^\ast-\underline{r}\le\varepsilon$.
    \end{enumerate}
\end{theorem}

\proof{Proof.}

By Theorem~\ref{thm:lp-ball-DD-free-valid}, the initial forward pass computes
$r_0=d_p(\bar{\boldsymbol{x}},\operatorname{Sol}(\dd{D}))$.
Since
$\bar{\boldsymbol{x}}\notin\operatorname{Sol}(\dd{D})$,
we have $r_0>0$.
By Proposition~\ref{prop:ball-contact-subgraph},
the recovery of a shortest path followed by coordinatewise projection in lines 1--2 yields
$\boldsymbol{z}
\in
\mathcal{Z}_p(\bar{\boldsymbol{x}},r_0)$.
Thus,
$\lVert
    \bar{\boldsymbol{x}}-\boldsymbol{z}
\rVert_p=r_0$
and
$\lVert\boldsymbol{u}\rVert_p=1$.
Then, \eqref{eq:radial-center} implies that
$\boldsymbol{c}_{\boldsymbol{z}}^{\boldsymbol{u}}(r_0)
=\bar{\boldsymbol{x}}$, and hence
$\mathcal{B}_p\bigl(
    \boldsymbol{c}_{\boldsymbol{z}}^{\boldsymbol{u}}(r_0),r_0
\bigr)
=
\mathcal{C}_p(\bar{\boldsymbol{x}})$.

We next verify the computation of $\mu$.
Consider an $\dd{r}$--$\dd{t}$ path
$\dd{q}=(\dd{a}_1,\ldots,\dd{a}_n)$.
By \eqref{eq:radial-expansion-arc-weight}, each
$\dd{w}^{+}(\dd{a}_i)$ is the maximum coordinate contribution to
the objective defining $\mu$ over
$\dd{l}(\dd{a}_i)$.
Since
$\mathcal{S}_{\dd{q}}
=\prod_{i\in[n]}\dd{l}(\dd{a}_i)$,
the sum
$\sum_{i\in[n]}\dd{w}^{+}(\dd{a}_i)$
is the maximum of
$\boldsymbol{a}^{\top}(\boldsymbol{s}-\boldsymbol{z})$
over $\mathcal{S}_{\dd{q}}$ when $p>1$, and of
$\Phi(\boldsymbol{s}-\boldsymbol{z})$
over $\mathcal{S}_{\dd{q}}$ when $p=1$.
Taking the maximum over all $\dd{r}$--$\dd{t}$ paths and using
$\operatorname{Sol}(\dd{D})
=
\bigcup_{\dd{q}\in\dd{P}(\dd{t})}
\mathcal{S}_{\dd{q}}$
shows that the maximum total path weight is exactly $\mu$.
In particular, $\mu\ge0$, because
$\boldsymbol{z}\in\operatorname{Sol}(\dd{D})$ gives objective
value zero.
Propositions~\ref{prop:radial-unbounded-smooth}
and~\ref{prop:radial-unbounded-l1} imply that
$\mu=0$ if and only if $r^\ast=+\infty$.

Suppose first that $\mu=0$.
Algorithm~\ref{alg:normal-ray-ball-expansion} then returns
$\mathcal{C}_{\Phi}$ when $p=1$ and
$\mathcal{H}_p$ when $p\in(1,\infty)$, which, by
Propositions~\ref{prop:radial-unbounded-l1}
and~\ref{prop:radial-unbounded-smooth}, respectively, are DD-free and contain every member of the radial family
$\bigl<\mathcal{B}_p\bigl(
    \boldsymbol{c}_{\boldsymbol{z}}^{\boldsymbol{u}}(r),r
\bigr)\bigr>_{r=r_0}^{\infty}$.
It therefore contains
$\mathcal{B}_p\bigl(
    \boldsymbol{c}_{\boldsymbol{z}}^{\boldsymbol{u}}(r_0),r_0
\bigr)
=
\mathcal{C}_p(\bar{\boldsymbol{x}})$.
The algorithm terminates before entering the bisection loop, which proves~\textup{(i)}.

Now suppose that $\mu>0$.
The coordinatewise maximizers selected along the longest path
$\dd{q}^{+}$ form a point
$\boldsymbol{s}^{+}\in
\mathcal{S}_{\dd{q}^{+}}
\subseteq\operatorname{Sol}(\dd{D})$
that attains $\mu$.
It follows from Theorem~\ref{thm:normal-ray-ball-expansion} that
$r_0\le r^\ast$, while
Proposition~\ref{prop:radial-expansion-upper-bound} gives
$r^\ast\le r^{\mathrm{ub}}<+\infty$.
Thus, at initialization,
$r_0\le\underline{r}\le r^\ast\le\overline{r}$.
We next show that this inequality chain is preserved by each bisection update in lines 22 and 24.

Let $r=(\underline{r}+\overline{r})/2$ be the trial radius.
By Theorem~\ref{thm:lp-ball-DD-free-valid}, the forward pass
computes
$d_r=d_p\bigl(
    \boldsymbol{c}_{\boldsymbol{z}}^{\boldsymbol{u}}(r),
    \operatorname{Sol}(\dd{D})
\bigr)$
exactly.
Since
$\boldsymbol{z}\in\operatorname{Sol}(\dd{D})$ and
$\lVert
    \boldsymbol{c}_{\boldsymbol{z}}^{\boldsymbol{u}}(r)
    -\boldsymbol{z}
\rVert_p=r$,
we have $d_r\le r$.
By Definition~\ref{def:DD-free-sets}, the $\ell_p$-ball
$\mathcal{B}_p\bigl(
    \boldsymbol{c}_{\boldsymbol{z}}^{\boldsymbol{u}}(r),r
\bigr)$
is DD-free if and only if $d_r\ge r$, and therefore if and only if
$d_r=r$.
Theorem~\ref{thm:normal-ray-ball-expansion}\textup{(i)} then implies that this
equality holds if and only if $r\le r^\ast$.
Consequently, replacing $\underline{r}$ with $r$ when $d_r=r$
preserves the inequality chain.
Otherwise, $d_r<r$ and $r>r^\ast$, so replacing
$\overline{r}$ with $r$ also preserves it.

Each iteration halves the interval length.
Since the initial length
$r^{\mathrm{ub}}-r_0$ is finite and $\varepsilon>0$, the loop terminates with
$\overline{r}-\underline{r}\le\varepsilon$.
The maintained inequality chain then gives
$0\le r^\ast-\underline{r}
\le\overline{r}-\underline{r}
\le\varepsilon$.
Theorem~\ref{thm:normal-ray-ball-expansion}\textup{(i)} shows that
the returned $\ell_p$-ball is DD-free, and
Lemma~\ref{lemma:radial-family-geometry} gives
$\mathcal{C}_p(\bar{\boldsymbol{x}})
=
\mathcal{B}_p\bigl(
    \boldsymbol{c}_{\boldsymbol{z}}^{\boldsymbol{u}}(r_0),r_0
\bigr)
\subseteq
\mathcal{B}_p\bigl(
    \boldsymbol{c}_{\boldsymbol{z}}^{\boldsymbol{u}}
    (\underline{r}),
    \underline{r}
\bigr)$.
This proves~\textup{(ii)} and completes the proof.
\Halmos
\endproof

\begin{proposition}
\label{prop:normal-ray-ball-expansion-time-complexity}
    Consider the setting of
    Theorem~\ref{thm:normal-ray-ball-expansion-valid}.
    If $\mu = 0$, then
    Algorithm~\ref{alg:normal-ray-ball-expansion}
    can be implemented in $\mathcal{O}(|\dd{A}|)$ time.
    If $\mu > 0$, let
    $N_{\varepsilon}\coloneqq
    \left\lceil
        \log_2\bigl(
            \max\{1,(r^{\mathrm{ub}}-r_0)/\varepsilon\}
        \bigr)
    \right\rceil$.
    Then the algorithm performs exactly $N_{\varepsilon}$
    bisection iterations and can be implemented in
    $\mathcal{O}\bigl(
        |\dd{A}|(1+N_{\varepsilon})
    \bigr)$ time.
\end{proposition}

\proof{Proof.}

It follows from Proposition~\ref{prop:lp-ball-DD-free-time-complexity} that the initial distance computation in line 1 can be performed in
$\mathcal{O}(|\dd{A}|)$ time.
The projection step in line 2 can be done in $\mathcal{O}(n)$.
In line 3, the arc weights in \eqref{eq:radial-expansion-arc-weight} can be computed in $\mathcal{O}(|\dd{A}|)$ time.
Because the IDD is layered and acyclic, a longest
$\dd{r}$--$\dd{t}$ path can be computed in
$\mathcal{O}(|\dd{U}|+|\dd{A}|)$ time.
Both path computations can retain predecessor information without changing these bounds, and all associated vector operations require $\mathcal{O}(n)$ time.
By \eqref{eq:DD-size-bounds},
$|\dd{U}|\le 2|\dd{A}|$ and $n\le|\dd{A}|$.
Thus, all work through the computation of $\mu$ in line 4 requires
$\mathcal{O}(|\dd{A}|)$ time.
By Propositions~\ref{prop:radial-unbounded-smooth}
and~\ref{prop:radial-unbounded-l1},
$\mu=0$ if and only if $r^\ast=+\infty$.

Suppose first that $\mu=0$, equivalently,
$r^\ast=+\infty$.
Then Algorithm~\ref{alg:normal-ray-ball-expansion}
terminates before entering the bisection loop, and the total running 
time is $\mathcal{O}(|\dd{A}|)$.

Now suppose that $\mu>0$, equivalently,
$r^\ast<+\infty$.
Constructing $\boldsymbol{s}^{+}$, computing $r^{\mathrm{ub}}$, and initializing the bisection interval in lines 12--17 require
$\mathcal{O}(n)$ time.
The initial interval length is 
$r^{\mathrm{ub}}-r_0$, and after $k$ bisection iterations it is 
$(r^{\mathrm{ub}}-r_0)/2^k$.
Consequently, the loop in lines 18--26 performs exactly $N_\varepsilon$ iterations, including no iteration when 
$r^{\mathrm{ub}}-r_0\le\varepsilon$.
Each bisection iteration requires a single distance computation in line 20, which takes $\mathcal{O}(|\dd{A}|)$ time by
Proposition~\ref{prop:lp-ball-DD-free-time-complexity}.
The remaining operations require $\mathcal{O}(n)$ time, so
$n\le|\dd{A}|$ implies that each iteration can be executed in
$\mathcal{O}(|\dd{A}|)$ time.
Therefore, the total running time is
$\mathcal{O}\bigl(|\dd{A}|(1+N_\varepsilon)\bigr)$.
\Halmos
\endproof

When $\mu>0$, the maximality conditions in
Theorems~\ref{thm:ball-maximal-characterization}
and~\ref{thm:l1-ball-maximal-characterization}
provide a sufficient criterion for early termination of the bisection search in
Algorithm~\ref{alg:normal-ray-ball-expansion}.
The following corollary establishes that the corresponding maximality test can be performed directly on the IDD 
in $\mathcal{O}(|\dd{A}|)$ time.

\begin{corollary}
\label{cor:radial-ball-maximality-test}
    Consider the setting of
    Proposition~\ref{prop:normal-ray-ball-expansion-time-complexity}
    with $\mu>0$.
    For any $\underline{r}\in[r_0,r^\ast]$,
    whether $\mathcal{B}_p(
        \boldsymbol{c}_{\boldsymbol{z}}^{\boldsymbol{u}}(\underline{r}),
        \underline{r}
    )$ is $\ell_p$-ball-maximal can be determined 
    in $\mathcal{O}(|\dd{A}|)$ time.
    If this ball is $\ell_p$-ball-maximal, 
    then $\underline{r}=r^\ast$.
\end{corollary}

\proof{Proof.}

Since $\underline{r}\le r^\ast$,
Theorem~\ref{thm:normal-ray-ball-expansion}\textup{(i)}
implies that
$\mathcal{B}_p(
    \boldsymbol{c}_{\boldsymbol{z}}^{\boldsymbol{u}}(\underline{r}),
    \underline{r}
)$ is DD-free.
By Lemma~\ref{lemma:radial-family-geometry},
$\boldsymbol{z}$ is a contact point of this ball, and hence
$d_p(
    \boldsymbol{c}_{\boldsymbol{z}}^{\boldsymbol{u}}(\underline{r}),
    \operatorname{Sol}(\dd{D})
)=\underline{r}$.
By Proposition~\ref{prop:ball-contact-subgraph},
we can construct a subgraph $\dd{D}^{\partial}$ consisting of all arcs lying on $\dd{r}$--$\dd{t}$ paths
whose induced hyperrectangles contain a contact point of this ball.

For $p\in(1,\infty)$, each retained path determines a unique contact point through the coordinatewise projection of 
$\boldsymbol{c}_{\boldsymbol{z}}^{\boldsymbol{u}}(\underline{r})$
onto its path-induced hyperrectangle, and every contact point arises in this way.
Since $\boldsymbol{z}$ is a contact point and each shortest-path hyperrectangle is at distance
$\underline{r}$ from 
$\boldsymbol{c}_{\boldsymbol{z}}^{\boldsymbol{u}}(\underline{r})$, 
this projection differs from $\boldsymbol{z}$ if and only if the hyperrectangle does not contain $\boldsymbol{z}$.
Thus, by Theorem~\ref{thm:ball-maximal-characterization},
the ball is $\ell_p$-ball-maximal if and only if
there exist $i\in[n]$ and a retained arc
$\dd{a}\in\dd{A}_i$ such that
$z_i\notin\dd{l}(\dd{a})$.

For $p=1$, the retained paths are precisely the active paths in $\dd{P}^*$ in
Theorem~\ref{thm:l1-ball-maximal-characterization}.
Hence, the coordinate bounds in
\eqref{eq:l1-ball-maximal-coordinate-bounds}
can be evaluated by taking the minimum lower endpoint and maximum upper endpoint of the retained arc labels in each layer.
Checking the corresponding strict inequalities for the coordinates of 
$\boldsymbol{c}_{\boldsymbol{z}}^{\boldsymbol{u}}(\underline{r})$ therefore determines whether the ball is $\ell_1$-ball-maximal.

By Proposition~\ref{prop:ball-contact-subgraph},
the subgraph $\dd{D}^{\partial}$ can be constructed in $\mathcal{O}(|\dd{A}|)$ time.
In either case, the above checks require
$\mathcal{O}(|\dd{A}|+n)$ time,
because each retained arc and each coordinate is processed at most a constant number of times.
Since $n\le|\dd{A}|$ by \eqref{eq:DD-size-bounds},
the maximality check requires $\mathcal{O}(|\dd{A}|)$ time.

Finally, 
suppose that the ball is $\ell_p$-ball-maximal.
If $\underline{r}<r^\ast$, then
Lemma~\ref{lemma:radial-family-geometry}
and Theorem~\ref{thm:normal-ray-ball-expansion}\textup{(i)}
imply that a strictly larger DD-free ball exists
in the same radial family, contradicting
$\ell_p$-ball-maximality.
Since $\underline{r}\le r^\ast$ by assumption,
we conclude that $\underline{r}=r^\ast$.
\Halmos
\endproof

While Corollary~\ref{cor:radial-ball-maximality-test} 
provides a sufficient condition for a lower bound $\underline{r}$ to equal $r^\ast$, thereby enabling early termination, one can also use the upper bound $\overline{r}$ to obtain an alternative early-termination criterion.
The idea is to compute a nearest point in $\operatorname{Sol}(\dd{D})$ to $\boldsymbol{c}_{\boldsymbol{z}}^{\boldsymbol{u}}(\overline{r})$ and then determine the radius at which the corresponding center is equidistant from this point and the initial contact point $\boldsymbol{z}$. The following proposition shows that the radius obtained from this equidistance condition yields a tighter upper bound on $r^\ast$ and coincides with $r^\ast$ if and only if the selected nearest point is also a contact point of the radial ball at radius $r^\ast$.

\begin{proposition}
\label{prop:radial-contact-point-update}
    Consider the setting of
    Theorem~\ref{thm:normal-ray-ball-expansion}
    with $r^\ast<+\infty$.
    Fix a finite upper bound $\overline{r}>r^\ast$,
    and let $\overline{\boldsymbol{z}}$ be any nearest point in 
    $\operatorname{Sol}(\dd{D})$ to
    $\boldsymbol{c}_{\boldsymbol{z}}^{\boldsymbol{u}}(\overline{r})$
    under the $\ell_p$-norm.
    Then, the equation
    \begin{equation}
    \label{eq:radial-contact-point-update}
        \left\|
            \boldsymbol{c}_{\boldsymbol{z}}^{\boldsymbol{u}}(r)
            -
            \overline{\boldsymbol{z}}
        \right\|_p
        =
        r
    \end{equation}
    has a unique solution
    $\widehat{r}\in[r^\ast,\overline{r})$.
    Moreover, $\widehat{r}=r^\ast$ if and only if
    $\overline{\boldsymbol{z}}
    \in
    \mathcal{Z}_p(
        \boldsymbol{c}_{\boldsymbol{z}}^{\boldsymbol{u}}(r^\ast),
        r^\ast
    )$.
    In particular, when $p=2$,
    \begin{equation}
    \label{eq:radial-contact-point-update-l2}
        \widehat{r}
        =
        \frac{
            \left\|
                \overline{\boldsymbol{z}}-\boldsymbol{z}
            \right\|_2^2
        }{
            2\boldsymbol{u}^{\top}
            (\overline{\boldsymbol{z}}-\boldsymbol{z})
        }.
    \end{equation}
\end{proposition}

\proof{Proof.}

For $r\in[r_0,\infty)$, define
$F(r)\coloneqq
\lVert
    \boldsymbol{c}_{\boldsymbol{z}}^{\boldsymbol{u}}(r)
    -
    \overline{\boldsymbol{z}}
\rVert_p-r$.
Since $\boldsymbol{c}_{\boldsymbol{z}}^{\boldsymbol{u}}(r)
=\boldsymbol{z}+r\boldsymbol{u}$,
the function $F$ is continuous and convex.
Moreover, for any $r_0\le r\le r'$, 
\[
    F(r')-F(r)
    \le
    \left\|
        \boldsymbol{c}_{\boldsymbol{z}}^{\boldsymbol{u}}(r')
        -
        \boldsymbol{c}_{\boldsymbol{z}}^{\boldsymbol{u}}(r)
    \right\|_p
    -(r'-r)
    =
    (r'-r)(\lVert\boldsymbol{u}\rVert_p-1)
    =0,
\]
where the inequality follows from the triangle inequality, and the last equality holds because $\lVert\boldsymbol{u}\rVert_p=1$.
Thus, $F$ is nonincreasing.

It follows from Theorem~\ref{thm:normal-ray-ball-expansion}\textup{(i)} that
$\mathcal{B}_p\bigl(
    \boldsymbol{c}_{\boldsymbol{z}}^{\boldsymbol{u}}(r^\ast),
    r^\ast
\bigr)$ is DD-free, while the member of the radial family
with radius $\overline{r}>r^\ast$ is not.
Since $\overline{\boldsymbol{z}}\in\operatorname{Sol}(\dd{D})$,
the DD-free property yields
$\lVert
    \boldsymbol{c}_{\boldsymbol{z}}^{\boldsymbol{u}}(r^\ast)
    -
    \overline{\boldsymbol{z}}
\rVert_p\ge r^\ast$,
and hence $F(r^\ast)\ge0$.
On the other hand, because $\overline{\boldsymbol{z}}$
is a nearest point in $\operatorname{Sol}(\dd{D})$ to
$\boldsymbol{c}_{\boldsymbol{z}}^{\boldsymbol{u}}(\overline{r})$,
we have
$F(\overline{r})
=
d_p\bigl(
    \boldsymbol{c}_{\boldsymbol{z}}^{\boldsymbol{u}}(\overline{r}),
    \operatorname{Sol}(\dd{D})
\bigr)-\overline{r}<0$.
By continuity, $F$ has a zero
$\widehat{r}\in[r^\ast,\overline{r})$.
To prove uniqueness, suppose that
$F(\delta_1)=F(\delta_2)=0$ for some $r_0\le \delta_1<\delta_2$.
Because $F$ is nonincreasing and $F(\overline{r})<0$,
we must have $\delta_2<\overline{r}$.
Since $\delta_1<\delta_2<\overline{r}$, 
we can write
$\delta_2=(1-\lambda)\delta_1+\lambda\overline{r}$,
where $\lambda\coloneqq
(\delta_2-\delta_1)/(\overline{r}-\delta_1)\in(0,1)$.
By convexity,
$F(\delta_2)
\le
(1-\lambda)F(\delta_1)+\lambda F(\overline{r})
=\lambda F(\overline{r})<0$,
contradicting $F(\delta_2)=0$.
Thus, \eqref{eq:radial-contact-point-update}
has the unique solution $\widehat{r}$.
It follows from the definition of the contact set that
$\overline{\boldsymbol{z}}\in
\mathcal{Z}_p(
    \boldsymbol{c}_{\boldsymbol{z}}^{\boldsymbol{u}}(r^\ast),
    r^\ast
)$ if and only if $F(r^\ast)=0$.
Since $\widehat{r}$ is the unique zero of $F$,
this is equivalent to $\widehat{r}=r^\ast$.

Finally, suppose that $p=2$.
Since $F(\overline{r})<0$, we have
$\overline{\boldsymbol{z}}\neq\boldsymbol{z}$.
Squaring \eqref{eq:radial-contact-point-update}
at $r=\widehat{r}$ and using
$\lVert\boldsymbol{u}\rVert_2=1$ gives
$2\widehat{r}\boldsymbol{u}^{\top}
(\overline{\boldsymbol{z}}-\boldsymbol{z})
=
\lVert\overline{\boldsymbol{z}}-\boldsymbol{z}\rVert_2^2
>0$.
Since $\widehat{r}\ge r_0>0$,
we obtain
$2\boldsymbol{u}^{\top}
(\overline{\boldsymbol{z}}-\boldsymbol{z})>0$.
Rearranging yields
\eqref{eq:radial-contact-point-update-l2}.
\Halmos
\endproof

To apply Proposition~\ref{prop:radial-contact-point-update},
consider a current upper bound $\overline{r}\ge r^\ast$.
The distance
$d_p\bigl(
    \boldsymbol{c}_{\boldsymbol{z}}^{\boldsymbol{u}}(\overline{r}),
    \operatorname{Sol}(\dd{D})
\bigr)$
can be computed through
Algorithm~\ref{alg:lp-ball-DD-free}.
If this distance equals $\overline{r}$,
then $\overline{r}=r^\ast$.
Otherwise, 
by Proposition~\ref{prop:distance-to-path},
a nearest point $\overline{\boldsymbol{z}}$ can be recovered by projecting the center coordinatewise onto the hyperrectangle induced by a shortest 
$\dd{r}$--$\dd{t}$ path of the weighted $\dd{D}$.
For $p\neq2$,
\eqref{eq:radial-contact-point-update}
can be solved numerically using a univariate
root-finding method.
After solving \eqref{eq:radial-contact-point-update}
for $\widehat{r}$,
another distance computation at
$\boldsymbol{c}_{\boldsymbol{z}}^{\boldsymbol{u}}(\widehat{r})$
determines whether $\widehat{r}=r^\ast$.
If the distance equals $\widehat{r}$,
then the threshold is attained;
otherwise, $\overline{r}$ can be updated to
$\widehat{r}$.
Thus, this approach provides an alternative to the
bisection search and can, in certain cases,
recover $r^\ast$ exactly.

\begin{figure}[tb]
    \centering
    \includegraphics[width=\linewidth]{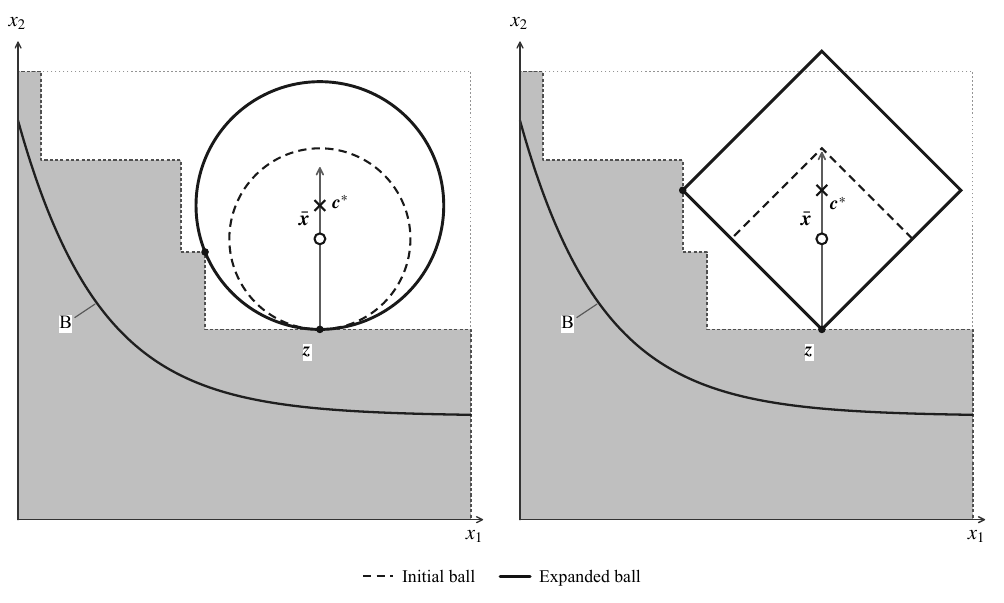}
    \caption{
        DD-free $\ell_p$-balls obtained by radial expansion,
        illustrated at $r=r^\ast$ for $p=2$ (left) and
        $p=1$ (right).
        The solid curve labeled $A$ depicts the nonlinear
        boundary of $\mathcal{G}^A$, with the feasible region
        lying on the side containing the origin.
        The shaded region represents
        $\operatorname{Sol}(\dd{D}^A)$.
        In each panel, the dashed boundary encloses the
        initial ball $\mathcal{C}_p(\bar{\boldsymbol{x}})$,
        whereas the thick solid boundary encloses the expanded
        ball $\mathcal{B}_p(\boldsymbol{c}^{\ast},r^\ast)$,
        where $\boldsymbol{c}^{\ast}\coloneqq
        \boldsymbol{c}_{\boldsymbol{z}}^{\boldsymbol{u}}(r^\ast)$.
        The gray arrow indicates the direction of center
        movement, the cross marks $\boldsymbol{c}^{\ast}$,
        and the two filled dots mark contact points with
        $\operatorname{Sol}(\dd{D}^A)$.
    }
    \label{fig:radial-ball-expansion}
\end{figure}

To illustrate the radial expansion procedure, we return to the example in
Section~\ref{subsec:illustrative-example},
focusing on constraint $A$ with the same IDD relaxation
$\dd{D}^A$ and point $\bar{\boldsymbol{x}}$.
Figure~\ref{fig:radial-ball-expansion} shows the
radial expansion of the initial DD-free $\ell_p$-ball
$\mathcal{C}_p(\bar{\boldsymbol{x}})$
considered in Algorithm~\ref{alg:normal-ray-ball-expansion}.
In each panel, the expanded ball is the largest DD-free member of the corresponding radial family, 
attained at $r=r^\ast$.

\subsubsection{Intersection Cuts from Expanded Balls}
\label{subsubsec:ball-IC}
We now return to the corner relaxation setting of
Section~\ref{subsec:IC-derivation}.
Let 
$\bar{\boldsymbol{x}}
\in
\mathcal{P}
\setminus
\operatorname{Sol}(\dd{D})$ 
be a basic feasible solution of $\mathcal{P}$,
and fix an associated corner relaxation 
$\mathcal{P}'
=\bar{\boldsymbol{x}}
+\operatorname{cone}\{
    \boldsymbol{v}^j:j\in R
\}$.
Recall that $R$ indexes the extreme rays of $\mathcal{P}'$, and every
$\boldsymbol{x}\in\mathcal{P}'$ admits a unique representation
$\boldsymbol{x}
=\bar{\boldsymbol{x}}
+\sum_{j\in R}
\lambda_j(\boldsymbol{x})\boldsymbol{v}^j$,
where $\lambda_j(\boldsymbol{x})\ge0$ for all $j\in R$.

Fix $p\in[1,\infty)$ and $\varepsilon>0$, and let
$\widehat{\mathcal{C}}$ denote the set returned by
Algorithm~\ref{alg:normal-ray-ball-expansion} at
$\bar{\boldsymbol{x}}$ with norm order $p$ and tolerance
$\varepsilon$.
Write $r_0\coloneqq
d_p(\bar{\boldsymbol{x}},\operatorname{Sol}(\dd{D}))>0$, and let
$\boldsymbol{z}$ and $\boldsymbol{u}$ be the contact point and
radial direction computed by the algorithm, so that
$\bar{\boldsymbol{x}}=\boldsymbol{z}+r_0\boldsymbol{u}$ and
$\lVert\boldsymbol{u}\rVert_p=1$.

Suppose first that $r^\ast<+\infty$.
In this case, the algorithm returns
$\widehat{\mathcal{C}}
=\mathcal{B}_p\bigl(
    \boldsymbol{c}_{\boldsymbol{z}}^{\boldsymbol{u}}
    (\underline{r}),\underline{r}
\bigr)$,
where $\underline{r}\in[r_0,r^\ast]$ is the radius lower bound
returned by Algorithm~\ref{alg:normal-ray-ball-expansion}.
The following proposition characterizes the gauge of
$\widehat{\mathcal{C}}-\bar{\boldsymbol{x}}$ in this case and gives an explicit expression when $p=2$.

\begin{proposition}
\label{prop:expanded-ball-gauge}
    Consider the setting of
    Theorem~\ref{thm:normal-ray-ball-expansion-valid},
    and suppose that $r^\ast<+\infty$.
    For every $\boldsymbol{w}\in\mathbb{R}^n$, we have
    $\gamma_{\widehat{\mathcal{C}}-\bar{\boldsymbol{x}}}
    (\boldsymbol{w})=\tau_p(\boldsymbol{w})$, where
    $\tau_p(\boldsymbol{w})$ is the unique nonnegative scalar satisfying
    \begin{equation}
    \label{eq:ball-root-equation}
        \left\|
            \boldsymbol{w}
            -\tau_p(\boldsymbol{w})
            (\underline{r}-r_0)\boldsymbol{u}
        \right\|_p
        =
        \tau_p(\boldsymbol{w})\underline{r}.
    \end{equation}
    Moreover,
    $0\le\tau_p(\boldsymbol{w})
    \le\lVert\boldsymbol{w}\rVert_p/r_0$,
    with $\tau_p(\boldsymbol{w})>0$ whenever
    $\boldsymbol{w}\neq\boldsymbol{0}$.
    In particular, when $p=2$,
    \begin{equation}
    \label{eq:ball-tau2-closed-form}
        \tau_2(\boldsymbol{w})
        =
        \frac{
            \sqrt{
                (\underline{r}-r_0)^2
                (\boldsymbol{u}^{\top}\boldsymbol{w})^2
                +
                r_0(2\underline{r}-r_0)
                \lVert\boldsymbol{w}\rVert_2^2
            }
            -
            (\underline{r}-r_0)
            \boldsymbol{u}^{\top}\boldsymbol{w}
        }{
            r_0(2\underline{r}-r_0)
        }.
    \end{equation}
\end{proposition}

\proof{Proof.}
Since
$\boldsymbol{c}_{\boldsymbol{z}}^{\boldsymbol{u}}(\underline{r})
=\boldsymbol{z}+\underline{r}\boldsymbol{u}$
by \eqref{eq:radial-center} and
$\bar{\boldsymbol{x}}=\boldsymbol{z}+r_0\boldsymbol{u}$,
we have
$\widehat{\mathcal{C}}-\bar{\boldsymbol{x}}
=\mathcal{B}_p\bigl(
    (\underline{r}-r_0)\boldsymbol{u},\underline{r}
\bigr)$.
Fix $\boldsymbol{w}\in\mathbb{R}^n$, and, for $t\ge0$, define
\[
    F_p(\boldsymbol{w},t)
    \coloneqq
    \lVert
        \boldsymbol{w}
        -t(\underline{r}-r_0)\boldsymbol{u}
    \rVert_p-t\underline{r}.
\]
For every $t>0$, the definition of $\widehat{\mathcal{C}}-\bar{\boldsymbol{x}}$ implies that $t(\widehat{\mathcal{C}}-\bar{\boldsymbol{x}}) = \mathcal{B}_p\bigl(
    t(\underline{r}-r_0)\boldsymbol{u},t\underline{r}
\bigr)$.
Therefore, $\boldsymbol{w}\in
t(\widehat{\mathcal{C}}-\bar{\boldsymbol{x}})$
if and only if
$\lVert
    \boldsymbol{w}
    -t(\underline{r}-r_0)\boldsymbol{u}
\rVert_p
\le t\underline{r}$,
which is equivalent to
$F_p(\boldsymbol{w},t)\le 0$.

The function $t\mapsto F_p(\boldsymbol{w},t)$ is continuous.
For $t'>t\ge0$, we have
\begin{align*}
    F_p(\boldsymbol{w},t')-F_p(\boldsymbol{w},t)
    &=
    \left\|
        \boldsymbol{w}
        -t'(\underline{r}-r_0)\boldsymbol{u}
    \right\|_p
    -
    \left\|
        \boldsymbol{w}
        -t(\underline{r}-r_0)\boldsymbol{u}
    \right\|_p
    -(t'-t)\underline{r} \\
    &\le
    \left\|
        (t'-t)(\underline{r}-r_0)\boldsymbol{u}
    \right\|_p
    -(t'-t)\underline{r} \\
    &=
    (t'-t)(\underline{r}-r_0)
    -(t'-t)\underline{r} \\
    &=
    -r_0(t'-t)
    <0.
\end{align*}
In the above equations, the inequality follows from the triangle inequality, while the second equality uses $\underline{r}\ge r_0$ and
$\lVert\boldsymbol{u}\rVert_p=1$.
The strict inequality follows from $r_0>0$ and $t'>t$.
Thus, $F_p(\boldsymbol{w},\cdot)$ is strictly decreasing on
$[0,\infty)$.

If $\boldsymbol{w}\neq\boldsymbol{0}$, then
$F_p(\boldsymbol{w},0)=\lVert\boldsymbol{w}\rVert_p>0$.
Moreover, for every $t\ge0$,
\begin{align*}
    F_p(\boldsymbol{w},t)
    &=
    \left\|
        \boldsymbol{w}
        -t(\underline{r}-r_0)\boldsymbol{u}
    \right\|_p
    -t\underline{r} \\
    &\le
    \lVert\boldsymbol{w}\rVert_p
    +
    \left\|
        t(\underline{r}-r_0)\boldsymbol{u}
    \right\|_p
    -t\underline{r} \\
    &=
    \lVert\boldsymbol{w}\rVert_p-r_0t,
\end{align*}
where the inequality follows from the triangle inequality, and the last equality uses $\underline{r}\ge r_0$ and
$\lVert\boldsymbol{u}\rVert_p=1$.
Consequently,
$F_p(\boldsymbol{w},
\lVert\boldsymbol{w}\rVert_p/r_0)\le0$.
By continuity, $F_p(\boldsymbol{w},\cdot)$ has a root in
$(0,\lVert\boldsymbol{w}\rVert_p/r_0]$,
and strict monotonicity implies that this root is unique.
If $\boldsymbol{w}=\boldsymbol{0}$, then
$F_p(\boldsymbol{0},t)=-r_0t$,
whose unique nonnegative root is $t=0$.
In either case, denote the unique root by
$\tau_p(\boldsymbol{w})$.

Since $t\mapsto F_p(\boldsymbol{w},t)$ is strictly decreasing with
unique root $\tau_p(\boldsymbol{w})$, we have
$F_p(\boldsymbol{w},t)\le0$ if and only if
$t\ge\tau_p(\boldsymbol{w})$.
Together with the equivalence established above between
$F_p(\boldsymbol{w},t)\le0$ and
$\boldsymbol{w}\in
t(\widehat{\mathcal{C}}-\bar{\boldsymbol{x}})$,
this implies that
$\boldsymbol{w}\in
t(\widehat{\mathcal{C}}-\bar{\boldsymbol{x}})$
if and only if
$t\ge\tau_p(\boldsymbol{w})$
for every $t>0$.
Hence, by the definition of the gauge,
$\gamma_{\widehat{\mathcal{C}}-\bar{\boldsymbol{x}}}
(\boldsymbol{w})=\tau_p(\boldsymbol{w})$.

It remains to derive the closed-form 
expression \eqref{eq:ball-tau2-closed-form} for $p=2$.
Squaring \eqref{eq:ball-root-equation} and using
$\lVert\boldsymbol{u}\rVert_2=1$ shows that
$\tau_2(\boldsymbol{w})$ satisfies
\[
    r_0(2\underline{r}-r_0)
    \tau_2(\boldsymbol{w})^2
    +
    2(\underline{r}-r_0)
    (\boldsymbol{u}^{\top}\boldsymbol{w})
    \tau_2(\boldsymbol{w})
    -
    \lVert\boldsymbol{w}\rVert_2^2
    =
    0.
\]
Solving this quadratic equation for its nonnegative root gives
\eqref{eq:ball-tau2-closed-form}.
\Halmos
\endproof

The following corollary gives the intersection cuts obtained from the sets returned by 
Algorithm~\ref{alg:normal-ray-ball-expansion}.

\begin{corollary}
\label{cor:ball-IC}
    Let $\widehat{\mathcal{C}}$ be the set returned by
    Algorithm~\ref{alg:normal-ray-ball-expansion} at
    $\bar{\boldsymbol{x}}$ with norm order $p$ and tolerance $\varepsilon$.
    Then the intersection cut obtained from
    $\widehat{\mathcal{C}}$ on $\mathcal{P}'$ is given as follows.
    \begin{enumerate}[label=\textup{(\roman*)}, ref=\textup{(\roman*)}]
    \item
    If $r^\ast<+\infty$, the intersection cut is
    \begin{equation}
    \label{eq:ball-IC}
        \sum_{j\in R}
        \tau_p(\boldsymbol{v}^j)\lambda_j(\boldsymbol{x})
        \ge1,
    \end{equation}
    where $\tau_p$ is defined in
    Proposition~\ref{prop:expanded-ball-gauge}.

    \item
    If $r^\ast=+\infty$ and $p\in(1,\infty)$, then
    $\widehat{\mathcal{C}}=\mathcal{H}_p$, where
    $\mathcal{H}_p$ is the halfspace defined in
    \eqref{eq:normal-ray-limit-halfspace}, namely,
    $\mathcal{H}_p
    =\{\boldsymbol{x}\in\mathbb{R}^n:
    \boldsymbol{a}^{\top}(\boldsymbol{x}-\boldsymbol{z})\ge0\}$,
    with
    $a_i=\operatorname{sgn}(u_i)|u_i|^{p-1}$
    for $i\in[n]$.
    The corresponding intersection cut is
    \begin{equation}
    \label{eq:ball-halfspace-IC}
        \sum_{j\in R}
        \frac{(-\boldsymbol{a}^{\top}\boldsymbol{v}^j)_+}{r_0}
        \lambda_j(\boldsymbol{x})
        \ge1.
    \end{equation}

    \item
    If $r^\ast=+\infty$ and $p=1$, then
    $\widehat{\mathcal{C}}=\mathcal{C}_{\Phi}$, where
    $\mathcal{C}_{\Phi}$ is the set defined in
    \eqref{eq:radial-l1-limit-set}, namely,
    $\mathcal{C}_{\Phi}
    =\{\boldsymbol{x}\in\mathbb{R}^n:
    \Phi(\boldsymbol{x}-\boldsymbol{z})\ge0\}$,
    with $\Phi$ defined in
    \eqref{eq:radial-l1-function}.
    The corresponding intersection cut is
    \begin{equation}
    \label{eq:ball-l1-limit-IC}
        \sum_{j\in R}
        \frac{(-\Phi(\boldsymbol{v}^j))_+}{r_0}
        \lambda_j(\boldsymbol{x})
        \ge1.
    \end{equation}
\end{enumerate}
    In each case, the inequality cuts off $\bar{\boldsymbol{x}}$
    and is valid for
    $\mathcal{P}\cap\operatorname{Sol}(\dd{D})$
    and hence for $\mathcal{P}\cap\mathcal{G}$.
\end{corollary}

\proof{Proof.}

Since
$\bar{\boldsymbol{x}}
\notin
\operatorname{Sol}(\dd{D})$,
the initial DD-free ball
$\mathcal{C}_p(\bar{\boldsymbol{x}})$
contains $\bar{\boldsymbol{x}}$ in its interior.
By Theorem~\ref{thm:normal-ray-ball-expansion-valid},
$\widehat{\mathcal{C}}$ is DD-free and contains
$\mathcal{C}_p(\bar{\boldsymbol{x}})$.
Hence,
$\bar{\boldsymbol{x}}
\in
\operatorname{int}(\widehat{\mathcal{C}})$.
We next identify the gauge of
$\widehat{\mathcal{C}}-\bar{\boldsymbol{x}}$
in each case.

Suppose first that $r^\ast<+\infty$.
By Proposition~\ref{prop:expanded-ball-gauge},
$\gamma_{
    \widehat{\mathcal{C}}
    -\bar{\boldsymbol{x}}
}
(\boldsymbol{w})=\tau_p(\boldsymbol{w})$
for every 
$\boldsymbol{w}\in\mathbb{R}^n$.

Suppose next that $r^\ast=+\infty$ and $p\in(1,\infty)$,
so that 
$\widehat{\mathcal{C}}
=\mathcal{H}_p
=\{
    \boldsymbol{x}\in\mathbb{R}^n:
    \boldsymbol{a}^{\top}(
        \boldsymbol{x}-\boldsymbol{z}
    )\ge0
\}$.
Since
$\bar{\boldsymbol{x}}
=
\boldsymbol{z}+r_0\boldsymbol{u}$
and 
$\boldsymbol{a}^{\top}\boldsymbol{u}=1$,
we have
$\mathcal{H}_p-\bar{\boldsymbol{x}}
=\{
    \boldsymbol{y}\in\mathbb{R}^n:
    \boldsymbol{a}^{\top}\boldsymbol{y}
    \ge
    -r_0
\}$.
Thus, for every 
$\boldsymbol{w}\in\mathbb{R}^n$ and $t>0$,
$\boldsymbol{w}
\in 
t(\mathcal{H}_p-\bar{\boldsymbol{x}})$
if and only if
$-\boldsymbol{a}^{\top}\boldsymbol{w}\le tr_0$.
By the definition of the gauge, taking the infimum over such $t>0$ yields
$\gamma_{\mathcal{H}_p-\bar{\boldsymbol{x}}}(\boldsymbol{w})
=(-\boldsymbol{a}^{\top}\boldsymbol{w})_+/r_0$.

Finally, suppose that $r^\ast=+\infty$ and $p=1$,
so that 
$\widehat{\mathcal{C}}
=\mathcal{C}_{\Phi}
=\{
    \boldsymbol{x}\in\mathbb{R}^n:
    \Phi(\boldsymbol{x}-\boldsymbol{z})
    \ge0
\}$.
By the definition of $\Phi$ in 
\eqref{eq:radial-l1-function},
$\Phi(\boldsymbol{y}+r_0\boldsymbol{u})
=\Phi(\boldsymbol{y})+r_0$
for every 
$\boldsymbol{y}\in\mathbb{R}^n$.
Together with
$\bar{\boldsymbol{x}}
=\boldsymbol{z}+r_0\boldsymbol{u}$, we have
$\mathcal{C}_{\Phi}-\bar{\boldsymbol{x}}
=\{
    \boldsymbol{y}\in\mathbb{R}^n:
    \Phi(\boldsymbol{y})
    \ge-r_0
\}$.
For every $\boldsymbol{w}\in\mathbb{R}^n$ and $t>0$,
it follows from \eqref{eq:radial-l1-function} that
$\Phi(\boldsymbol{w}/t)=\Phi(\boldsymbol{w})/t$.
Thus,
$\boldsymbol{w}\in
t(\mathcal{C}_{\Phi}-\bar{\boldsymbol{x}})$
if and only if
$\Phi(\boldsymbol{w}/t)\ge-r_0$,
or equivalently,
$-\Phi(\boldsymbol{w})\le tr_0$.
It follows from the definition of the gauge that, taking the infimum over such $t>0$ gives
$\gamma_{
    \mathcal{C}_{\Phi}
    -\bar{\boldsymbol{x}}
}
(\boldsymbol{w})
=
(-\Phi(\boldsymbol{w}))_+/r_0$.

Substituting these gauge representations into
\eqref{eq:IC-gauge} yields the stated cuts.
By the validity argument in
Section~\ref{subsec:IC-derivation},
each inequality is valid for
$\mathcal{P}\cap\operatorname{Sol}(\dd{D})$
and hence for $\mathcal{P}\cap\mathcal{G}$.
Finally, $\lambda_j(\bar{\boldsymbol{x}})=0$
for every $j\in R$, so $\bar{\boldsymbol{x}}$ violates
each of the corresponding inequalities.
\Halmos
\endproof

When $r^\ast<+\infty$, 
Proposition~\ref{prop:expanded-ball-gauge} shows that,
for each $j\in R$,
the coefficient $\tau_p(\boldsymbol{v}^j)$ is the unique solution of 
\eqref{eq:ball-root-equation} in
$(0,\lVert\boldsymbol{v}^j\rVert_p/r_0]$ and can be computed to any prescribed accuracy by bisection.
As a result, after any number of bisection iterations, the updated upper bound can be used in place of $\tau_p(\boldsymbol{v}^j)$ in \eqref{eq:ball-IC} to obtain a valid intersection cut.

The inclusion
$\mathcal{C}_p(\bar{\boldsymbol{x}})\subseteq\widehat{\mathcal{C}}$
yields the following comparison with the initial
$\ell_p$-ball intersection cut on the same corner relaxation.

\begin{corollary}
\label{cor:ball-IC-strength}
    Fix $p\in[1,\infty)$, and let
    $r_0\coloneqq
    d_p(\bar{\boldsymbol{x}},\operatorname{Sol}(\dd{D}))$.
    Let $\widehat{\mathcal{C}}$ be the set returned by
    Algorithm~\ref{alg:normal-ray-ball-expansion} at
    $\bar{\boldsymbol{x}}$ with norm order $p$ and any tolerance
    $\varepsilon>0$.
    Then, for every $j\in R$,
    $\gamma_{\widehat{\mathcal{C}}-\bar{\boldsymbol{x}}}
    (\boldsymbol{v}^j)\le\lVert\boldsymbol{v}^j\rVert_p/r_0$.
    Consequently, the corresponding cut in
    \eqref{eq:ball-IC}, \eqref{eq:ball-halfspace-IC}, or
    \eqref{eq:ball-l1-limit-IC} is at least as strong on the corner relaxation $\mathcal{P}'$ as the $\ell_p$-ball intersection cut
    in \eqref{eq:lp-ball-IC}.
\end{corollary}

\proof{Proof.}

By Theorem~\ref{thm:normal-ray-ball-expansion-valid},
$\mathcal{B}_p(\bar{\boldsymbol{x}},r_0)
=\mathcal{C}_p(\bar{\boldsymbol{x}})
\subseteq\widehat{\mathcal{C}}$.
Translating both sets by $-\bar{\boldsymbol{x}}$ and using
gauge monotonicity under set inclusion gives
$\gamma_{\widehat{\mathcal{C}}-\bar{\boldsymbol{x}}}
(\boldsymbol{v}^j)
\le\gamma_{\mathcal{B}_p(\boldsymbol{0},r_0)}
(\boldsymbol{v}^j)
=\lVert\boldsymbol{v}^j\rVert_p/r_0$
for every $j\in R$.
Since $\lambda_j(\boldsymbol{x})\ge0$ for every $j\in R$
and $\boldsymbol{x}\in\mathcal{P}'$, every point of
$\mathcal{P}'$ satisfying the corresponding intersection cut
also satisfies the initial $\ell_p$-ball cut.
\Halmos
\endproof

We now return to the cut closures defined in
\eqref{eq:rank-t-cut-closure} and reinstate the constraint index.
Fix a norm order $p\in[1,\infty)$. For each $k\in K$, fix any tolerance $\varepsilon_k>0$.
For
$\boldsymbol{x}\notin\operatorname{Sol}(\dd{D}^k)$,
let $\mathcal{C}_{\mathrm{ball}}^k(\boldsymbol{x})$ denote
any fixed output set returned by
Algorithm~\ref{alg:normal-ray-ball-expansion} applied to
$\dd{D}^k$ at $\boldsymbol{x}$ with norm order $p$ and
tolerance $\varepsilon_k$.
For $\boldsymbol{x}\in\operatorname{Sol}(\dd{D}^k)$,
set $\mathcal{C}_{\mathrm{ball}}^k(\boldsymbol{x})
\coloneqq\emptyset$.

\begin{corollary}
\label{cor:ball-cut-closure-convergence}
    For every $k\in K$, the mapping
    $\mathcal{C}_{\mathrm{ball}}^k(\cdot)$ defined above 
    is a DD-free-set oracle
    for $\dd{D}^k$, and the resulting family
    $\{\mathcal{C}_{\mathrm{ball}}^k(\cdot)\}_{k\in K}$
    satisfies the local covering condition
    \eqref{eq:local-covering-condition}.
    Consequently, if
    $\{\mathcal{P}_t\}_{t\in\mathbb{W}}$
    denotes the cut closures in
    \eqref{eq:rank-t-cut-closure} induced by this oracle family,
    then
    \[
        \mathcal{P}_\infty
        =
        \operatorname{conv}\left(
            \mathcal{P}
            \cap
            \bigcap_{k\in K}
            \operatorname{Sol}(\dd{D}^k)
        \right).
    \]
\end{corollary}

\proof{Proof.}

Fix $k\in K$ and
$\boldsymbol{x}\notin\operatorname{Sol}(\dd{D}^k)$.
By Proposition~\ref{prop:lp-ball-DD-free},
$\mathcal{C}_p^k(\boldsymbol{x})$ is DD-free and contains
$\boldsymbol{x}$ in its interior.
Applying Theorem~\ref{thm:normal-ray-ball-expansion-valid}
shows that $\mathcal{C}_{\mathrm{ball}}^k(\boldsymbol{x})$
is DD-free and contains $\mathcal{C}_p^k(\boldsymbol{x})$.
Therefore,
$\boldsymbol{x}\in\operatorname{int}
(\mathcal{C}_{\mathrm{ball}}^k(\boldsymbol{x}))$.
Together with the definition on
$\operatorname{Sol}(\dd{D}^k)$, this proves that
$\mathcal{C}_{\mathrm{ball}}^k(\cdot)$ is a DD-free-set oracle
for $\dd{D}^k$.

Moreover,
$\mathcal{C}_p^k(\boldsymbol{x})
\subseteq\mathcal{C}_{\mathrm{ball}}^k(\boldsymbol{x})$
for every $k\in K$ and $\boldsymbol{x}\in\mathbb{R}^n$.
Corollary~\ref{cor:lp-ball-enlargement-cut-closure-convergence}
therefore applies with the fixed $p$, yielding the local
covering condition and the stated cut-closure equality.
\Halmos
\endproof

\section{Computational Results}
\label{sec:computational-results}

In this section, we evaluate an IDD-based outer-approximation procedure augmented with the proposed intersection cuts on 13 benchmark instances from MINLPLib \cite{minlplibWebsite}. 
These instances are selected from `open' problem classes, for which the dual bounds reported by modern solvers and the best-known primal bounds still exhibit notable gaps. Their computational difficulty makes them well suited for evaluating the effectiveness of the proposed framework relative to existing solvers.
We compare the resulting dual bounds with those obtained by Gurobi 13.0 and SCIP 10.1.0 at the root node.
All experiments were performed on a machine equipped with a 2.40 GHz Intel Core Ultra 7 processor and 64 GB of RAM, running Windows 11 Education (64-bit).
The implementation uses Julia v1.10.11 via JuMP v1.11.1 \cite{dunning2017jump}.

\subsection{Implementation Settings}
\label{subsec:implementation-settings}

Each MINLPLib instance is reformulated to fit the problem form described in \eqref{eq:problem-def}, with each equality constraint replaced by two inequalities.
For variables without explicit bounds, we derive valid finite bounds from the model constraints.

\paragraph{IDD construction.}
We construct an IDD relaxation for each nonlinear constraint using Algorithm~\ref{alg:IDD-construction}.
We use uniform subdomain partitions with partition parameter $m$.
For each continuous variable $x_i$, $i\in J$, we divide $\mathcal{D}_i$ into $m$ equal-length subintervals.
For each integer variable $x_i$, $i\in I$, we divide $\mathcal{D}_i$ into $m_i=\min\{m,|\mathcal{D}_i|\}$ subdomains of consecutive integers whose cardinalities differ by at most one.
The default settings are $m=50$ and an IDD width limit of $\omega=10000$.
For the four instances with the \textbf{eg\_} prefix, whose nonlinear constraints involve products of exponential terms, we use $m=5$ and $\omega=50$.
We use the merging policy in 
\cite[Algorithm~5]{davarnia2026graphical}.
Whenever a node layer exceeds this limit, the range of its state values is divided into $\omega$ sub-ranges of equal length, and nodes assigned to the same sub-range are merged.
The look-ahead states are computed using the monotonicity-based bounding rules described in 
\cite[Section~3.4]{davarnia2026graphical}.

\paragraph{Cutting-plane procedure.}
We incorporate the proposed intersection cuts into an LP-based outer-approximation procedure, denoted by \emph{DD-IC}.
The procedure also uses the IDD cuts described in Section~\ref{subsec:outer-approximation}.
The initial LP relaxation consists of the variable bounds and all linear constraints of the reformulated instance.
At each iteration, we solve the current LP relaxation to obtain an optimal basic feasible solution, add cuts generated from the individual IDDs, and reoptimize the LP.
To generate intersection cuts, we first compute a DD-free $\ell_\infty$-ball using Algorithm~\ref{alg:lp-ball-DD-free} with $p=\infty$ and then apply Algorithm~\ref{alg:box-expansion} to expand this initial box, processing variables in increasing index order.
The intersection cuts are generated from the resulting box-maximal DD-free sets using \eqref{eq:box-IC}.
To generate IDD cuts, we use the subgradient-type separation method in 
\cite[Algorithm~6]{davarnia2026graphical}, with a constant step size of 1, the origin as the starting point, and a limit of 20 iterations per separation call.
The outer-approximation procedure terminates when no additional cuts are generated or when the LP relaxation bound improves by no more than $0.01\%$ over the last 20 LP solves.

\paragraph{Solver configurations.}
To strengthen the dual bounds at the root node, we set Gurobi's \texttt{Cuts} parameter to 3 and use SCIP's aggressive separation setting, enabling both solvers to employ their full arsenal of cutting-plane tools and techniques to improve the root-node dual bound.
To ensure a fair comparison, for each instance, we set the time limit for both solvers equal to the total running time of DD-IC.
Both solvers are set to process the root node of the B\&B tree.
All other solver parameters are left at their default values.

\subsection{Results}
\label{subsec:computational-results}

Table~\ref{tab:computational-results} summarizes the nonlinear structures and computational results for the 13 MINLPLib instances.
The column ``Primal'' reports the best-known primal bounds from MINLPLib, with signs adjusted to match the maximization formulation in \eqref{eq:problem-def}.
The columns ``DD-IC'', ``Gurobi'', and ``SCIP'' report the corresponding dual bounds, while the adjacent columns $\Delta_1$ and $\Delta_2$ report the relative gap improvements achieved by DD-IC over Gurobi and SCIP, respectively.
For instance, the relative gap improvement is computed as
$\frac{\mathrm{Gurobi}-\zeta^{\mathrm D}}
{\mathrm{Gurobi}-\zeta^{\mathrm P}}$ for Gurobi, 
where $\zeta^{\mathrm P}$ denotes the primal bound from MINLPLib and $\zeta^{\mathrm D}$ is the dual bound obtained by DD-IC.
The column ``Time'' reports the total running time of DD-IC in seconds.

\begin{table}[tbp]
\TABLE
{Computational results for MINLPLib instances.
\label{tab:computational-results}}
{%
\small
\setlength{\tabcolsep}{2.5pt}
\renewcommand{\arraystretch}{1.12}
\begin{tabular*}{\textwidth}{@{\extracolsep{\fill}}llrrrrrrrr@{}}
\toprule
Instance
& \shortstack[l]{Nonlinear\\structures}
& Primal
& DD-IC
& Time
& Gurobi
& $\Delta_1$
& SCIP
& $\Delta_2$\\
\midrule
ball\_mk4\_15 & poly, sqr
& Infeasible & Infeasible & $15.74$
& $77$ & -- & $33.85$ & -- \\
cesam2log & \shortstack[l]{exp, log,\\ poly}
& $-0.51$ & $14.35$ & $37.27$
& $639.1$ & $97.68\%$ & $468.17$ & $96.83\%$ \\
eg\_all\_s & exp, sqr
& $-7.66$ & $-4.96$ & $41.28$
& $2.67$ & $73.86\%$ & $3.91$ & $76.66\%$ \\
eg\_disc\_s & exp, sqr
& $-5.76$ & $-2.44$ & $43.01$
& $9.71$ & $78.54\%$ & $6.86$ & $73.69\%$ \\
eg\_disc2\_s & exp, sqr
& $-5.64$ & $0.96$ & $56.46$
& $9.77$ & $57.17\%$ & $6.64$ & $46.25\%$ \\
eg\_int\_s & exp, sqr
& $-6.45$ & $-0.57$ & $34.14$
& $3.72$ & $42.18\%$ & $4.26$ & $45.10\%$ \\
ex6\_2\_5 & frac, log
& $70.75$ & $2359.37$ & $36.08$
& $16528.73$ & $86.09\%$ & $15786.79$ & $85.44\%$ \\
ex6\_2\_7 & log
& $0.16$ & $25.98$ & $38.42$
& $130.99$ & $80.26\%$ & $172.84$ & $85.05\%$ \\
maxmin & sqr, sqrt
& $0.37$ & $1.03$ & $26.26$
& $1.15$ & $15.38\%$ & $1.08$ & $7.04\%$ \\
pricing050 & \shortstack[l]{exp, power,\\ sqr}
& $-1813.83$ & $-1706.47$ & $25.43$
& $-1117.24$ & $84.59\%$ & $-1382.54$ & $75.11\%$ \\
fct & \shortstack[l]{abs, mod,\\ sin, sqr}
& $0$ & $0$ & $39.15$
& -- & -- & -- & -- \\
quantum & \shortstack[l]{frac, gamma,\\ vpower, sqr}
& $-0.8$ & $-0.74$ & $25.68$
& -- & -- & -- & -- \\
worst & \shortstack[l]{frac, errorf,\\ exp, log, sqr}
& $-20762609.21$ & $-19891838.51$ & $39.51$
& -- & -- & -- & -- \\
\bottomrule
\end{tabular*}
}
{%
\textit{Notes.}
The labels ``sqr'' and ``power'' denote squared terms and integer powers, respectively, whereas ``poly'' denotes other polynomial terms.
The label ``vpower'' denotes powers whose exponents depend on decision variables.
The label ``frac'' denotes fractional expressions whose denominators depend on decision variables.
The labels ``mod'', ``gamma'', and ``errorf'' denote the remainder operation, gamma function, and standard normal cumulative distribution function, respectively.
The entry ``Infeasible'' indicates proven infeasibility.
The symbol ``--'' indicates an unavailable dual bound or an undefined gap.
}
\end{table}

For all nine instances for which all three methods report finite dual bounds, DD-IC yields tighter bounds than both Gurobi and SCIP.
The relative gap improvements range from $15.38\%$ to $97.68\%$, with an average of $68.42\%$ over Gurobi, and from $7.04\%$ to $96.83\%$, with an average of $65.68\%$ over SCIP, demonstrating substantial improvements in gap closure.
For example, on \textbf{pricing050}, DD-IC achieves relative gap improvements of $84.59\%$ and $75.11\%$ over Gurobi and SCIP, respectively.
Moreover, DD-IC proves the infeasibility of \textbf{ball\_mk4\_15}, which neither solver detects within the time limit at the root node.
The procedure also obtains dual bounds for \textbf{fct}, \textbf{quantum}, and \textbf{worst}, for which no dual bounds are reported by either solver because these instances contain unsupported nonlinear operators.
In particular, the dual bound obtained for \textbf{fct} matches the primal bound.
These results highlight the potential of the proposed DD-IC outer-approximation procedure to strengthen MINLP relaxations across a wide range of challenging nonlinear structures.

\section{Conclusion}
\label{sec:conclusion}

In this paper, we develop a DD-based framework for constructing $S$-free sets and deriving intersection cuts for bounded MINLPs with constraints of general nonlinear functional form.
We first construct an IDD relaxation for each nonlinear constraint and exploit the resulting graph structure to efficiently compute DD-free $\ell_p$-balls.
We show that, under a local covering condition, the sequence of cut closures generated using DD-free sets for the individual IDDs converges to the convex hull of the intersection of the initial polyhedral relaxation with all IDD solution sets.
We use DD-free sets together with corner relaxations to derive intersection cuts that strengthen outer approximations beyond the convexification of individual IDD solution sets. We then study maximality of DD-free sets and characterize the computational complexity of verifying maximality in general. To obtain more tractable constructions and verification procedures, we introduce the notion of maximality within structured families. In particular, we first consider the family of DD-free boxes and develop an efficient sequential facet-expansion algorithm for constructing box-maximal DD-free sets. We then consider the family of DD-free $\ell_p$-balls for $p\in[1,\infty)$, establish maximality conditions, and develop a radial expansion procedure for enlarging the initial balls. Computational experiments on MINLPLib instances demonstrate that combining the proposed intersection cuts with existing IDD cuts can strengthen dual bounds beyond those obtained by state-of-the-art global solvers.




\bibliographystyle{informs2014} 
\bibliography{references} 





  



\end{document}